\documentclass[3p,times]{elsarticle}

\usepackage{ecrc}
\usepackage{amscd}

\volume{00}

\firstpage{1}

\journalname{}

\runauth{}

\jid{procs}

\jnltitlelogo{}

\usepackage{amssymb}

\usepackage[figuresright]{rotating}
\usepackage{mhchem}
\usepackage{esint}
\usepackage{mathtools}
\usepackage{bm}
\usepackage{amsmath}
\allowdisplaybreaks
\numberwithin{equation}{section}
\numberwithin{figure}{section}
\numberwithin{table}{section}
 
\usepackage{graphicx }
\usepackage{subfigure,txfonts}
\usepackage{amsthm}
\usepackage{amssymb}
\usepackage{cancel}
\usepackage{mathdots}
\usepackage{stmaryrd}
\usepackage{stackrel}
\usepackage{graphicx}
\usepackage{booktabs}
\usepackage{lineno}
\usepackage{nameref}
\usepackage{xcolor}
\usepackage[colorlinks,
linkcolor=blue,
anchorcolor=yellow,
citecolor=red,
]{hyperref}

 \allowdisplaybreaks[4]
\usepackage{version}
\usepackage{cases}
\usepackage{color}

\def\H{\boldsymbol{H}}
	\def\L{\boldsymbol{L}}
	
    \def\X{\boldsymbol{X}}
    \def\Y{\boldsymbol{Y}}
    
    \def\I{\boldsymbol{I}}
    \def\P{\boldsymbol{P}}

    \def\B{\boldsymbol{B}}
    \def\C{\boldsymbol{C}}

    \def\u{\boldsymbol{u}}
    \def\v{\boldsymbol{v}}
    \def\w{\boldsymbol{w}}

    \def\0{\boldsymbol{0}}
    \def\x{\boldsymbol{x}}
    
    \def\n{\boldsymbol{n}}
    
\usepackage{lineno}
\makeatletter

\theoremstyle{plain}
\ifx\thesection\undefined
\newtheorem{thm}{\protect\theoremname}
\else

\fi
\theoremstyle{remark}
\ifx\thesection\undefined
\newtheorem{rem}{\protect\remarkname}

\else

\fi
\theoremstyle{plain}
\ifx\thesection\undefined
\newtheorem{lem}{\protect\lemmaname}
\else

\fi
\theoremstyle{plain}
\ifx\thesection\undefined
\newtheorem{Example}{\protect\lemmaname}
\else

\fi

\graphicspath{{./figures/}}

\begin{document}
\newtheorem{The}{Theorem}[section]
\newtheorem{lemma}{Lemma}[section]
\newtheorem{Remark}{Remark}[section]
\newtheorem{Assumption}{Assumption}[section]
\newtheorem{Definition}{Definition}[section]
\newtheorem{Proposition}{Proposition}[section]
 
\newcommand{\sign}{\mathrm{sign}}

\begin{frontmatter}



\dochead{}
\title{Optimal error estimates for a fully
discrete, highly efficient decoupled scheme for the 2D/3D  diffuse interface  two-phase MHD flows\tnoteref{t1}}

\tnotetext[t1]{This work is partly supported by the NSF of China (No. 12471392) support.
}

\author[zzu]{Ke Zhang}
\ead{zhangkemath@139.com, zkmath@stu.xju.edu.cn}

\author[zzu]{Haiyan Su\corref{cor1}}%
\ead{shymath@126.com }
\cortext[cor1]{Corresponding author. Tel./fax number: 86 9918582482.}

\address[zzu]{College of Mathematics and System Sciences,
Xinjiang University, Urumqi 830046, P.R. China}
  
\begin{abstract}
In this paper, we derive optimal $\L^{2}$- and $\H^{1}$-norm error estimates for a fully discrete convex-splitting decoupled finite element method (FEM) for the two-phase diffuse interface magnetohydrodynamics (MHD) system. We use the semi-implicit backward Euler scheme in time and employ the standard inf–sup stable Taylor--Hood/Mini elements to discretize the velocity and pressure. Furthermore, we apply a pressure-correction scheme to decouple the velocity from the pressure. The optimal error estimates are obtained via novel Ritz and Stokes quasi-projection techniques. In addition, the unconditional energy stability of the proposed scheme is ensured. Numerical examples are presented to validate the theoretical analysis.

\end{abstract}


\end{frontmatter}
 
\section{Introduction}

The two-phase MHD focuses on the dynamic behavior of two incompressible and immiscible conducting fluids under an external electromagnetic field. The governing model consists of the Cahn--Hilliard equations (describing the free interface), the Navier--Stokes equations (describing the hydrodynamics), and the Maxwell equations (describing the magnetic field), which are coupled  through convection, stresses, and Lorentz forces.  It has extensive application prospects in the fields of nuclear fusion, metallurgy, liquid metal magnetic pumps, aluminum electrolysis, as well as in addressing problems encountered in other fields  \cite{2010An, Jean2006Mathematical, 2000Liquid}.

In this paper, we mainly consider the following two-phase diffuse interface MHD model \cite{2019A, 2023Error}:
\begin{equation}\label{model}
\begin{aligned}
&\phi_{t}+\nabla \phi \cdot\u=\gamma \nabla\cdot (M \nabla \omega),  \\
&\omega=- \gamma\Delta\phi+\gamma^{-1}f(\phi),  \\
&\u_{t}- \nabla\cdot (\nu  \nabla\u)+(\u\cdot \nabla) \u+ \frac{1}{\mu}\B\times\nabla\times\B+\nabla p= \lambda \omega \nabla \phi + \bm{f}, \\
&\nabla \cdot \u=0,  \\
&\B_{t}+ \frac{1}{\mu }\nabla\times(\frac{1}{\sigma }\nabla\times\B)- \nabla\times(\u\times\B) = \0,  \\
&\nabla \cdot\B=0,  \\
& \frac{\partial\phi}{\partial \n}|_{\partial\Omega}=0, \, \frac{\partial w}{\partial \n}|_{\partial\Omega}=0, \, \u|_{\partial\Omega}=\0,\,\B\times \n|_{\partial\Omega}=\0,\\
&\phi|_{t=0}=\phi_0,\, \u|_{t=0}=\u_0,\, \B|_{t=0}=\B_0,
\end{aligned}
\end{equation}
for $(\boldsymbol{x}, t )\in  \Omega\times (0,T]$,  $\Omega$ is a bounded smooth polyhedral domain in $R^{d}$, d=2, 3, and T indicates the final time. The phase field $\phi$  expresses the mixture of two immiscible, incompressible fluids. The two different conducting fluids can be labeled by
\begin{equation}\label{2-1c}
\phi (x, t)=\left\{
\begin{aligned}
-1, \qquad \mathrm{fluid\ 1},\\
1, \qquad \mathrm{fluid\ 2}.
\end{aligned}
\right.
\end{equation}

The function $f(\phi)$ is the derivative of the Ginzburg--Landau double-well potential function $F(\phi)=\frac{1}{4}(\phi^{2}-1)^{2}$ with respect to $\phi$ \cite{Novick2008}. The unknown variables ($\boldsymbol{u}, p$) represent the velocity field and pressure field, $\omega$ denotes the chemical potential, and $\boldsymbol{B}$ represents the magnetic field. Several positive parameters are introduced, such as the interfacial width $\gamma$ between the two phases,  the mobility parameter $M$, the kinematic viscosity $\nu$ (inverse of the Reynolds number), the magnetic permeability $\mu$, the electric conductivity $\sigma$, and the capillary coefficient $\lambda$.
 
Recent studies have focused on developing efficient numerical schemes to handle the challenges posed by the strong nonlinearity and coupled effects. The diffuse interface two-phase MHD model based on Cahn--Hilliard dynamics was first proposed and analyzed in \cite{2019A}. Later work aimed to  create energy-stable numerical schemes for this system.  As for the first-order schemes, fully decoupled invariant energy quadratization (IEQ) scheme was proposed in \cite{2022Highly}, the semi-implicit stabilization scheme was presented in 
 \cite{ZHANG2023107477, chen2022unconditional}, and the convex splitting schemes were shown in \cite{2019A, 2023Error}.
Moreover, the second-order schemes \cite{2023Energy, wang2024fully}  were also presented to handle the considered model.  For convenience, we employ the convex splitting scheme in this paper, which was proposed in \cite{1998UnconditionallyE} and has been popularized in \cite{2014ExistenceHAN, 2010UnconditionallyWISE}. 


Additionally, theoretical investigations into two-phase MHD system remain an active area. Recently, the first-order semi-implicit stabilization scheme was developed in \cite{chen2022unconditional}. The convergence analysis for the MINI finite element pair used for the velocity field and pressure field, and the $P_{1}$ element for other variables, is carried out as follows:
\begin{equation}\label{zhang}
\| \phi^{k+1}-\phi_{h}^{k+1} \|+\|\nabla(\phi^{k+1}-\phi_{h}^{k+1} )\|+ \|\nabla(\u^{k+1}-\u_{h}^{k+1})\|+  \|\nabla(\B^{k+1}-\B_{h}^{k+1})\|\leq C_{0}(\Delta t+ h).
\end{equation}

The first-order Euler semi-implicit discretization based on a convex-splitting scheme was given in \cite{2023Error}. This discretization employs the standard inf-sup stable Taylor--Hood finite element  pair $(\boldsymbol{u}_{h}^{k+1}, p_{h}^{k+1})\in \boldsymbol{X}_{h}^{r+1}\times \mathring{S}_{h}^{r}$, with $\boldsymbol{B}_{h}^{k+1}\in \boldsymbol{Y}_{h}^{r+1}$, and $\phi_{h}^{k+1}, \omega_{h}^{k+1}\in S_{h}^{r+1}$. The error estimates in \cite{2023Error} are  
\begin{equation}\label{qiu}
\|\nabla(\phi^{k+1}-\phi_{h}^{k+1} )\|+ \|\u^{k+1}-\u_{h}^{k+1}\|+  \|\B^{k+1}-\B_{h}^{k+1}\|\leq C_{0}(\Delta t+ h^{r+1}).
\end{equation}
The second-order modified Crank--Nicolson type fully discrete scheme  was designed for solving the two-phase MHD model \cite{wang2024convergence}. Specifically, it achieved the same spatial convergence order as that presented in  (\ref{qiu}). 

For numerical velocity, the $\H^1$-norm error estimate was presented only in \eqref{zhang}. The $\L^2$-norm error estimate in \eqref{qiu} is one order lower than the optimal one. 
On the orther hand, the  $H^{1}$-norm error estimate of the phase field $\phi$ is optimal in equations \eqref{zhang}-(\ref{qiu}), whereas the $L^{2}$-norm is not.  The primary reason for this is the artificial pollution resulting from the approximation of the phase field $\phi$, which affects the accuracy of numerical velocity field and magnetic field  analysis. The standard projection operators in the traditional sense may not be valid due to this pollution. The novel elliptic Ritz and Stokes quasi-projections were provided in  \cite{2023Optimalwang} to avoid the above artificial pollution. In addition,  we should note that the  $\L^{2}$-norm error estimates for the velocity field and magnetic field (\ref{qiu}) are one order lower in accuracy than the optimal estimates. Recently, error estimates based on the above quasi-projections were obtained in \cite{DUAN202643}. However, in the existing literature, few studies provide optimal error analysis for decoupled schemes.  

Thus, our goal is to  derive the optimal $L^{2}/\L^2$-norm error estimates of the pressure-correction scheme for the phase field, velocity field, and magnetic field by employing the elliptic Ritz, Stokes quasi-projections and Maxwell projection. For the uniform MINI type finite element  $(\phi_{h}^{k+1}$, $\omega_{h}^{k+1}$, $\boldsymbol{u}_{h}^{k+1}$, $p_{h}^{k+1}, \boldsymbol{B}_{h}^{k+1})\in S_{h}^{1}\times S_{h}^{1}\times \boldsymbol{X}_{h}^{1b}\times \mathring{S}_{h}^{1}\times \boldsymbol{Y}_{h}^{1}$, we obtain the optimal $L^{2}/\L^2$ error estimates 
\begin{equation}\label{L1}
\|\phi^{k+1}-\phi_{h}^{k+1} \|+\|\u^{k+1}-\u_{h}^{k+1}\| +\|\B^{k+1}-\B_{h}^{k+1}\|\leq C_{0}(\Delta t+h^{2}).
\end{equation}
For the inf-sup stable Taylor--Hood type finite element pair $(\phi_{h}^{k+1}, \omega_{h}^{k+1}, \u_{h}^{k+1}, p_{h}^{k+1}, \B_{h}^{k+1})\in S_{h}^{r}\times S_{h}^{r}\times \X_{h}^{r+1}\times \mathring{S}_{h}^{r}\times \Y_{h}^{r+1}$, we obtain the following main error estimates 
\begin{equation}\label{Taylor-Hood}
\left\{
\begin{aligned} 
\|\phi^{k+1}-\phi_{h}^{k+1} \|&\leq C_{0}(\Delta t+h^{r+1}),\\
\|\u^{k+1}-\u_{h}^{k+1}\|  &\leq C_{0}(\Delta t+\beta_{h}),\\
\|\B^{k+1}-\B_{h}^{k+1}\|&\leq C_{0}(\Delta t+\beta_{h}),
\end{aligned}
\right.
\quad 
\text{where}
\quad  
\beta_{h}= 
\left\{
\begin{aligned}
h^{r+2},\quad  &r\geq2,\\
h^{r+1},\quad  &r=1.
\end{aligned}
\right.
\end{equation}

The error estimates for the velocity field and magnetic field are one order lower than the  interpolation error when $r$=1, while our numerical results confirm that the estimates are also optimal. For $r\geq$2, the results in conclusion (\ref{Taylor-Hood}) indicate that the error estimates for the phase field, velocity field, and magnetic field are optimal. The more details can be seen in Theorem \ref{theorem2-1}. 

The rest of this work is organized as follows: In Section \ref{sec-main}, we state  the main results. In Section \ref{sec-projection}, we present the Ritz and Stokes quasi-projections.  Furthermore, Theorem \ref{theorem2-1} is derived in Section \ref{sec-error}.   Numerical examples are conducted in Section \ref{sec-examples} to confirm our theoretical analysis and demonstrate the efficiency of the method. The concluding remarks are summarized in Section \ref{sec-conclusion}.
 
\section{Main results}\label{sec-main}

In this section, we consider a fully discrete finite element convex-splitting algorithm for the two-phase MHD model (\ref{model}) and discuss the main error estimates.
 
\subsection{Preliminaries and weak formulation}
Let $W^{k,p}(\Omega)$ denote the standard Sobolev spaces, equipped with the standard Sobolev norms $\|\cdot\|_{W^{k,p}}$, for $k\geq0$, $1\leq p\leq \infty$. As usual, we write $H^{k}(\Omega)$=$W^{k,2}(\Omega)$ and $L^{p}(\Omega)$=$W^{0,p}(\Omega)$. Furthermore, we denote the norms of $H^{k}(\Omega)$ and $L^{p}(\Omega)$ by $\|\cdot\|_{H^{k}}$ and $\|\cdot\|_{L^{p}}$, respectively. Specifically,  the inner product and norm in $L^{2}(\Omega)^{d}$ are denoted by $(\cdot, \cdot)$ and $\|\cdot\|$.  The standard Sobolev spaces are as follows:
\begin{equation*}
\begin{aligned}
&H^{1}(\Omega)=\{\psi \in L^{2}(\Omega): \nabla\psi \in \L^{2}(\Omega) \},\quad \H_{0}^{1}(\Omega)=\{\v\in  H^{1}(\Omega)^{d}: \boldsymbol{v}|_{\partial\Omega}=\0 \},\\
& L_{0}^{2}(\Omega)=\{p\in L^{2}(\Omega): \int_{\Omega}p\,\mathrm{d}\x=0\},\quad \H_{\tau}^{1}(\Omega)=\{\B\in H^{1}(\Omega)^{d}: \n\times \B|_{\partial\Omega}=\0 \}.
\end{aligned}
\end{equation*}


We assert that the weak solution of the two-phase MHD model satisfies the following regularity assumption and variational
formulation. For brevity, we set the parameters $\gamma=M =\nu =\mu=\lambda=\sigma $=1.

\begin{Definition}\label{def21}
We suppose that the solution to the considered model (\ref{model}) exists and satisfies:

1. The regularity assumption ($r\geq1$),
\begin{equation}\label{regularity}
\begin{aligned}
& \phi \in H^{2}(0,T; L^{2}(\Omega))\cap H^{1}(0,T; H^{r+1}(\Omega))\cap C(0,T; W^{2, 4}(\Omega)), \quad  \omega\in H^{1}(0,T; H^{r+1}(\Omega)),\\
&\u\in H^{2}(0,T; \L^{2}(\Omega))\cap H^{1}(0,T; \H^{r+2}(\Omega)),\quad p\in L^{2}(0,T; H^{r+1}(\Omega)\cap L_{0}^{2}(\Omega)),\\
&\B\in H^{2}(0,T; \L^{2}(\Omega))\cap H^{1}(0,T; \H^{r+2}(\Omega)).
\end{aligned}
\end{equation}

2. We find $(\phi, \omega, \u, p, \B)\in H^{1}(\Omega)\times H^{1}(\Omega)\times \H_{0}^{1}(\Omega)\times L_{0}^{2}(\Omega)\times \H_{\tau}^{1}(\Omega)$
from 
\begin{subequations}\label{weak-form}
\begin{align}
&(\phi_{t}, \xi)+(\nabla \phi \cdot\u, \xi) +(\nabla \omega, \nabla \xi) =  0,   \label{mdel1}\\
&(\omega, \psi) = (\nabla\phi, \nabla \psi) +(f(\phi), \psi),    \label{model2}\\
&(\u_{t}, \v)+(\nabla\u, \nabla\v)+b(\u, \u,\v) + (\nabla\times\B, \v\times\B)-( p, \nabla\cdot\v)= (\omega \nabla \phi, \v) + (\bm{f}, \v),  \label{model3}\\
&(\nabla \cdot \u, q)=0, \label{model4}\\
&(\B_{t}, \C)+ (\nabla\times\B, \nabla\times \C)- (\u\times\B, \nabla\times\C)+(\nabla \cdot\B, \nabla \cdot\C)  =0, \label{model5}
\end{align}
\end{subequations}
for any test function $(\xi, \psi, \v, q, \C)\in H^{1}(\Omega)\times H^{1}(\Omega)\times \H_{0}^{1}(\Omega)\times L_{0}^{2}(\Omega)\times \H_{\tau}^{1}(\Omega)$. And here we define  $b(\u, \v, \w)=\frac{1}{2}[( \u\cdot \nabla \v, \w)- ( \u\cdot \nabla  \w, \v) ]$, see in \cite{1975Sobolev}.
\end{Definition}
 
\begin{Remark} Let ($\phi, \omega, \u, p, \B$) be the weak solution of the  two-phase MHD model (\ref{model}). Then, for any $t\in(0,T]$, the mass is conserved,
\begin{equation}
\int_{\Omega}\phi \mathrm{d}\x=\int_{\Omega}\phi_{0} \mathrm{d}\x,
\end{equation}
and the system energy is stable,
\begin{equation}\label{model-energy}
\frac{dE(\phi, \u, \B)}{dt} =-(\|\nabla\omega\|^{2}+ \|\nabla\u \|^{2}+ \| \nabla\times \B\|^{2} )\leq 0,
\end{equation}
where the total energy  is given by
\begin{equation}\label{system-energy}
E(\phi, \u, \B)= \frac{1}{2}\|\nabla\phi\|^{2}+\int_{\Omega}F(\phi) dx+\frac{1}{2}\|\u\|^{2}+\frac{1}{2 }\|\B\|^{2},
\end{equation}
  and set $\bm{f}=\0$ without loss of generality.
\end{Remark}

\subsection{Numerical scheme and main results}
 Let $\zeta_{h}$ be a regular and quasi-uniform partition of $\Omega$ with
mesh size $h=\max\limits_{1\leq j\leq M}{ \rm diam} (K_{j})$.  To approximate the unknowns variables, we employ the finite element spaces as follows,
\begin{equation*}
\begin{aligned}
& S_{h}^{r}=\{\phi_{h}\in C(\Omega): \phi_{h}|_{K_{j}}\in P_{r}(K_{j}), \forall K_{j} \in \zeta_{h}\},\\
&\mathring{S}_{h}^{r}=S_{h}^{r}\cap L_{0}^{2}(\Omega),\\
& \X_{h}^{r+1}=\{\v_{h}\in H_{0}^{1}(\Omega)^{d}: \v_{h}|_{K_{j}}\in P_{r+1}(K_{j})^{d}, \forall K_{j} \in \zeta_{h}\},\\
& \X_{h}^{1b}= (S_{h}^{1}\oplus B_{d+1})^{d}\cap \H_{0}^{1}(\Omega),\\
& \Y_{h}^{r+1}=\{\B_{h}\in \H_{\tau}^{1}(\Omega): \B_{h}|_{K_{j}}\in P_{r+1}(K_{j})^{d}, \forall K_{j} \in \zeta_{h}\},
\end{aligned}
\end{equation*}
where $r\geq1$, $P_{r}(K_{j})$ is the space of polynomials of total degree $r$ on $K_{j}$, $B_{3}$ and $B_{4}$ are the spaces of cubic bubbles and quartic bubbles for 2D and 3D \cite{2023Optimalwang}, respectively. According to the classical finite element theory, we have the following discrete inf-sup condition \cite{1986FiniteG, 1991BF}.
\begin{Assumption}\label{a1}
The Taylor-Hood element $\X_{h}^{r+1}\times \mathring{S}_{h}^{r}$ and MINI element  $\X_{h}^{1b}\times \mathring{S}_{h}^{r}$ satisfy the following  inf-sup condition:
\begin{equation}\label{infsup}
\begin{aligned}
&\inf_{q_{h}\in (\mathring{S}_{h}^{r}/\mathring{S}_{h}^{1})\backslash\{0\} }\sup_{\v_{h}\in (\X_{h}^{r+1}/\X_{h}^{1b})\backslash\{\0\}}\frac{(\nabla\cdot \v_{h}, q_{h})}{\|q_{h}\|\|\nabla \v_{h}\|}\geq \beta_{0},
\end{aligned}
\end{equation}
where $\beta_{0}$ is a positive constants depending only on $\Omega$.
\end{Assumption} 
 
\begin{lemma}
Based on the  $Poincar\acute{e}$  inequalities and embedding inequalities in \cite{1975Sobolev, Jean2006Mathematical, 1986On, 2019A}, we have
\begin{subequations}
\begin{align}
&\| \v\|_{L^{q}}\leq C_{0}\|\nabla \v\|, \quad  \forall \v\in \boldsymbol{H}_{0}^{1}(\Omega),\, 2\leq q\leq 6, \label{i1}\\
&\|\psi\|_{L^{q}}\leq C_{0}\|\psi\|_{H^{1}},\quad \forall \psi\in H^{1}(\Omega),\, 2\leq q \leq 6,\label{i2}\\
&\|\nabla\times \C\|\leq C_{0}\|\nabla \C\|, \quad \forall \C\in \H_{\tau}^{1}(\Omega), \label{i3}\\
& \|v_{h}\|_{W^{m,s}}\leq C_{0}h^{n-m+\frac{d}{s}-\frac{d}{q}}\|v_{h}\|_{W^{n,q}},\quad v_{h}\in S_{h}^{r}, \mathring{S}_{h}^{r}, \label{i4} \X_{h}^{r+1}/\X_{h}^{1b},\quad 0\leq n\leq m\leq1,\quad  1\leq q\leq s\leq\infty,
\end{align}
\end{subequations}
where $C_{0}$ denotes a generic positive constant independent of $\Delta t$, $h$, and discretization parameters, which may take different values at different places.
\end{lemma}
 
For simplicity of notation, we denote $\phi^{k+1}=\phi (\x, t_{k+1})$ and $\u^{k+1}=\u (\x, t_{k+1})$, and introduce the following notation: 
\begin{equation*}
	D_{t}\phi^{k+1}=\frac{\phi^{k+1}-\phi^{k}}{\Delta t},
\end{equation*}
where $\Delta t$=$\frac{T}{K}$ is the time step, for arbitrary and  fixed finial time $T>0$, and $K$ is a positive integer.
And the space notations as:
\begin{equation}\label{space}
	\boldsymbol \chi_{h}^{r}:=\left\{
	\begin{aligned}
		& S_{h}^{r}\times S_{h}^{r}\times (\X_{h}^{r+1}+ \nabla \mathring{S}_{h}^{r}) \times \mathring{S}_{h}^{r}\times \Y_{h}^{r+1}, \quad &&r\geq2,\\
		&S_{h}^{1}\times S_{h}^{1}\times (\X_{h}^{2} + \nabla \mathring{S}_{h}^{1})\times \mathring{S}_{h}^{1}\times \Y_{h}^{2},  \quad &&r=1, \\ 
		&S_{h}^{1}\times S_{h}^{1}\times (\X_{h}^{1b}+ \nabla \mathring{S}_{h}^{1}) \times \mathring{S}_{h}^{1}\times \Y_{h}^{1},  \quad &&r=1.
	\end{aligned}
	\right.
\end{equation}

We find ($\phi_{h}^{k+1}, \omega_{h}^{k+1}, \u_{h}^{k+1}, p_{h}^{k+1}, \B_{h}^{k+1})$ from the following fully discrete, decoupled  finite element method for the considered two-phase MHD model (\ref{model}),

\textbf{Step 1.} Find ($\phi_{h}^{k+1}, \omega_{h}^{k+1})\in   S_{h}^{r}\times S_{h}^{r}$,  $\tilde{\u}_{h}^{k+1}\in \X_{h}^{r+1}/ \X_{h}^{1b}$,  $\B_{h}^{k+1} \in   \Y_{h}^{r+1}$ such that for all ($\xi_{h}, \psi_{h})\in S_{h}^{r}\times S_{h}^{r}$, $\v_{h}\in \X_{h}^{r+1}/ \X_{h}^{1b}$, $\C_{h}^{k+1} \in   \Y_{h}^{r+1}$ there hold 
\begin{subequations}\label{numerical-scheme1}
	\begin{align}
		&(D_{t}\phi_h^{k+1}, \xi_{h}) +(\nabla\phi_{h}^{k}\cdot \tilde{\u}_{h}^{k+1},  \xi_{h})+ ( \nabla\omega_{h}^{k+1}, \nabla\xi_{h})=0, \label{eqn123}\\
		&(\omega_{h}^{k+1}, \psi_{h})= (\nabla\phi_{h}^{k+1}, \nabla\psi_{h})+ (f(\phi_{h}^{k+1}), \psi_{h}), \label{22}\\
		&(D_{t}\B_h^{k+1}, \C_{h}) +  ( \nabla\times \B_{h}^{k+1}, \nabla\times \C_{h})-(\tilde{\u}_{h}^{k+1}\times \B_{h}^{k}, \nabla\times \C_{h})+   ( \nabla\cdot\B_{h}^{k+1}, \nabla\cdot \C_{h}) =0,\label{B-1} \\
		&(\frac{\tilde{\u}_{h}^{k+1}- \u_{h}^{k}}{\Delta t},\v_{h}) + ( \nabla\tilde{\u}_{h}^{k+1},\nabla\v_{h} )+b(\u_{h}^{k}, \tilde{\u}_{h}^{k+1},\v_{h})+(  \v_{h},  \nabla p_{h}^{k}) +(\nabla\times \B_{h}^{k+1}, \v_{h}\times \B_{h}^{k})-(\omega_{h}^{k+1}  \nabla \phi_{h}^{k}, \v_{h})=0,\label{equau} 
	\end{align}
\end{subequations}
where $f(\phi_{h}^{k+1})$=$(\phi_{h}^{k+1})^{3}-\phi_{h}^{k}$. 

\textbf{Step 2.} Find  $p_{h}^{k+1}\in \mathring{S}_{h}^{r}$ such that for all $q_{h} \in   \mathring{S}_{h}^{r}$   by 
\begin{equation}\label{numerical-scheme2}
	(\nabla p_{h}^{k+1}, \nabla q_{h})=-\frac{ 1}{\Delta t}(\nabla\cdot \tilde{\u}_{h}^{k+1},  q_{h}) + (\nabla p_{h}^{k}, \nabla q_{h}),  
\end{equation}

\textbf{Step 3.} Find $\u_{h}^{k+1} \in \X_{h}^{r+1}/ \X_{h}^{1b} + \nabla \mathring{S}_{h}^{r}$   by 
	\begin{align}\label{numerical-scheme3}
		\u_{h}^{k+1}= \tilde{\u}_{h}^{k+1}-\Delta t \nabla(p_{h}^{k+1}-p_{h}^{k}).
		\end{align}


\begin{Remark}[The weakly divergence free ]
	Taking the $L^2$ inner product  of the equation \eqref{numerical-scheme3} with $\nabla q_h \in \nabla \mathring{S}_{h}^{r}$, we obtain
	\begin{equation*}
		(u_{h}^{k+1}, \nabla q_h)= (\tilde{\u}_{h}^{k+1}, \nabla q_h)-\Delta t (\nabla(p_{h}^{k+1}-p_{h}^{k}), \nabla q_h).
	\end{equation*}
With the help of the equation \eqref{numerical-scheme2}, we have the weakly divergence free for $u_{h}^{k+1}$, namely 
\begin{equation}\label{div-free}
(u_{h}^{k+1}, \nabla q_h)=0. 
\end{equation}
\end{Remark}

%
 
\begin{The} \label{theorem2-1}
The two-phase MHD model (\ref{model}) has a unique solution $(\phi, \omega, \u, p, \B)$, which satisfies the regularity assumption (\ref{regularity}). Then the solution $(\phi_{h}^{k+1}, \omega_{h}^{k+1}, \u_{h}^{k+1}, p_{h}^{k+1}, \B_{h}^{k+1})\in \boldsymbol \chi_{h}^{r}$ of the fully discrete scheme \eqref{numerical-scheme1}-\eqref{numerical-scheme3} satisfies the following error estimates:
\begin{subequations}
\begin{align}
&\max\limits_{0\leq k\leq K-1}\|\phi^{k+1}-\phi_{h}^{k+1} \| +\left(\Delta t \sum_{k=0}^{K-1}\|\omega^{k+1}-\omega_{h}^{k+1}\|^{2}\right)^{\frac{1}{2}}\leq C_{0}(\Delta t+h^{r+1}),\\
&\max\limits_{0\leq k\leq K-1}\|\nabla(\phi^{k+1}-\phi_{h}^{k+1}) \|\leq C_{0}(\Delta t+h^{r}),\\
&\max\limits_{0\leq k\leq K-1}\left(\|\u^{k+1}-\u_{h}^{k+1}\| +\|\B^{k+1}-\B_{h}^{k+1}\| \right)\leq C_{0}(\Delta t+\beta_{h}),\\
&\left(\Delta t\sum_{k=0}^{K-1}\left(  \|\nabla \cdot (\B^{k+1}-\B_{h}^{k+1})\|^{2} +\|\nabla\times(\B^{k+1}-\B_{h}^{k+1} )\|^{2} \right)\right)^{\frac{1}{2}}\leq C_{0}(\Delta t+\beta_{h}^{\star}),
\end{align}
\end{subequations}
where
\begin{equation}\label{bta}
\beta_{h}=
\left\{
\begin{aligned}
h^{r+2},\quad  &r\geq2,\\
h^{r+1},\quad  &r=1,
\end{aligned}
\right.
\qquad
\beta_{h}^{\star}=\left\{
\begin{aligned}
h^{r+1},\quad  &r\geq2,\\
h^{r+1},\quad  &r=1 \quad ( (\u_{h}^{k+1},p_{h}^{k+1})\in \X_{h}^{2}\times \mathring{S}_{h}^{1} ),\\
h^{r},\quad    &r=1 \quad ( (\u_{h}^{k+1},p_{h}^{k+1})\in \X_{h}^{1b}\times \mathring{S}_{h}^{1} ).
\end{aligned}
\right.
\end{equation}
\end{The}

\subsection{Unconditional energy stability} 
In this subsection, we derive the unconditional energy stability of the proposed shceme. 
\begin{The} Without loss of generality, we set $\bm{f}=\0$. The scheme \eqref{numerical-scheme1}-\eqref{numerical-scheme3} satisfies  the following  discrete energy law for any $k\geq0$,
	\begin{equation}\label{energy}
		E(\phi_{h}^{k+1}, \u_{h}^{k+1}, \B_{h}^{k+1}) -E(\phi_{h}^{k}, \u_{h}^{k}, \B_{h}^{k})\leq 0,
	\end{equation}
	where the discrete  energy is given by
	\begin{equation}\label{algorithm-energy}
		E(\phi_{h}^{k+1}, \u_{h}^{k+1}, \B_{h}^{k+1})= \frac{1}{2}\|\nabla\phi_{h}^{k+1}\|^{2}+\frac{1}{4 } \| (\phi_{h}^{k+1})^2-1\|^2+\frac{1}{2}\|\u_{h}^{k+1}\|^{2}+\frac{1}{2}\|\B_{h}^{k+1}\|^{2}.
	\end{equation}
\end{The}

\begin{proof}
	Taking $\varphi_{h}= \Delta t \omega_{h}^{k+1}$ in \eqref{eqn123}, we have
	\begin{align}\label{ener:1}
		& \Delta t\big(D_{t}\phi_{h}^{k+1}, \omega_{h}^{k+1}\big) 
		+ \Delta t\big(\nabla\phi_{h}^{k}\cdot \tilde{\u}_{h}^{k+1},  \omega_{h}^{k+1}\big) 
		+  \Delta t \|\nabla\omega_{h}^{k+1}\|^2=0.
	\end{align}
	Taking $ \psi_{h}= \Delta t D_{t}\phi_{h}^{k+1}$ in  \eqref{22}, we have
	\begin{align}\label{ener:2}
		& \Delta t\big(\nabla\phi_{h}^{k+1}, \nabla  D_{t}\phi_{h}^{k+1}\big)
		+  \Delta t\big((\phi_{h}^{k+1})^{3}-\phi_{h}^{k},  D_{t}\phi_{h}^{k+1} \big)
		- \Delta t\big(\omega_{h}^{k+1},  D_{t}\phi_{h}^{k+1} \big) = 0. 
	\end{align}
	With the help of the identity 
	\begin{align}
		&2(a-b, a)=a^2-b^2+(a-b)^2, \label{equality:1}\\
		&(a^3-b)(a-b)=\frac{1}{4}[(a^2-1)^2-(b^2-1)^2]+ \frac{1}{4} (a^2-b^2)^2+\frac{1}{2}a^2(a-b)^2+\frac{1}{2}(a-b)^2,
	\end{align}
	combining the equations \eqref{ener:1}-\eqref{ener:2}, we obtain:
	\begin{equation}\label{ener:phi}
		\begin{aligned}
			& \Delta t \|\nabla\omega_{h}^{k+1}\|^2+\frac{1}{2}(\|\nabla\phi_{h}^{k+1}\|^2- \|\nabla\phi_{h}^{k}\|^2 +\|\nabla(\phi_{h}^{k+1}-\phi_{h}^{k})\|^2)+\frac{1}{4 }[\| (\phi_{h}^{k+1})^2-1\|^2-\| (\phi_{h}^{k})^2-1)\|^2] \\   
			&\,+ \frac{1}{4 } \|(\phi_{h}^{k+1})^2-(\phi_{h}^{k})^2\|^2 +\frac{1}{2 }\|\phi_{h}^{k+1}\|^2\|\phi_{h}^{k+1}-\phi_{h}^{k}\|^2+ \frac{1}{2 } \|\phi_{h}^{k+1}-\phi_{h}^{k}\|^2+ \Delta t\big(\nabla\phi_{h}^{k}\cdot \tilde{\u}_{h}^{k+1},  \omega_{h}^{k+1}\big) =0,
		\end{aligned}
	\end{equation}
	Applying $\C_h=\Delta t \B_h^{k+1}$ in \eqref{B-1} and using \eqref{equality:1}, we have 
	\begin{equation}
		\frac{1}{2}(\| \B_{h}^{k+1}\|^2- \|\B_{h}^{k}\|^2 +\| \B_{h}^{k+1}-\B_{h}^{k} \|^2) +  \Delta t  \|\nabla\times \B_{h}^{k+1}\|^2 -\Delta t (\tilde{\u}_{h}^{k+1}\times \B_{h}^{k}, \nabla\times \B_h^{k+1})+  \Delta t \|\nabla\cdot\B_{h}^{k+1}\|^2 =0. 
	\end{equation}	
	Letting $\v_{h}=\Delta t \tilde{\u}_{h}^{k+1}$ in \eqref{equau} and applying \eqref{equality:1}, we derive 
	\begin{equation}\label{err:umid}
		\frac{1}{2}( \|\tilde{\u}_{h}^{k+1}\|^2-\|\u_{h}^{k}\|^2+ \|\tilde{\u}_{h}^{k+1}-\u_{h}^{k}\|^2 ) +\Delta t \|\nabla\tilde{\u}_{h}^{k+1}\|^2 +\Delta t( \tilde{\u}_{h}^{k+1}, \nabla p_{h}^{k}) +\Delta t (\nabla\times \B_{h}^{k+1}, \tilde{\u}_{h}^{k+1}\times \B_{h}^{k})-\Delta t(\omega_{h}^{k+1}  \nabla \phi_{h}^{k}, \tilde{\u}_{h}^{k+1})=0.
	\end{equation}
	Equation \eqref{numerical-scheme3} can be rewritten as the following form
	\begin{align}
		\u_{h}^{k+1}+\Delta t \nabla p_{h}^{k+1}= \tilde{\u}_{h}^{k+1}+\Delta t \nabla p_{h}^{k}, 
	\end{align}
then taking the $L^2$ inner product with itself, we get   
	\begin{equation}\label{sum:u}
		\begin{aligned}
			\Delta t(\tilde{\u}_{h}^{k+1}, \nabla  p_{h}^{k})=\frac{1}{2}(\|\u_{h}^{k+1}\|^2-\|\tilde{\u}_{h}^{k+1}\|^2)+  \frac{1}{2} \Delta t^2(\|\nabla p_{h}^{k+1}\|^2-\|\nabla p_{h}^{k}\|^2 ). 
		\end{aligned}
	\end{equation}
	
	By making the summations of \eqref{ener:phi}-\eqref{sum:u}, we obtain:
	\begin{equation}\label{ener:sum}
		\begin{aligned}
			&\frac{1}{2}(\|\nabla\phi_{h}^{k+1}\|^2- \|\nabla\phi_{h}^{k}\|^2 +\|\nabla(\phi_{h}^{k+1}-\phi_{h}^{k})\|^2)
			+\frac{1}{2}( \|\u_{h}^{k+1}\|^2-\|\u_{h}^{k}\|^2+ \|\tilde{\u}_{h}^{k+1}-\u_{h}^{k}\|^2 )\\
			&\, +\frac{1}{4 }[\| (\phi_{h}^{k+1})^2-1\|^2-\| (\phi_{h}^{k})^2-1\|^2] + \frac{1}{4 } \|(\phi_{h}^{k+1})^2-(\phi_{h}^{k})^2\|^2+  \Delta t \|\nabla\omega_{h}^{k+1}\|^2+\Delta t  \|\nabla\times \B_{h}^{k+1}\|^2 + \Delta t \|\nabla\cdot\B_{h}^{k+1}\|^2\\ 
			&\,+\frac{1}{2 }\|\phi_{h}^{k+1}\|^2\|\phi_{h}^{k+1}-\phi_{h}^{k}\|^2 + \frac{1}{2 } \|\phi_{h}^{k+1}-\phi_{h}^{k}\|^2   +  \Delta t\|\nabla \tilde{\u}_{h}^{k+1}\|^2  +\frac{1}{2}(\| \B_{h}^{k+1}\|^2- \|\B_{h}^{k}\|^2 +\| \B_{h}^{k+1}-\B_{h}^{k} \|^2) 
			= 0. 
		\end{aligned}
	\end{equation} 
	This concludes the proof. 
\end{proof}

\begin{Remark}[\cite{2013AnalysisDIEGEL}, the boundedness of $\| \nabla \phi_{h}^{k}\|_{L^4}$ ]
The solution $(\phi_{h}^{k+1}, \omega_{h}^{k+1}, \u_{h}^{k+1}, p_{h}^{k+1}, \B_{h}^{k+1})\in \boldsymbol \chi_{h}^{r}$ of the fully discrete scheme \eqref{numerical-scheme1}-\eqref{numerical-scheme3} satisfies the following error estimate for all $\Delta t, h>0$, 
\begin{equation}
\max\limits_{0<k<K-1}\| \nabla \phi_{h}^{k}\|_{L^4}\leq C_0. 
\end{equation}
 
\end{Remark}

\section{Projections and their properties}\label{sec-projection}



(1) The classic Ritz projection $R_{h}: H^{1}(\Omega)\rightarrow S_{h}^{r}$ is defined by \cite{wheeler1973priori},
\begin{equation*}
(\nabla (\varphi-R_{ h}\varphi), \nabla\psi_{h})=0,
\end{equation*}
for all $\psi_{h}\in S_{h}^{r}$ and $\int_{\Omega}(\varphi-R_{ h}\varphi) \mathrm{d}\x$=0.  And the Ritz projection   is equipped with following estimates:
\begin{equation*}
\begin{aligned}
&\|\varphi-R_{ h}\varphi\|_{L^{s}}+h\|\varphi-R_{  h}\varphi\|_{W^{1, s}}\leq C_{0}h^{r+1}\|\varphi\|_{W^{r+1, s}},\\
&\|\varphi-R_{ h}\varphi\|_{H^{-1}}\leq C_{0}\beta_{h}\|\varphi\|_{H^{r+1}},\\
&\|d_{t}(\varphi^{k+1}-R_{  h}\varphi^{k+1})\|+h\|d_{t}(\varphi^{k+1}-R_{  h}\varphi^{k+1})\|_{H^{1}}\leq C_{0}h^{r+1}\|d_{t}\varphi^{k+1}\|_{H^{r+1}},\\
&\|d_{t}(\varphi^{k+1}-R_{ h}\varphi^{k+1})\|_{H^{-1}}\leq C_{0}\beta_{h}\|d_{t}\varphi^{k+1}\|_{H^{r+1}},
\end{aligned}
\end{equation*}
for $s\in[2, \infty]$, $k=0, 1, 2, \cdots, K-1$, and $\beta_{h}$ is defined in equation \eqref{bta}.

(2) The Ritz quasi-projection $\tilde{R}_{h}: H^{1}(\Omega)\rightarrow S_{h}^{r}$ is defined by \cite{2023Optimalwang},
\begin{equation*}
(\nabla (\omega-\tilde{R}_{h}\omega), \nabla\xi_{h})+  (\nabla (\phi-R_{ h}\phi)\cdot \u,   \xi_{h})=0,
\end{equation*}
for all $\xi_{h}\in S_{h}^{r}$ and $\int_{\Omega}(\omega-\tilde{R}_{h}\omega) \mathrm{d}\x$=0. And the Ritz quasi-projection   is equipped with the following estimates:
\begin{equation*}
\begin{aligned}
&\|\omega-\tilde{R}_{h}\omega\| +h\|\nabla(\omega-\tilde{R}_{h}\omega)\|\leq C_{0}h^{r+1}(\|\boldsymbol{u}\|_{L^{\infty}}\|\phi\|_{H^{r+1}}+ \|\omega\|_{H^{r+1}}),\\
&\|\omega-\tilde{R}_{h}\omega\|_{H^{-1}}\leq C_{0}\beta_{h}(\|\u\|_{W^{1, 4}}\|\phi\|_{H^{r+1}}+ \|\omega\|_{H^{r+1}}),\\
&\|\nabla(d_{t}(\omega^{k+1}-\tilde{R}_{h}\omega^{k+1}) )\|\leq C_{0}h^{r}(  \|\u^{k+1}\|_{L^{\infty}}\|d_{t}\phi^{k+1}\|_{H^{r}}+ \|d_{t}\u^{k+1}\|_{L^{\infty}}\|\phi^{k}\|_{H^{r}} +\|d_{t}\omega^{k+1}\|_{H^{r+1}}   ),\\
&\|d_{t}(\omega^{k+1}-\tilde{R}_{h}\omega^{k+1})\|_{H^{-1}}\leq C_{0}\beta_{h}(  \|\u^{k+1}\|_{W^{1,4}}\|d_{t}\phi^{k+1}\|_{H^{r+1}}+ \|d_{t}\u^{k+1}\|_{W^{1, 4}}\|\phi^{k}\|_{H^{r+1}} +\|d_{t}\omega^{k+1}\|_{H^{r+1}}   ),
\end{aligned}
\end{equation*}
for $k=0, 1, 2, \cdots, K-1$.

(3) The Stokes quasi-projection $(\P_{h}, P_{h}): \H_{0}^{1}(\Omega)\times L_{0}^{2}(\Omega)\rightarrow (\X_{h}^{r+1}+ \nabla \mathring{S}_{h}^{1}) \times \mathring{S}_{h}^{r}/ (\X_{h}^{1b}+ \nabla \mathring{S}_{h}^{1})\times \mathring{S}_{h}^{1}$ is defined by \cite{2023Optimalwang},
\begin{equation*}
\begin{aligned}
&(\nabla (\u-\P_{h}(\u, p)), \nabla \v_{h})- (p-P_{h}(\u, p), \nabla\cdot \v_{h}  )=( \omega \nabla( \phi- R_{  h}\phi), \v_{h} ),\\
&(\nabla\cdot (\u-\P_{h}(\u, p)), q_{h})=0,
\end{aligned}
\end{equation*}
for all $(\v_{h}, q_{h})\in (\X_{h}^{r+1}+ \nabla \mathring{S}_{h}^{1}) \times \mathring{S}_{h}^{r}/ (\X_{h}^{1b}+ \nabla \mathring{S}_{h}^{1}) \times \mathring{S}_{h}^{1}$. After this, we denote $\P_{h}\u:=\P_{h}(\u, p)$ and $P_{h}p:=P_{h}(\u, p)$ for simplicity. And the  Stokes quasi-projection has the following estimates: 
\begin{equation}\label{Stokes1}
\|\u-\P_{h}\u\|\leq \left\{
\begin{aligned}
&C_{0}h^{r+2} (\|\u\|_{H^{r+2}}+  \|p\|_{H^{r+1}} +  \|\phi\|_{H^{r+1}}\|\omega\|_{H^{2}}),\quad &&r\geq2      ,\\
&C_{0}h^{r+1} (\|\u\|_{H^{r+1}}+  \|p\|_{H^{r}} +  \|\phi\|_{H^{r+1}}\|\omega\|_{H^{2}}),\quad &&r=1,
\end{aligned}
\right.
\end{equation}
\begin{equation*}\label{Stokes2}
\|\nabla (\u-\P_{h}\u)\|+\|p-P_{h}p\|\leq\left\{
\begin{aligned}
&C_{0}h^{r+1} (\|\u\|_{H^{r+2}}+  \|p\|_{H^{r+1}} +  \|\phi\|_{H^{r+1}}\|\omega\|_{W^{1, 4}}),\quad &&r\geq2,\\
&C_{0}h^{r+1} (\|\u\|_{H^{r+2}}+  \|p\|_{H^{r+1}} +  \|\phi\|_{H^{r+1}}\|\omega\|_{W^{1, 4}}),\quad  &&r=1 \quad (\X_{h}^{2}\times \mathring{S}_{h}^{1}),\\
&C_{0}h^{r} (\|\u\|_{H^{r+1}}+  \|p\|_{H^{r}} +  \|\phi\|_{H^{r}}\|\omega\|_{W^{1, 4}}),\quad  &&r=1\quad (\X_{h}^{1b}\times \mathring{S}_{h}^{1}),
\end{aligned}
\right.
\end{equation*}
\begin{equation*}\label{Stokes3}
\|d_{t} (\u^{k+1} -\P_{h}\u^{k+1} )\| \leq \left\{
\begin{aligned}
&C_{0}h^{r+2} (\|d_{t}\u^{k+1} \|_{H^{r+2}}+  \|d_{t}p^{k+1} \|_{H^{r+1}} +  \|d_{t}\phi^{k+1} \|_{H^{r+1}}\|\omega^{k+1} \|_{H^{2}}+\|\phi^{k} \|_{H^{r+1}}\|d_{t}\omega^{k+1} \|_{H^{2}}   ),\quad  &&r\geq2,\\
&C_{0}h^{r+1}(\|d_{t}\u^{k+1} \|_{H^{r+1}}+  \|d_{t}p^{k+1} \|_{H^{r}} +  \|d_{t}\phi^{k+1} \|_{H^{r+1}}\|\omega^{k+1} \|_{H^{2}}+\|\phi^{k} \|_{H^{r+1}}\|d_{t}\omega^{k+1} \|_{H^{2}}   ),\quad &&r=1 ,
\end{aligned}
\right.
\end{equation*}

\begin{Remark}\label{UMAX}
Based on the estimates in (\ref{Stokes1}), we have the following boundedness:
\begin{equation*}
\|\P_{h}\u\|_{L^{\infty}}+\|\P_{h}\u\|_{W^{1, 3}}\leq C_{0}(\|\u\|_{H^{2}}+ \|p\|_{H^{2}} +\|\phi\|_{H^{2}}\|\omega\|_{H^{2}}  ).
\end{equation*}
\end{Remark}

(4) The $L^{2}$ projection $I_{h}: L^{2}(\Omega)\rightarrow S_{h}^{r}$, and $\I_{h}: \L^{2}(\Omega)\rightarrow (\X_{h}^{r+1}+ \nabla \mathring{S}_{h}^{1})/(\X_{h}^{1b}+ \nabla \mathring{S}_{h}^{1})$ are defined as follows,
\begin{equation*}
\begin{aligned}
& (v -I_{h}v, \xi_{h})=0, \quad \forall\, \xi_{h}\in S_{h}^{r},\\
& (\v-\I_{h}\v, \v_{h})=0, \quad \forall\, \v_{h}\in \X_{h}^{r+1}/\X_{h}^{1b}.\\
\end{aligned}
\end{equation*}
Based on the above classic $L^{2}$ projection, the following estimates hold,
\begin{equation*}
\begin{aligned}
& \|v -I_{h}v \| +h \|\nabla (v -I_{h}v)\|\leq C_{0} h^{r+1}\|v\|_{H^{r+1}},\\
& \|\v -\I_{h}\v \| +h \|\nabla (\v -\I_{h}\v)\|\leq C_{0} h^{r+2}\|\v\|_{H^{r+2}},\quad && \text{if}\, \I_{h}\v \in \X_{h}^{r+1},  \\
& \|\v -\I_{h}\v\| +h \|\nabla (\v-\I_{h}\v)\|\leq C_{0} h^{2}\|\v\|_{H^{2}},\quad && \text{if}\, \I_{h}\v\in \X_{h}^{1b}.
\end{aligned}
\end{equation*}

(5) The Maxwell projection $\bm{\Pi}_{h}: \H_{\tau}^{1}(\Omega)\rightarrow \Y_{h}^{r+1}$ is defined as follows:
\begin{equation*}
(\nabla\times (\B-\bm{\Pi}_{h}\B), \nabla\times \C_{h})+(\nabla\cdot (\B-\bm{\Pi}_{h}\B), \nabla\cdot \C_{h} )=0, \quad \forall\, \C_{h}\in \Y_{h}^{r+1}.
\end{equation*}
The following estimates hold for the Maxwell projection:
\begin{equation*}
\begin{aligned}
& \|\B-\bm{\Pi}_{h}\B\|+h\| \B-\bm{\Pi}_{h}\B\|_{H^{1}}\leq C_{0}h^{r+2}\| \B\|_{H^{r+2}}.
\end{aligned}
\end{equation*}

\section{The Proof of Theorem \ref{theorem2-1} }\label{sec-error}
The well-posedness of the convex-splitting algorithm for the two-phase MHD model \eqref{model} has been given in \cite{2019A}. In this section, we present the proof of Theorem \ref{theorem2-1}. To this end, we shall introduce the following discrete Gronwall inequality \cite{heywood1990finite}.

\begin{lemma}\label{Gronwall}
Let $\alpha_{n}, \beta_{n}, c_{n}, \gamma_{n}$ and $g_{0}$ be a sequence of nonnegative numbers for integers $n\geq 0$ such that
\begin{equation*}
\alpha_{n}+\Delta t\sum_{j=0}^{n}\beta_{j}\leq \Delta t\sum_{j=0}^{n}\gamma_{j}\alpha_{j}+\Delta t\sum_{j=0}^{n}c_{j}+g_{0}.
\end{equation*}
Assume that $\gamma_{j}\Delta t< 1$ for all $j$, and set $\sigma_{j}=(1-\gamma_{j}\Delta t)^{-1}$. Then, for all $n\geq 0$,
\begin{equation*}
\alpha_{n}+\Delta t\sum_{j=0}^{n}\beta_{j}\leq \rm{exp}\left(\Delta t\sum_{j=0}^{n}\sigma_{j} \gamma_{j} \right)\left(  \Delta t\sum_{j=0}^{n}c_{j}+g_{0}\right).
\end{equation*}
\end{lemma}

For simplicity, we introduce the following notations,
\begin{equation*}
\begin{aligned}
&e_{\phi}^{k+1}:=R_{  h}\phi^{k+1}-\phi^{k+1}_{h},\;  e_{\omega}^{k+1}:=\tilde{R}_{h}\omega^{k+1}-\omega^{k+1}_{h}, \; e_{p}^{k+1}:=P_{h}p^{k+1}-p^{k+1}_{h},\;\\ &\tilde{e}_{\u}^{k+1}:=\P_{h}\u^{k+1}-\tilde{\u}^{k+1}_{h}, \; e_{\u}^{k+1}:=\P_{h}\u^{k+1}-\u^{k+1}_{h}, \; e_{\B}^{k+1}:=\bm{\Pi}_{h}\B^{k+1}-\B^{k+1}_{h}.
\end{aligned}
\end{equation*}

With the help of projection operators defined in the previous section, we subtract (\ref{mdel1})-(\ref{model5}) from (\ref{eqn123})-(\ref{equau}) and  \eqref{numerical-scheme3} to get the following error equations for ($e_{\phi}^{k+1}, e_{\omega}^{k+1}, e_{\u}^{k+1}, e_{p}^{k+1}, e_{\B}^{k+1}$),

\begin{subequations}\label{error-equation}
\begin{align}
\left(D_{t}e_{\phi}^{k+1}, \xi_{h}  \right) +\left(\nabla e_{\omega}^{k+1}, \nabla  \xi_{h} \right)=\left(D_{t}(R_{  h}\phi^{k+1}- \phi^{k+1}  ), \xi_{h}\right)& +\left(\nabla \phi_{h}^{k}\cdot\tilde{\u}_{h}^{k+1},   \xi_{h}  \right)-\left( \nabla R_{ h}\phi^{k+1}\cdot\u^{k+1},   \xi_{h} \right)+\left(R_{1}^{k+1}, \xi_{h}\right), \label{error:phi}\\
\left( \nabla e_{\phi}^{k+1}, \nabla\psi_{h}  \right)-\left(e_{\omega}^{k+1}, \psi_{h} \right)+\frac{1}{2}\left( Z^{k+1}(e_{\phi}^{k+1}-e_{\phi}^{k}  ),  \psi_{h}    \right)=   &-\left(Z^{k+1}\bar{e}_{\phi}^{k+\frac{1}{2}},  \psi_{h} \right)+\left(e_{\phi}^{k},  \psi_{h}\right)-\left((\phi^{k+1})^{3}-(R_{ h} \phi^{k+1})^{3} ,  \psi_{h}\right)\nonumber\\
 &+\left(\phi^{k}-R_{ h}\phi^{k}, \psi_{h}\right)+\left(\omega^{k+1}-\tilde{R}_{h}\omega^{k+1}, \psi_{h}\right)+\left(R_{2}^{k+1} ,\psi_{h} \right), \label{error:w}\\
 \left(D_{t}e_{\B}^{k+1}, \C_{h}\right)+ \left(\nabla\times e_{\B}^{k+1},  \nabla\times  \C_{h}  \right)+\left( \nabla\cdot e_{\B}^{k+1},   \nabla \cdot \C_{h}  \right)=&\left(D_{t}(  \bm{\Pi}_{h}\B^{k+1}-\B^{k+1}), \C_{h} \right)+\left(R_{3}^{k+1} ,\C_{h}\right) \nonumber \\
 &-\left( ( \tilde{\u}_{h}^{k+1}\times \B_{h}^{k}, \nabla\times \C_{h})-  (\u^{k+1}\times \B^{k}, \nabla\times \C_{h})\right), \label{error:B}\\
\left(\frac{\tilde{e}_{\u}^{k+1}-e_{\u}^{k  } }{\Delta t}, \v_{h}\right)+\left(\nabla  \tilde{e}_{\u}^{k+1}, \nabla\v_{h} \right)+ \left(  \v_{h}, \nabla e_{p}^{k }\right)=&\left (b(\u_{h}^{k}, \tilde{\u}_{h}^{k+1}, \v_{h} )-b(\u^{k+1}, \u^{k+1}, \v_{h} )\right)\nonumber \\
&   -\left( (\nabla \phi_{h}^{k}\cdot \v_{h}, \omega_{h}^{k+1}) -(\nabla R_{ h}\phi^{k+1}\cdot \v_{h}, \omega^{k+1})  \right)\nonumber  \\
&  + (D_{t}(\P_{h}\u^{k+1}-\u^{k+1}), \v_{h})    +\left( R_{4}^{k+1} ,\boldsymbol{v}_{h}  \right) \nonumber\\
&   +\left((\nabla\times \B_{h}^{k+1}, \v_{h}\times\B_{h}^{k}  )  -(\nabla\times \B^{k+1}, \v_{h}\times\B^{k}  )  \right) ,\label{error:u}\\
\frac{e_{\u}^{k+1}}{\Delta t}+\nabla e_{p}^{k+1}=&\frac{\tilde{e}_{\u}^{k+1}}{\Delta t}+\nabla e_{p}^{k}+  R_{5}^{k+1}, \label{error:uend} 
\end{align}
\end{subequations}
where $(\xi_{h}, \psi_{h}, \v_{h}, q_{h},  \C_{h})\in \boldsymbol \chi_{h}^{r}$  and $k=0, 1, \cdots, K-1$. We need to pay attention to
\begin{equation*}
\begin{aligned}
&\bar{e}_{\phi}^{k+\frac{1}{2}}:=\frac{1}{2}(e_{\phi}^{k+1}+e_{\phi}^{k} ),\quad e_{\phi}^{k+1}=\bar{e}_{\phi}^{k+\frac{1}{2}}+\frac{1}{2}(e_{\phi}^{k+1}-e_{\phi}^{k}).
\end{aligned}
\end{equation*}
And we have 
\begin{align*}
	&(R_h\phi^{k+1})^{3}-(\phi_{h}^{k+1})^{3}\\
	&=3e_{\phi}^{k+1}\int_{0}^{1}\left((1-\theta)\phi_{h}^{k+1} +\theta R_h \phi^{k+1}   \right)^{2}d\theta   ,\\
	&= \frac{3}{2}(e_{\phi}^{k+1}-e_{\phi}^{k} ) \int_{0}^{1}\left((1-\theta)\phi_{h}^{k+1} +\theta R_h \phi^{k+1}   \right)^{2}d\theta+3\bar{e}_{\phi}^{k+\frac{1}{2}}  \int_{0}^{1}\left((1-\theta)\phi_{h}^{k+1} +\theta R_h \phi^{k+1}   \right)^{2}d\theta\\
	&= \frac{1}{2}(e_{\phi}^{k+1}-e_{\phi}^{k} )Z^{k+1}+\bar{e}_{\phi}^{k+\frac{1}{2}} Z^{k+1},
\end{align*}
where $Z^{k+1}:=3\int_{0}^{1}\left((1-\theta)\phi_{h}^{k+1} +\theta R_h \phi^{k+1}   \right)^{2}d\theta.$

In addition,  $R_{1}^{k+1}, R_{2}^{k+1}, R_{3}^{k+1}, R_{4}^{k+1}$, and $R_{5}^{k+1}$ are the truncation errors satisfying
\begin{equation*}
\begin{aligned}
(R_{1}^{k+1}, \xi_{h})=
&\big(D_{t}\phi^{k+1}, \xi_{h}\big) -  \big(\phi_t^{k+1}, \xi_{h}\big),\\
(R_{2}^{k+1}, \psi_{h})= & \big(\phi^k, \psi_{h}\big) -\big(\phi^{k+1}, \psi_{h}\big),	\\
(R_{3}^{k+1}, \C_{h})= 
&\big(D_{t}\B^{k+1}, \C_{h}\big) - \big(\B_{t}^{k+1}, \C_{h}\big) 
- \big(\u^{k+1}\times \B^{k} , \nabla\times \C_{h}\big) + \big(\u^{k+1}\times \B^{k+1}, \nabla\times \C_{h}\big),\\ 
(R_{4}^{k+1}, \v_{h})=
&\big(D_{t}\u^{k+1}, \v_{h}\big) - \big(\u_{t}^{k+1}, \v_{h}\big)+ \big(\nabla\times \B^{k+1}, \v_{h}\times\B^{k} \big) - \big(\nabla\times \B^{k+1}, \v_{h}\times\B^{k+1} \big)+ \big(\v_{h}, \nabla (p^{k}-p^{k+1}) \big), \\
R_{5}^{k+1} =&   \nabla (e_{p}^{k+1}-e_{p}^{k}). 
\end{aligned}
\end{equation*}
By Taylor expansion, we have the following estimates:
\begin{equation}\label{Taylor}
\begin{aligned}
&\left(\Delta t\sum_{k=0}^{K-1}(\|R_{1}^{k+1}\|^{2}+ \|R_{2}^{k+1}\|^{2} +\|R_{3}^{k+1}\|^{2}+\|R_{4}^{k+1}\|^{2} +\|R_{5}^{k+1}\|^{2}   )     \right)^{\frac{1}{2}}\leq C_{0}\Delta t,
\end{aligned}
\end{equation}
where $C_{0}$ is a generic positive constant. 

The following lemma will be  used in subsequent proof. 
Next, we will give the error estimates of the numerical solutions below \cite{2023Optimalwang}. 

\begin{lemma}\label{lemmaphi}
The fully discrete scheme \eqref{numerical-scheme1}-\eqref{numerical-scheme3} is uniquely solvable and unconditionally stable, and satisfies the following estimate:
\begin{align}\label{phi}
(i) \quad \|\nabla e_{\phi}^{k+1}\|^{2}+\Delta t\sum_{m=0}^{k}\|\nabla e_{\omega}^{m+1}\|^{2}\leq & C_{\varepsilon}\Delta t\sum_{m=0}^{k}\left(\|\tilde{e}_{\u}^{m+1}\|^{2} + \|\nabla e_{\phi}^{m+1}\|^{2} + \|  e_{\phi}^{m}\|^{2}    \right)+C_{\varepsilon} \left(\beta_{h}^{2}+\Delta t^{2}  \right) +\varepsilon \Delta t \sum_{m=0}^{k} \|  e_{\omega}^{m+1}\|^{2}        \nonumber  \\
&   + C_{0}\Delta t \sum_{m=0}^{k}\|D_{t}Z^{m+1}\|_{L^{\frac{3}{2}}}\| e_{\phi}^{m}\|^{2}_{H^{1}}+ \left(2+\varepsilon \right)\|e_{\phi}^{k+1}\|^{2}.
\end{align}
We should note that $\varepsilon$ is a very small generic positive constant, and $C_{\varepsilon}$ is a generic positive constant that depends on $\varepsilon$.  They may take different values at different places.
By taking $\xi_{h}=(-\Delta_{h})^{-1}e_{\phi}^{k+1}$
and $\psi_{h}=e_{\phi}^{k+1}-\frac{1}{|\Omega|}( e_{\phi}^{k+1}, 1)$ in equations (\ref{error:phi})-(\ref{error:w}), there exists a positive constant $\Delta t_{1}$ such that when $\Delta t \leq\Delta t_{1}$,  the $L^{2}$-norm estimates for $e_{\phi}^{k+1}$ and $e_{\omega}^{k+1}$ are as follows:
\begin{subequations}
\begin{align}
&(ii)\quad \|e_{\omega}^{k+1}\|\leq C_{0}\left(\|\nabla e_{\omega}^{k+1}\|+  \|\nabla e_{\phi}^{k+1}\| +\|\nabla e_{\phi}^{k}\| +\beta_{h}+\Delta t \right),\\
&(iii)\quad \| e_{\phi}^{k+1}\|^{2}\leq \varepsilon \|\nabla e_{\phi}^{k+1}\|^{2}+C_{\varepsilon}\Delta t\sum_{m=0}^{k}\|e_{\textbf{u}}^{k+1}\|^{2}+ C_{\varepsilon} (\beta_{h}^{2}+\Delta t^{2}),
\end{align}
\end{subequations}
where $\beta_{h}$ is defined in  Theorem \ref{theorem2-1}. Here, the discrete Laplacian operator $\Delta_{h}:\mathring{S}_{h}^{r}\rightarrow \mathring{S}_{h}^{r}$ is denoted by
	\begin{equation*}
		(-\Delta_{h}\psi_{h}, \xi_{h})=(\nabla\psi_{h}, \nabla\xi_{h}),\quad \forall\, \psi_{h}, \xi_{h}\in \mathring{S}_{h}^{r}.
	\end{equation*}
	In addition, if $v_{h}$ is constant, we define  $(-\Delta_{h})^{\frac{1}{2}}v_{h}= (-\Delta_{h})^{-\frac{1}{2}}v_{h}:=0$.

\end{lemma}

\subsection{Estimates for $e_{\textbf{B}}^{k+1}$.}

Taking $\C_{h}=e_{\B}^{k+1}$ in equation (\ref{error:B}), we have
\begin{align}\label{error-B}
	\begin{aligned}
		&\frac{1}{2}D_{t}\|e_{\B}^{k+1}\|^{2}+\frac{1}{2\Delta t}\|e_{\B}^{k+1}-e_{\B}^{k}\|^{2}+\|\nabla \times e_{\B}^{k+1}\|^{2}+\|\nabla \cdot e_{\B}^{k+1}\|^{2} \\
		&=\left(D_{t}(  \bm{\Pi}_{h}\B^{k+1}-\B^{k+1}), e_{\B}^{k+1} \right)+\left(R_{3}^{k+1} ,e_{\B}^{k+1}\right)     -\left( (\tilde{\u}_{h}^{k+1}\times \B_{h}^{k}, \nabla\times e_{\B}^{k+1})-  (\u^{k+1}\times \B^{k}, \nabla\times e_{\B}^{k+1})\right)\\
		&:=\sum_{i=1}^{3}Q_{i}.
	\end{aligned}
\end{align}
According to the regularity assumption (\ref{regularity}), the Taylor expansion, and the estimates of the Maxwell projection, we obtain 
\begin{equation}\label{Q1}
	\begin{aligned}
		Q_{1}&=\left(D_{t}(  \bm{\Pi}_{h}\B^{k+1}-\B^{k+1}), e_{\B}^{k+1} \right)\\
		&\leq \|D_{t}(  \bm{\Pi}_{h}\B^{k+1}-\B^{k+1})\| \|e_{\B}^{k+1}\| \\
		&\leq C_{0}\left((h^{r+2})^{2}+  \|e_{\B}^{k+1}\|^{2}\right), \\
		Q_{2}&=\left(R_{3}^{k+1} ,e_{\B}^{k+1}\right)  \leq C_{0}\left( \Delta t^{2}+ \|e_{\B}^{k+1}\|^{2}\right).
	\end{aligned}
\end{equation}
Thus, combining estimates in \eqref{Q1}, we derive 
\begin{align}\label{Q12}
	Q_{1}+Q_{2}\leq C_{0}\left((h^{r+2})^{2}+ \Delta t^{2}+ \|e_{\B}^{k+1}\|^{2}\right).
\end{align}
By further expanding $Q_{3}$, we have
\begin{align*}
	Q_{3}\leq& \left(\u^{k+1}\times (\B^{k}-\bm{\Pi}_{h}\B^{k}   ), \nabla\times e_{\B}^{k+1}  \right) +\left(\u^{k+1}\times e_{\B}^{k}, \nabla\times e_{\B}^{k+1}  \right)  &(\text{term}\quad Q_{3,1})\nonumber \\
	&+ \left( (\u^{k+1}-\P_{ h}\u^{k+1}  )\times \B_{h}^{k}, \nabla\times e_{\B}^{k+1}  \right) &(\text{term}\quad Q_{3,2})\nonumber\\
	&+\left(\tilde{e}_{\u}^{k+1}\times\B_{h}^{k}, \nabla\times e_{\B}^{k+1}\right).  
\end{align*}

With the help of the regularity assumption (\ref{regularity}) and Remark \ref{UMAX}, we consider the following estimate  
\begin{align*}
	Q_{3,1}=&\left(\u^{k+1}\times (\B^{k}-\bm{\Pi}_{h}\B^{k}   ), \nabla\times e_{\B}^{k+1}  \right) +\left(\u^{k+1}\times e_{\B}^{k}, \nabla\times e_{\B}^{k+1}  \right) \\
	&\leq \|\u^{k+1}\|_{L^{\infty}} \|\B^{k}-\bm{\Pi}_{h}\B^{k}\| \|\nabla\times e_{\B}^{k+1}  \| +\|\u^{k+1}|\|_{L^{\infty}} \|e_{\B}^{k}\| \| \nabla\times e_{\B}^{k+1}\|\\
	&\leq C_{\varepsilon}\left(\| e_{\B}^{k}\|^{2}+ \beta_{h}^{2} \right)+\varepsilon \|\nabla\times e_{\B}^{k+1} \|^{2},\\
	Q_{3,2}=&-\left( (\u^{k+1}-\P_{ h}\u^{k+1}  )\times e_{\B}^{k}, \nabla\times e_{\B}^{k+1}  \right)
	+\left((\u^{k+1}-\P_{ h}\u^{k+1}  )\times (\bm{\Pi}_{h}\B^{k}-\B^{k}), \nabla\times e_{\B}^{k+1}\right) \nonumber \\
	&+\left( (\u^{k+1}-\P_{ h}\u^{k+1}  )\times\B^{k}, \nabla\times e_{\B}^{k+1} \right)\nonumber \\
	\leq & \|\u^{k+1}-\P_{ h}\u^{k+1} \|_{L^{\infty}}\left(\| e_{\B}^{k}\|+\|\bm{\Pi}_{h}\B^{k}-\B^{k}\|  \right)\|\nabla\times e_{\B}^{k+1} \|+ \|\u^{k+1}-\P_{ h}\u^{k+1} \|\,\|\B^{k}\|_{L^{\infty}}\|\nabla\times e_{\B}^{k+1}\| \nonumber\\
	\leq &C_{\varepsilon}\left(\| e_{\B}^{k}\|^{2}+ \beta_{h}^{2} \right)+\varepsilon \|\nabla\times e_{\B}^{k+1} \|^{2}.
\end{align*}
And obviously, we estimate $Q_{3}$  as 
\begin{align}\label{Q3}
	Q_{3}\leq C_{\varepsilon}\left(\| e_{\B}^{k}\|^{2}+ \beta_{h}^{2} \right)+\varepsilon \|\nabla\times e_{\B}^{k+1} \|^{2}- \left(\nabla\times e_{\B}^{k+1}\times\B_{h}^{k}, \tilde{e}_{\u}^{k+1} \right).
\end{align}
For a sufficiently small $\varepsilon$, by combining the inequalities (\ref{Q12}) and (\ref{Q3}),  the equation (\ref{error-B}) reduces to
\begin{align*}
	D_{t}\|e_{\B}^{k+1}\|^{2} +\|\nabla \times e_{\B}^{k+1}\|^{2} +\|\nabla \cdot e_{\B}^{k+1}\|^{2} \leq C_{\varepsilon}\left(\| e_{\B}^{k}\|^{2}+\|e_{\B}^{k+1}\|^{2}+ \beta_{h}^{2} + \Delta t^{2}\right)  - \left(\nabla\times e_{\B}^{k+1}\times\B_{h}^{k}, \tilde{e}_{\u}^{k+1} \right).
\end{align*}
Then, by summing the result from time step $t_{0}$ to $t_{k}$, we get
\begin{align}
	\begin{aligned}
		&\|e_{\B}^{k+1}\|^{2} +\Delta t\sum_{m=0}^{k}\left(\|\nabla \times e_{\B}^{m+1}\|^{2} +\|\nabla \cdot e_{\B}^{m+1}\|^{2} \right)  \\
		& \leq C_{\varepsilon}\Delta t\sum_{m=0}^{k} \| e_{\B}^{m+1}\|^{2}+C_{\varepsilon}\left( \beta_{h}^{2} + \Delta t^{2} \right)  - \Delta t\sum_{m=0}^{k}\left(\nabla\times e_{\B}^{m+1}\times\B_{h}^{m}, \tilde{e}_{\u}^{m+1} \right).
	\end{aligned}
\end{align}

\subsection{Estimates for $\tilde{e}_{\textbf{u}}^{k+1}$. }

Taking $\v_{h}=\tilde{e}_{\u}^{k+1}$  in equations \eqref{error:u}, we have
\begin{align}\label{error-u}
	\begin{aligned}
&\frac{1}{2\Delta t}\left( \|\tilde{e}_{\u}^{k+1}\|^2-\|e_{\u}^{k}\|^2+\|\tilde{e}_{\u}^{k+1}- e_{\u}^{k}\|^2 \right)  +\|\nabla \tilde{e}_{\u}^{k+1}\|^{2} + (\tilde{e}_{\u}^{k+1}, \nabla e_p^{k} )\\
&= \left (b(\u_{h}^{k}, \tilde{\u}_{h}^{k+1}, \tilde{e}_{\u}^{k+1} )-b(\u^{k+1}, \u^{k+1}, \tilde{e}_{\u}^{k+1})\right) & (\text{term}\quad I_{1}) \\
& \quad  -\left( (\nabla \phi_{h}^{k}\cdot \tilde{e}_{\u}^{k+1}, \omega_{h}^{k+1}) -(\nabla R_{  h}\phi^{k+1}\cdot \tilde{e}_{\u}^{k+1}, \omega^{k+1})  \right) & (\text{term}\quad I_{2})    \\
&\quad    +\left((D_{t}(\P_{h}\u^{k+1}-\u^{k+1}), \tilde{e}_{\u}^{k+1})  \right) +\left( R_{4}^{k+1} , \tilde{e}_{\u}^{k+1} \right) & (\text{term}\quad I_{3})  \\
&\quad  +\left((\nabla\times \B_{h}^{k+1}, \tilde{e}_{\u}^{k+1}\times\B_{h}^{k}  )  -(\nabla\times \B^{k+1},  \tilde{e}_{\u}^{k+1} \times\B^{k} )  \right). & (\text{term}\quad I_{4})
\end{aligned}
\end{align}


We give the estimate of term $I_1$ as
\begin{equation}\label{err:I1}
\begin{aligned}
	I_1=&  b(\u_{h}^{k}, \tilde{\u}_{h}^{k+1}, \tilde{e}_{\u}^{k+1} )-b(\u^{k+1}, \u^{k+1}, \tilde{e}_{\u}^{k+1}) \\
	=&b(\u^{k+1}-\u^{k}, \u^{k+1}, \tilde{e}_{\u}^{k+1})
	+b(\u^{k}-\P_h \u^{k}, \u^{k+1}, \tilde{e}_{\u}^{k+1})
	+b(e_{\u}^{k}, \P_h\u^{k+1}, \tilde{e}_{\u}^{k+1})\\
	&  + b(\P_h\u^{k},  \u^{k+1}-\P_h\u^{k+1}, \tilde{e}_{\u}^{k+1}) 
	+ b( \u_h^{k},  \tilde{e}_{\u}^{k+1}, \tilde{e}_{\u}^{k+1})\\
	=& b(\u^{k+1}-\u^{k}, \u^{k+1}, \tilde{e}_{\u}^{k+1})
	+b(\u^{k}-\P_h \u^{k}, \u^{k+1}, \tilde{e}_{\u}^{k+1})
	+b(e_{\u}^{k}, \P_h\u^{k+1}, \tilde{e}_{\u}^{k+1})+ b(\P_h\u^{k},  \u^{k+1}-\P_h\u^{k+1}, \tilde{e}_{\u}^{k+1})\\
	\leq & \Delta t b(D_t\u^{k+1}, \u^{k+1}, \tilde{e}_{\u}^{k+1})+\|\u^{k}-\P_h \u^{k}\|\, \|\u^{k+1}\|_{L^{\infty}} \, \|\tilde{e}_{\u}^{k+1}\|
	+\|e_{\u}^{k}\|\, (\|\P_h\u^{k+1}\|_{W^{1,3}} +\|\P_h\u^{k+1}\|_{L^{\infty}})\,\|\tilde{e}_{\u}^{k+1}\| \\
	&+(\|\P_h\u^{k+1}\|_{W^{1,3}} +\|\P_h\u^{k+1}\|_{L^{\infty}})\, \|\u^{k+1}-\P_h\u^{k+1}\| \, \|\tilde{e}_{\u}^{k+1}\| \\
	\leq &C_{\varepsilon} (\Delta t^2 +\beta_{h}^2+\|e_{\u}^{k}-\tilde{e}_{\u}^{k+1}\|  )
	+\varepsilon \|\nabla \tilde{e}_{\u}^{k+1}\|^2,
	\end{aligned} 
\end{equation} 
where we consider $\|e_{\u}^{k}\| \leq \|e_{\u}^{k}-\tilde{e}_{\u}^{k+1}\|+ \| \tilde{e}_{\u}^{k+1}\|$.
Then, we have the estimate of term $I_2$,
\begin{equation}\label{err:I2}
	\begin{aligned}
		I_2&=-\left( (\nabla \phi_{h}^{k}\cdot \tilde{e}_{\u}^{k+1}, \omega_{h}^{k+1}) -(\nabla R_{  h}\phi^{k+1}\cdot \tilde{e}_{\u}^{k+1}, \omega^{k+1})  \right)\\
		&= (\nabla \phi_{h}^{k}\cdot \tilde{e}_{\u}^{k+1}, e_\omega^{k+1}) 
		+  (\nabla \phi_{h}^{k}\cdot \tilde{e}_{\u}^{k+1}, \omega^{k+1}-\tilde{R}_{h}\omega^{k+1} ) 
		+(\nabla e_{\phi}^{k}\cdot \tilde{e}_{\u}^{k+1}, \omega^{k+1}) 
		+ (\nabla (R_{h}\phi^{k+1}-R_{h}\phi^{k}   )\cdot \tilde{e}_{\u}^{k+1}, \omega^{k+1})\\
		&\leq C_0 \|\nabla \phi_{h}^{k}\|_{L^4}\|\tilde{e}_{\u}^{k+1}\|_{L^4} \|e_\omega^{k+1}\|+ C_{\varepsilon} \left(\|\nabla e_{\phi}^{k} \|^2+ \beta_{h}^2\right)+ \varepsilon\|\nabla \tilde{e}_{\u}^{k+1}\|^2
		+C_0\|\nabla e_{\phi}^{k}\|\,\|\tilde{e}_{\u}^{k+1}\|\,\|\omega^{k+1} \|_{L^{\infty}}\\
		&\quad + C_0 \|\nabla (R_{h}\phi^{k+1}-R_{h}\phi^{k}   )\|\, \|\tilde{e}_{\u}^{k+1}\|\,\|\omega^{k+1}\|_{L^{\infty}}\\
		&\leq C_{\varepsilon}\left(\|\nabla \tilde{e}_{\u}^{k+1}\|^{2} + \|\nabla e_{\phi}^{k}\|^{2} + \beta_{h}^{2}+\Delta t^{2}    \right)+  \varepsilon\left( \| e_{\omega}^{k+1}\|^{2}+  \|\nabla \tilde{e}_{\u}^{k+1}\|^{2}\right),
	\end{aligned}
\end{equation}
where we use the known result  of $\| \nabla \phi_{h}^{k}\|_{L^4}\leq C_0$. Next, we consider the estimate of term $I_3$ 
\begin{equation}
	\begin{aligned}\label{err:I3}
		I_3 &= \left((D_{t}(\P_{h}\u^{k+1}-\u^{k+1}), \tilde{e}_{\u}^{k+1})  \right) +\left( R_{4}^{k+1} , \tilde{e}_{\u}^{k+1} \right)  \\
		&\leq C_{\varepsilon}\left( \beta_{h}^{2}  +\Delta t^2  \right)+\varepsilon \|\tilde{e}_{\u}^{k+1}\|^2. 
	\end{aligned}
\end{equation}

For the last term $I_4$, we consider 
\begin{equation}\label{I5}
\begin{aligned}
I_{4}&=\left(\nabla\times e_{\B}^{k+1}\times \B_{h}^{k},  \tilde{e}_{\u}^{k+1}\right)+  \left(\nabla\times (\B^{k+1}-\bm{\Pi}_{h} \B^{k+1} )\times \B_{h}^{k},  \tilde{e}_{\u}^{k+1} \right) \\
&\quad +\left( \nabla\times \B^{k+1}\times  e_{\B}^{k}, \tilde{e}_{\u}^{k+1}\right) +\left(\nabla\times\B^{k+1}\times(\B^{k}-\bm{\Pi}_{h} \B^{k}), \tilde{e}_{\u}^{k+1}\right)\\
&:=\sum_{j=1}^{4}I_{4, j}.
\end{aligned}
\end{equation}
By the inequations (\ref{i1})-(\ref{i4}) and the estimates of Maxwell projection, we obtain 
\begin{equation}\label{eq:I42}
\begin{aligned}
I_{4,2}= &\left(\nabla\times (\B^{k+1}-\bm{\Pi}_{h} \B^{k+1} )\times (\B_{h}^{k}-\bm{\Pi}_{h} \B^{k} ), \tilde{e}_{\u}^{k+1}\right) + \left(\nabla\times (\B^{k+1}-\bm{\Pi}_{h} \B^{k+1} )\times (\bm{\Pi}_{h} \B^{k}- \B^{k} ), \tilde{e}_{\u}^{k+1} \right)\\
& +\left(\nabla\times (\B^{k+1}-\bm{\Pi}_{h} \B^{k+1} )\times \B^{k}, \tilde{e}_{\u}^{k+1}\right) \\
\leq &\|\nabla\times (\B^{k+1}-\bm{\Pi}_{h} \B^{k+1} )\|\,\|e_{\B}^{k}\|_{L^{3}}\|\tilde{e}_{\u}^{k+1}\|_{L^{6}}+ \|\nabla\times (\B^{k+1}-\bm{\Pi}_{h} \B^{k+1} )\|\,\|\bm{\Pi}_{h} \B^{k}- \B^{k}\|_{L^{3}}\|\tilde{e}_{\u}^{k+1}\|_{L^{6}}\\
&+C_{0}\| \B^{k+1}-\bm{\Pi}_{h}\B^{k+1}  \|\,\|\nabla\B^{k}\|_{L^{\infty}}\|\tilde{e}_{\u}^{k+1}\|+C_{0}\| \B^{k+1}-\bm{\Pi}_{h} \B^{k+1}  \|\,\|\B^{k}\|_{L^{\infty}}\|\nabla \tilde{e}_{\u}^{k+1}\|\\
\leq &C_{0}h^{-\frac{d}{6}} \|\nabla  (\B^{k+1}-\bm{\Pi}_{h} \B^{k+1} )\|\|e_{\B}^{k}\|\, \|\nabla \tilde{e}_{\u}^{k+1}\| +C_{0}\|\nabla (\B^{k+1}-\bm{\Pi}_{h} \B^{k+1} )\|\,\|\bm{\Pi}_{h} \B^{k}- \B^{k}\|_{L^{3}}\|\nabla \tilde{e}_{\u}^{k+1}\| \\
&+C_{\varepsilon} (h^{r+2})^{2}+ \varepsilon\|\nabla \tilde{e}_{\u}^{k+1}\|^{2}\\
\leq & C_{0}h^{1-\frac{d}{6}}\|e_{\B}^{k}\|\, \|\nabla e_{\u}^{k+1}\| +C_{\varepsilon}(h^{r+2})^{2}+\varepsilon \|\nabla \tilde{e}_{\u}^{k+1}\|^{2}\\
\leq & C_{\varepsilon}(h^{r+2})^{2}+C_{\varepsilon}\|e_{\B}^{k}\|^{2}+\varepsilon\|\nabla \tilde{e}_{\u}^{k+1}\|^{2}.
\end{aligned}
\end{equation}
\begin{equation}\label{eq:I43}
	\begin{aligned}
		I_{4,3}&=\left( \nabla\times \B^{k+1}\times  e_{\B}^{k}, \tilde{e}_{\u}^{k+1}\right)\\
		&\leq \|\nabla\times \B^{k+1}\|_{L^3}  \|e_{\B}^{k} \| \|\tilde{e}_{\u}^{k+1}\|_{L^6} \\
		&\leq C_{\varepsilon}\|e_{\B}^{k}\|^{2}+\varepsilon\|\nabla \tilde{e}_{\u}^{k+1}\|^{2}.
	\end{aligned} 
\end{equation}
The estimates of  $I_{4,4}$ follows a similar procedure as the one described above \eqref{eq:I42},  so we omit the detailed steps here. By substituting the  estimates \eqref{eq:I42} and \eqref{eq:I43} into equation (\ref{I5}), we obtain 
\begin{equation}\label{I5END}
I_{4}\leq \left(\nabla\times e_{\B}^{k+1}\times \B_{h}^{k},  \tilde{e}_{\u}^{k+1}\right)+ C_{\varepsilon}(h^{r+2})^{2}+C_{\varepsilon}\|e_{\B}^{k}\|^{2}+\varepsilon\|\nabla \tilde{e}_{\u}^{k+1}\|^{2}.
\end{equation}
Combining the equations \eqref{err:I1}, \eqref{err:I2},\eqref{err:I3}, and  the above inequality (\ref{I5END}), and for a sufficiently small $\varepsilon$, the equation (\ref{error-u}) reduces to
\begin{equation}\label{error:u-mid}
	\begin{aligned}
&\frac{1}{2\Delta t}\left( \|\tilde{e}_{\u}^{k+1}\|^2-\|e_{\u}^{k}\|^2+\|\tilde{e}_{\u}^{k+1}- e_{\u}^{k}\|^2 \right)  +\|\nabla \tilde{e}_{\u}^{k+1}\|^{2}+ (\tilde{e}_{\u}^{k+1}, \nabla e_p^{k} )\\
&\leq   C_{\varepsilon}\left(\|\nabla \tilde{e}_{\u}^{k+1}\|^{2} + \|e_{\u}^{k}- \tilde{e}_{\u}^{k+1}\|^{2} + \|\nabla e_{\phi}^{k}\|^{2}+\|e_{\B}^{k}\|^{2}+\beta_{h}^{2} +\Delta t^{2}  \right)  +\varepsilon \|e_{\omega}^{k+1}\|^{2}+\left(\nabla\times e_{\B}^{k+1}\times \B_{h}^{k},  \tilde{e}_{\u}^{k+1}\right).
\end{aligned}
\end{equation}

\subsection{Estimates for $e_{\textbf{u}}^{k+1}$. }
We rewrite \eqref{error:uend} to derive 
\begin{equation}
e_{\u}^{k+1} +\Delta t\nabla e_{p}^{k+1}=  \tilde{e}_{\u}^{k+1} +\Delta t\nabla e_{p}^{k}+  \Delta t  R_{5}^{k+1}.
\end{equation}
By  taking the $L^2$ inner product of equation \eqref{error:uend} on both sides and using \eqref{div-free}, we obtain
\begin{equation}\label{error:u-end}
	\begin{aligned}
		&\frac{1}{2\Delta t}(\|e_{\u}^{k+1}\|^2 -\|\tilde{e}_{\u}^{k+1}\|^2) + \frac{\Delta t}{2} \left(\|\nabla  e_{p}^{k+1}\|^2- \|\nabla e_{p}^{k} \|^2 \right) \\
		&\leq C_{\varepsilon} \Delta t^2
		+ ( \tilde{e}_{\u}^{k+1}, \nabla e_{p}^{k})
		+ ( \tilde{e}_{\u}^{k+1},   R_{5}^{k+1})
		+\Delta t ( R_{5}^{k+1}, \nabla e_{p}^{k}) \\
		&\leq \varepsilon \left( \|\tilde{e}_{\u}^{k+1}\|^2   +\|e_{\u}^{k+1}\|^2\right) + C_{\varepsilon} \left(\Delta t^2 +  \Delta t ^2\|\nabla e_{p}^{k}\|^2  \right )+ ( \tilde{e}_{\u}^{k+1}, \nabla e_{p}^{k}). 
	\end{aligned} 
\end{equation}
Combine \eqref{error:u-mid} and \eqref{error:u-end}, we have
\begin{equation} 
	\begin{aligned}
		&\frac{1}{2\Delta t}\left( \|e_{\u}^{k+1}\|^2-\|e_{\u}^{k}\|^2+\|\tilde{e}_{\u}^{k+1}- e_{\u}^{k}\|^2 \right)  +\|\nabla \tilde{e}_{\u}^{k+1}\|^{2}  + \frac{\Delta t}{2} \left(\|\nabla  e_{p}^{k+1}\|^2- \|\nabla e_{p}^{k} \|^2 \right)\\
		&\leq    C_{\varepsilon}\left(\|\nabla \tilde{e}_{\u}^{k+1}\|^{2} + \|e_{\u}^{k}- \tilde{e}_{\u}^{k+1}\|^{2} + \|\nabla e_{\phi}^{k}\|^{2}+ \|e_{\B}^{k}\|^{2} +  \Delta t ^2\|\nabla e_{p}^{k}\|^2 +\beta_{h}^{2} +\Delta t^{2}\right)   +\varepsilon \|e_{\omega}^{k+1}\|^{2}+\left(\nabla\times e_{\B}^{k+1}\times \B_{h}^{k},  \tilde{e}_{\u}^{k+1}\right).
	\end{aligned}
\end{equation}

By summing up the above estimate from time step $t_{0}$ to $t_{k}$, we  obtain
\begin{align}\label{u}
	\begin{aligned}
		&\frac{1}{2}\|e_{\u}^{k+1}\|^{2} +\frac{1}{2 } \|\tilde{e}_{\u}^{m+1}- e_{\u}^{m}\|^2 + \Delta t \|\nabla \tilde{e}_{\u}^{m+1}\|^{2}+ \frac{\Delta t^2}{2} \|\nabla  e_{p}^{m+1}\|^2 \\
		&\leq C_{\varepsilon}\Delta t \sum_{m=0}^{k} \left(\|\nabla \tilde{e}_{\u}^{m+1}\|^{2} + \|  \tilde{e}_{\u}^{m+1}-e_{\u}^{m}\|^{2}  + \|\nabla e_{\phi}^{m+1}\|^{2}+  \|e_{\B}^{m+1}\|^{2} +  \Delta t ^2\|\nabla e_{p}^{m}\|^2 \right)+C_{\varepsilon} \left(\beta_{h}^{2} +\Delta t^{2}\right)    \\
		&\quad +\varepsilon \Delta t \sum_{m=0}^{k} \|e_{\omega}^{m+1}\|^{2}+\Delta t \sum_{m=0}^{k}  \left(\nabla\times e_{\B}^{m+1}\times \B_{h}^{m},  \tilde{e}_{\u}^{m+1}\right),
	\end{aligned}
\end{align}
noting that $\|e_{\u}^{0}\|^{2}\leq C_{0}\beta_{h}$.

By virtue of (ii)-(iii) in Lemma \ref{lemmaphi}, by adding inequalities  (\ref{phi}), (\ref{u}) and the above inequality, we have the following estimate
\begin{align*}
&\|\nabla e_{\phi}^{k+1}\|^{2}+\|e_{\u}^{k+1}\|^{2}+\|e_{\B}^{k+1}\|^{2}+\Delta t\sum_{m=0}^{k}\left(\|\nabla e_{\omega}^{m+1}\|^{2} +\|\nabla \times e_{\B}^{m+1}\|^{2} +\|\nabla \cdot e_{\B}^{m+1}\|^{2} \right)\nonumber \\
& \leq C_{\varepsilon}\Delta t\sum_{m=0}^{k}\left(\|e_{\u}^{m+1}\|^{2} + \|\nabla e_{\phi}^{m+1}\|^{2} + \|e_{\B}^{m+1}\|^{2}    \right)+ C_{\varepsilon} \left(\beta_{h}^{2}+\Delta t^{2}  \right)
+ C_{0}\Delta t \sum_{m=0}^{k}\|D_{t}Z^{m+1}\|_{L^{\frac{3}{2}}}\| e_{\phi}^{m}\|^{2}_{H^{1}}.
\end{align*}
By using the Lemma \ref{Gronwall} ( discrete Gronwall's inequality), there exists a positive constant $\Delta t_{2}$ such that, if $\Delta t\leq \min\{ \Delta t_{1}, \Delta t_{2} \}$,  
\begin{align}
&\|\nabla e_{\phi}^{k+1}\|^{2}+\|e_{\u}^{k+1}\|^{2}+\|e_{\B}^{k+1}\|^{2}+\Delta t\sum_{m=0}^{k}\left(\|\nabla e_{\omega}^{m+1}\|^{2} +\|\nabla \times e_{\B}^{m+1}\|^{2} +\|\nabla \cdot e_{\B}^{m+1}\|^{2} \right) \leq   C_{\varepsilon} \left(\beta_{h}^{2}+\Delta t^{2}  \right).
\end{align}
By employing the above inequality, (ii)-(iii) in  Lemma \ref{lemmaphi} and the estimates of projection operators in Section \ref{sec-projection}, we give the error estimates in Theorem \ref{theorem2-1}.

\section{Numerical examples}\label{sec-examples}

In this section, we conduct several 2D/3D numerical examples to verify the theoretical analysis using the
finite element software FreeFem \cite{Hecht+2012+251+266}. We solve the scheme \eqref{numerical-scheme1}-\eqref{numerical-scheme3} with the two examples of finite element spaces as follows:
\begin{table}[hpt]
\setlength{\belowcaptionskip}{0.1cm} \caption{The finite element spaces.} \label{elements}
\centering
\begin{tabular}{ccccccccccccccccccccccccccccccc}
\toprule
 & $ \phi_h $ &    $ \u_h$   & $p_h$ & $\B_h$ &\\ 
\midrule
 case I \;  & $P_{1}$   &    $P_{1}^{b}$    &  $P_{1}$   &   $P_{1}$    & \\
 case II\;  & $P_{1}$   &   $P_{2}$    &  $P_{1}$   &   $P_{2}$    & \\
\bottomrule
\end{tabular}
\end{table}

\subsection{2D/3D convergence of the scheme}

In this subsection, we simulate the 2D convergence results of the numerical scheme in a square domain $\Omega=[0,1]^{d}$ using the following smooth exact solutions:
\begin{eqnarray*}
\left\{
\begin{aligned}
\phi&= \rm \cos(t)\cos^2(\pi x)\cos^2(\pi y),\\
\u&=\rm \cos(t)\Big(\pi \sin(2\pi y) \sin^2(\pi x) ,\ -\rm \pi \sin(2\pi x) \sin^2(\pi y) \Big),\\
p&=\rm \cos(t)(2x-2) (2y -1),\\
\B&=\rm \cos(t)\Big(\sin(\pi x)\cos(\pi y),  -\rm \sin(\pi y) \cos(\pi x)\Big), \\
\end{aligned}
\right.
\end{eqnarray*}
3D smooth solutions as
\begin{eqnarray*}
\left\{
\begin{aligned}
\phi&= \rm \exp(-2t)\sin^2(\pi x)\sin^2(\pi y) \sin^2(\pi z),\\
\u&=\rm \exp(t)\Big(y(1-y)z(1-z) ,\ x(1-x)z(1-z),\ x(1-x)y(1-y) \Big),\\
p&=\rm \exp(t)(2x-1) (2y -1) (2z-1),\\
\B&=\rm \exp(t)\Big(\sin(\pi y)\sin(\pi z),\  \rm \sin(\pi x) \sin(\pi z),\ \rm \sin(\pi x)\sin(\pi y) \Big). \\
\end{aligned}
\right.
\end{eqnarray*}
The parameters are chosen as
\begin{equation*}
\gamma=1,\quad M=1,\quad \nu=1,\quad \mu=1,\quad \lambda=1,\quad \sigma=1.
\end{equation*}

For simplicity, we verify the time and  space convergence orders at the end time $T=1$ using the relationship  $\Delta t=O(h^{2})$ between the time step and space step for case I, and  $\Delta t=O(h^{3})$  for case II. The numerical results for case I in 2D and 3D are given in  Tables  \ref{case-I}  and \ref{case-I-3D}.  Table \ref{case-II} shows the results of case II in 2D case.   From the numerical results, the convergence orders of the scheme are consistent with the theoretical results, which are shown in Theorem \ref{theorem2-1}. 
\begin{table}[hpt]
\setlength{\belowcaptionskip}{0.1cm}
\caption{ Convergence results  with case I  in 2D.}\label{case-I}
\centering
\begin{tabular}{ccccccccccccccccccccccccccccccc}
\toprule
$h$ & $\|\phi-\phi_{h}\|$ & rate & $\|\nabla (\phi-\phi_{h})\|$ & rate & $\|\boldsymbol{u}-\boldsymbol{u}_{h}\|$ & rate & $\|\nabla (\boldsymbol{u}-\boldsymbol{u}_{h})\|$ & rate &\\ 
\midrule
 1/8 &2.66e-01 & 1.48  &3.62e-01 &  1.20   &1.30e-01 & 1.66 &3.02e-01 & 0.92  & \\
 1/16 &7.48e-02 & 1.83 &1.57e-01 & 1.21  &3.44e-02 & 1.92   &1.52e-01 & 0.99  & \\
 1/32 &1.93e-02 &1.95  &7.36e-02 & 1.09  &8.67e-03 & 1.99 &7.58e-02 & 1.00  &\\
 1/64 &4.86e-03 & 1.99 &3.61e-02 & 1.03  &2.17e-03 & 2.00  &3.78e-02 & 1.00  &\\
\bottomrule
\toprule
$h$ &$\|\boldsymbol{B}-\boldsymbol{B}_{h}\|$ & rate & $\|\nabla (\boldsymbol{B}-\boldsymbol{B}_{h})\|$ & rate & $\|p-p_{h}\|$ & rate & \\ 
\midrule
  1/8 &1.52e-02  & 1.58  & 1.95e-01 &  0.97  & 3.28e-00 &1.45  & \\
  1/16 &4.09e-03  & 1.89   &9.80e-02  & 0.99   & 1.04e-00 & 1.65  & \\
  1/32 &1.04e-03  & 1.97   &4.91e-02 &  1.00 &3.36e-01 & 1.64 &\\
  1/64 &2.62e-04 &1.99 &2.45e-02 & 1.00  &1.13e-01  & 1.57  & \\
\bottomrule
\end{tabular}
\end{table}

\begin{table}[hpt]
\setlength{\belowcaptionskip}{0.1cm}
\caption{  Convergence results  with case II  in 2D.}\label{case-II}
\centering
\begin{tabular}{ccccccccccccccccccccccccccccccc}
\toprule
$h$ & $\|\phi-\phi_{h}\|$ & rate & $\|\nabla (\phi-\phi_{h})\|$ & rate & $\|\boldsymbol{u}-\boldsymbol{u}_{h}\|$ & rate & $\|\nabla (\boldsymbol{u}-\boldsymbol{u}_{h})\|$ & rate &\\ 
\midrule
 1/8 &2.66e-01 & 1.48  &3.62e-01 &  1.20   &5.92e-03 & 3.11 &4.47e-02  & 1.87  & \\
 1/16 &7.48e-02 & 1.83 &1.57e-01 & 1.21  &9.41e-04 & 2.65   &1.15e-02 & 1.96  & \\
 1/32 &1.93e-02 &1.95  &7.36e-02 & 1.09  &2.04e-04  & 2.20 &2.89e-03 & 1.99  &\\
 1/64 &4.86e-03 & 1.99 &3.61e-02 & 1.03  & 4.90e-05  & 2.05  &7.23e-04 & 2.00  &\\
\bottomrule
\toprule
$h$ &$\|\boldsymbol{B}-\boldsymbol{B}_{h}\|$ & rate & $\|\nabla (\boldsymbol{B}-\boldsymbol{B}_{h})\|$ & rate & $\|p-p_{h}\|$ & rate & \\ 
\midrule
  1/8 &3.50e-04  & 3.45  &1.50e-02 &  1.96  & 1.73e-00 &1.16  & \\
  1/16 & 5.49e-05    & 2.67   &3.79e-03  & 1.99   &5.23e-01 & 1.73  & \\
  1/32 & 1.31e-05  & 2.07  &9.50e-04  &  2.00 &1.38e-01  & 1.93 &\\
  1/64 & 3.27e-06   & 2.00 &2.38e-04  & 2.00  &3.48e-02   & 1.98  & \\
\bottomrule
\end{tabular}
\end{table}

\begin{table}[hpt] 
\setlength{\belowcaptionskip}{0.1cm}
\caption{ Convergence results  with case I in 3D.}\label{case-I-3D}
\centering
\begin{tabular}{ccccccccccccccccccccccccccccccc}
\toprule
$h$ & $\|\phi-\phi_{h}\|$ & rate & $\|\nabla (\phi-\phi_{h})\|$ & rate & $\|\boldsymbol{u}-\boldsymbol{u}_{h}\|$ & rate & $\|\nabla (\boldsymbol{u}-\boldsymbol{u}_{h})\|$ & rate &\\ 
\midrule
1/4 &1.09e-00  &    &9.59e-01 &      &1.30e-01 &   &4.10e-01  &    & \\
 1/8 &4.37e-01 & 1.32 &4.54e-01 & 1.08  &3.58e-02  & 1.86   &2.07e-01 & 0.99  & \\
 1/12 &2.17e-01 &1.73  &2.74e-01 & 1.25  &1.68e-02 & 1.87 &1.37e-01 & 1.02  &\\
 1/16 &1.27e-01& 1.86 &1.93e-01 & 1.21  &9.66e-03 & 1.92  &1.02e-01 & 1.02  &\\
\bottomrule
\toprule
$h$ &$\|\boldsymbol{B}-\boldsymbol{B}_{h}\|$ & rate & $\|\nabla (\boldsymbol{B}-\boldsymbol{B}_{h})\|$ & rate & $\|p-p_{h}\|$ & rate & \\ 
\midrule
  1/4 &1.73e-01 &    &3.83e-01 &     & 2.25e-00&  & \\
  1/8 &4.90e-02 & 1.82   &1.95e-01 & 0.97  &7.17e-01 & 1.65  & \\
  1/12 &2.23e-02  & 1.94   &1.31e-01  &  0.99 &3.32e-01  & 1.90 &\\
  1/16 &1.27e-02 & 1.97 &9.81e-02 & 1.00  &1.90e-01 & 1.95  & \\
\bottomrule
\end{tabular}
\end{table}
 
\subsection{Spinodal decomposition}
The spinodal decomposition is a phase separation phenomenon that occurs in binary or multi-component alloys, polymer blends and liquid crystals \cite{2023Energy, shi2024structure}.
The computational domain is $\Omega=[0,1]^{2}$.  The initial values read as  
\begin{eqnarray}\label{fai0}
\phi_{0}=-0.05+0.001\mathrm{ rand }(x),\quad
\u_{0}=\0,\quad
p_{0}=0,\quad
\B_{0}=\0,
\end{eqnarray}
where $\mathrm{rand}(x)$ is a uniformly distributed random function in $[-1,1]$ with zero mean. We select finite element pairs  cases I  to test the spinodal decomposition phenomenon. The parameters are given as
\begin{equation*}
\gamma=1/100,\quad M=1,\quad \nu=1,\quad \mu=1,\quad \lambda=1,\quad \sigma=1.
\end{equation*}

\begin{figure}[htbp]
	\centering
\subfigure[$t=0.0001$]{
		\begin{minipage}[t]{0.2\linewidth}
			\centering
			\includegraphics[width=\textwidth]{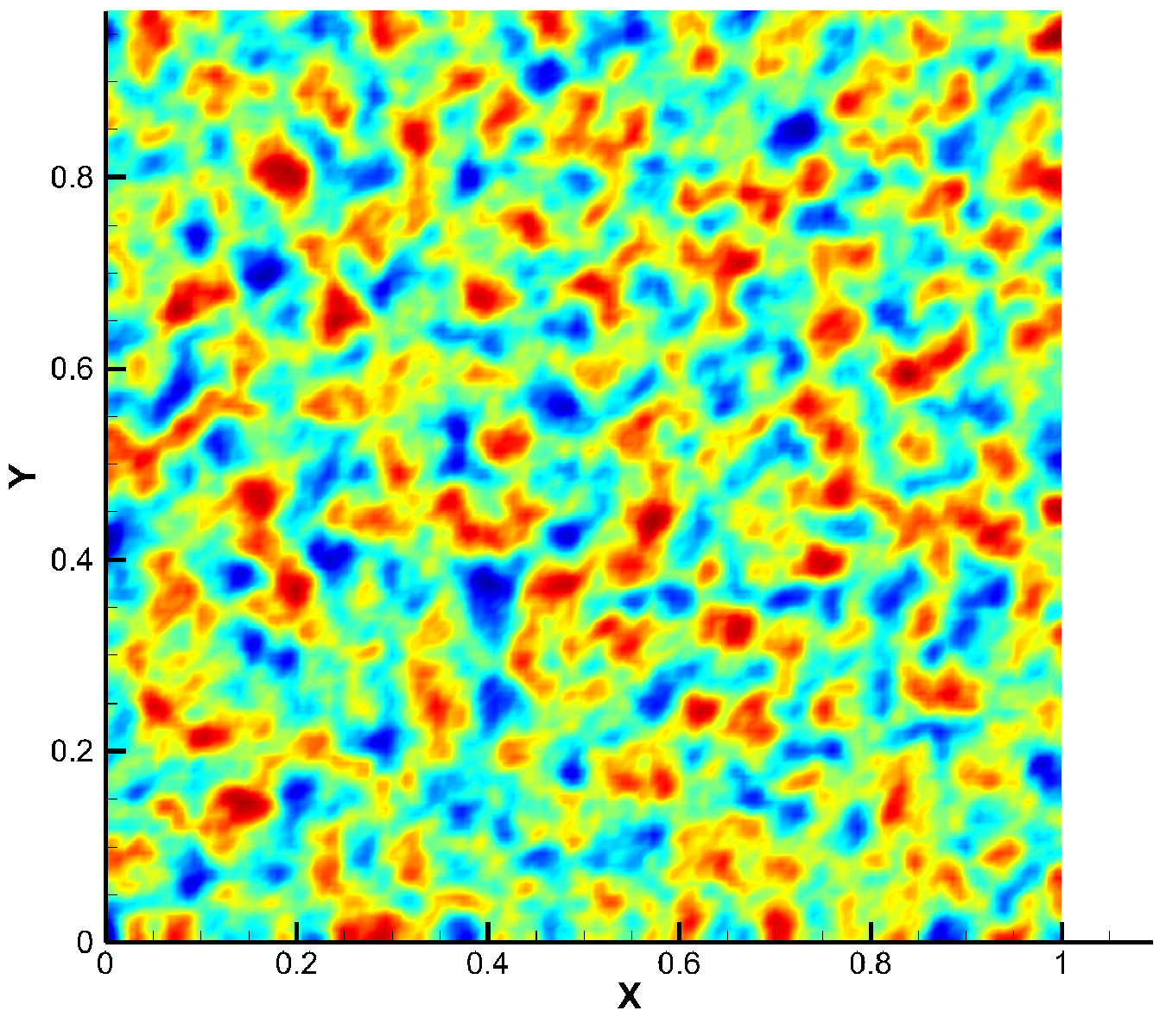}
		\end{minipage}
	}%
	\subfigure[$t=0.05$]{
		\begin{minipage}[t]{0.2\linewidth}
			\centering
			\includegraphics[width=\textwidth]{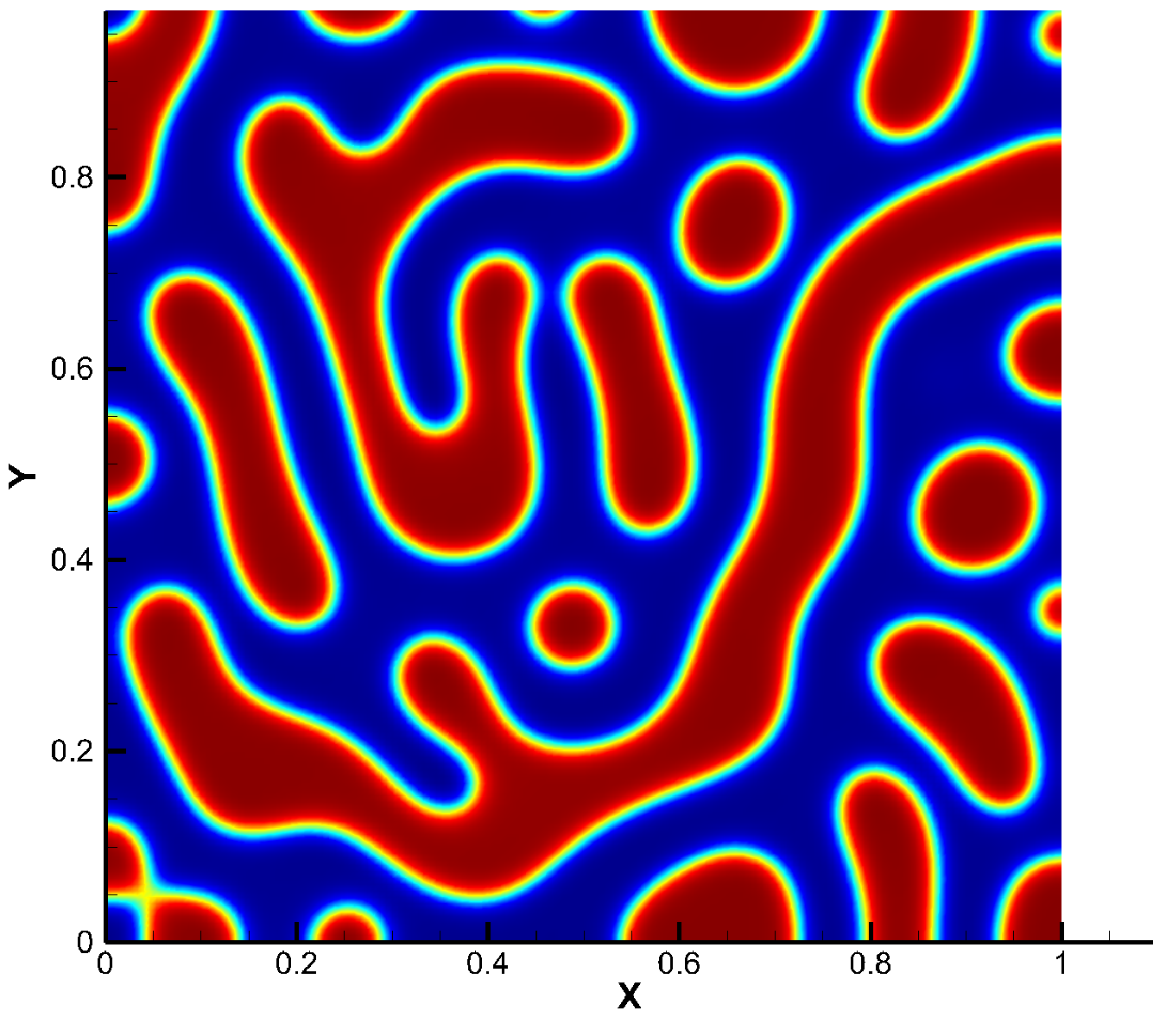}
		\end{minipage}
	}%
	\subfigure[$t=0.5$]{
		\begin{minipage}[t]{0.2\linewidth}
			\centering
			\includegraphics[width=\textwidth]{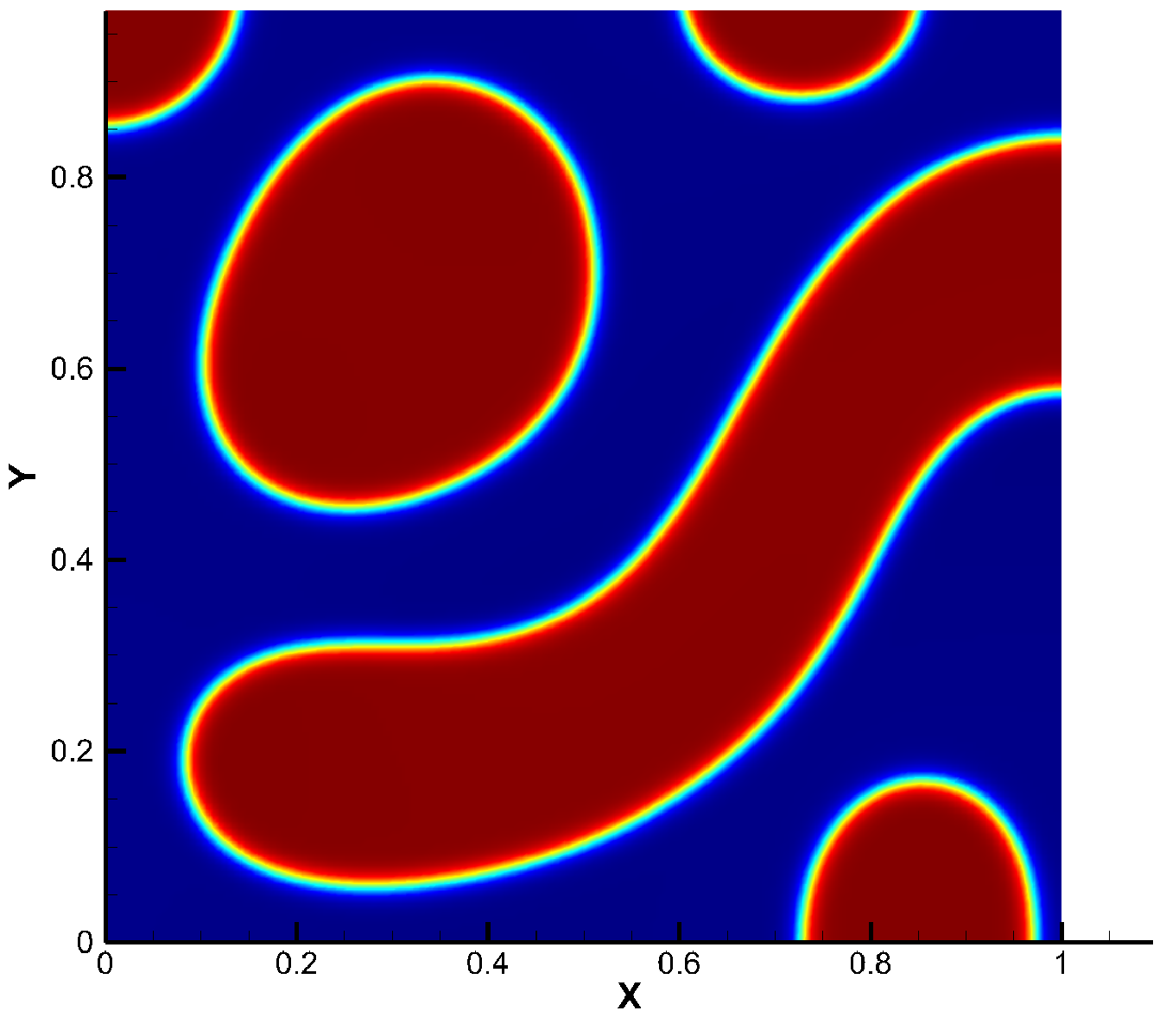}
		\end{minipage}
	}%
	\subfigure[$t=2.5$]{
		\begin{minipage}[t]{0.2\linewidth}
			\centering
			\includegraphics[width=\textwidth]{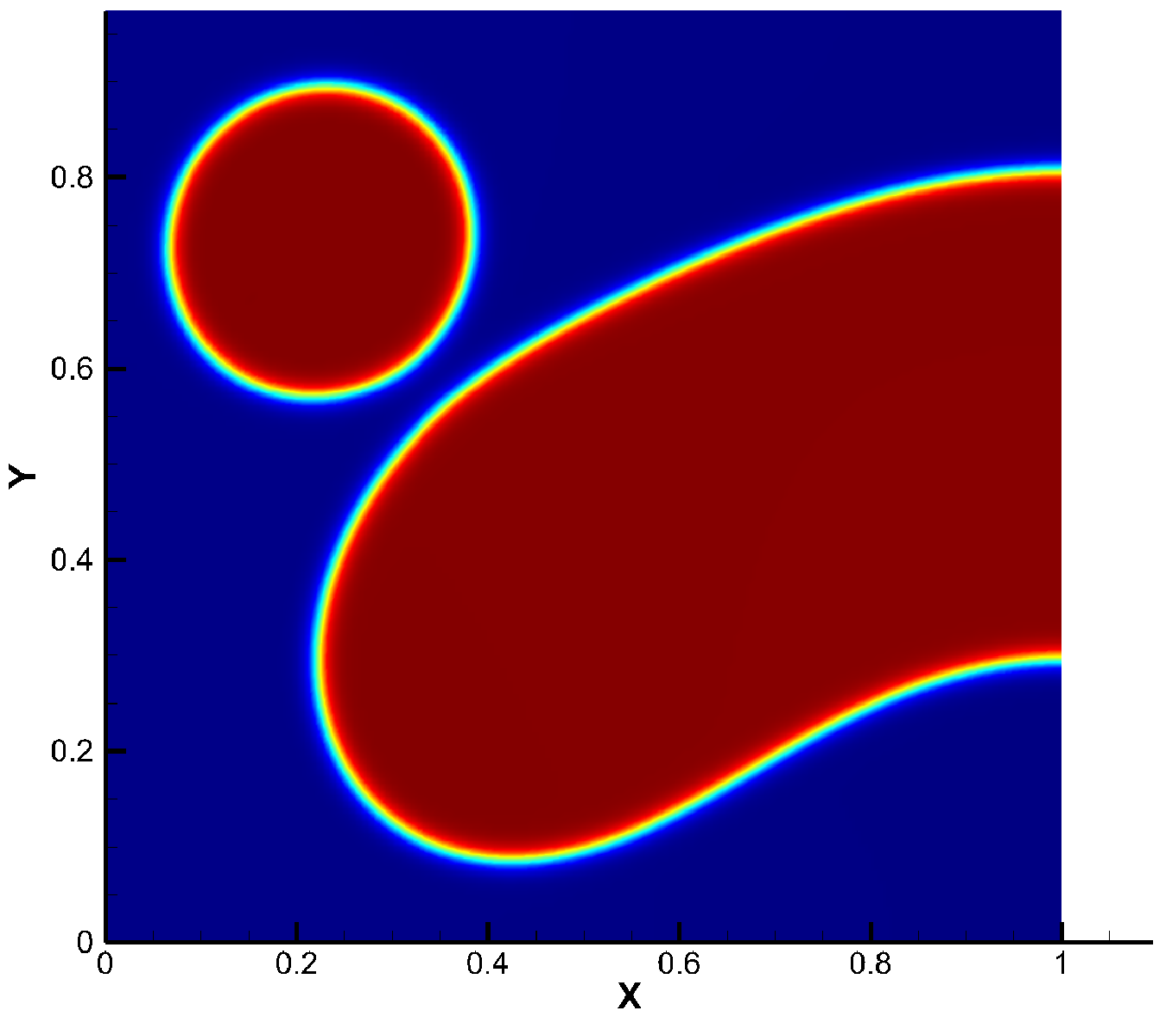}
		\end{minipage}
	}%
   \subfigure[$t=4$]{
		\begin{minipage}[t]{0.2\linewidth}
			\centering
			\includegraphics[width=\textwidth]{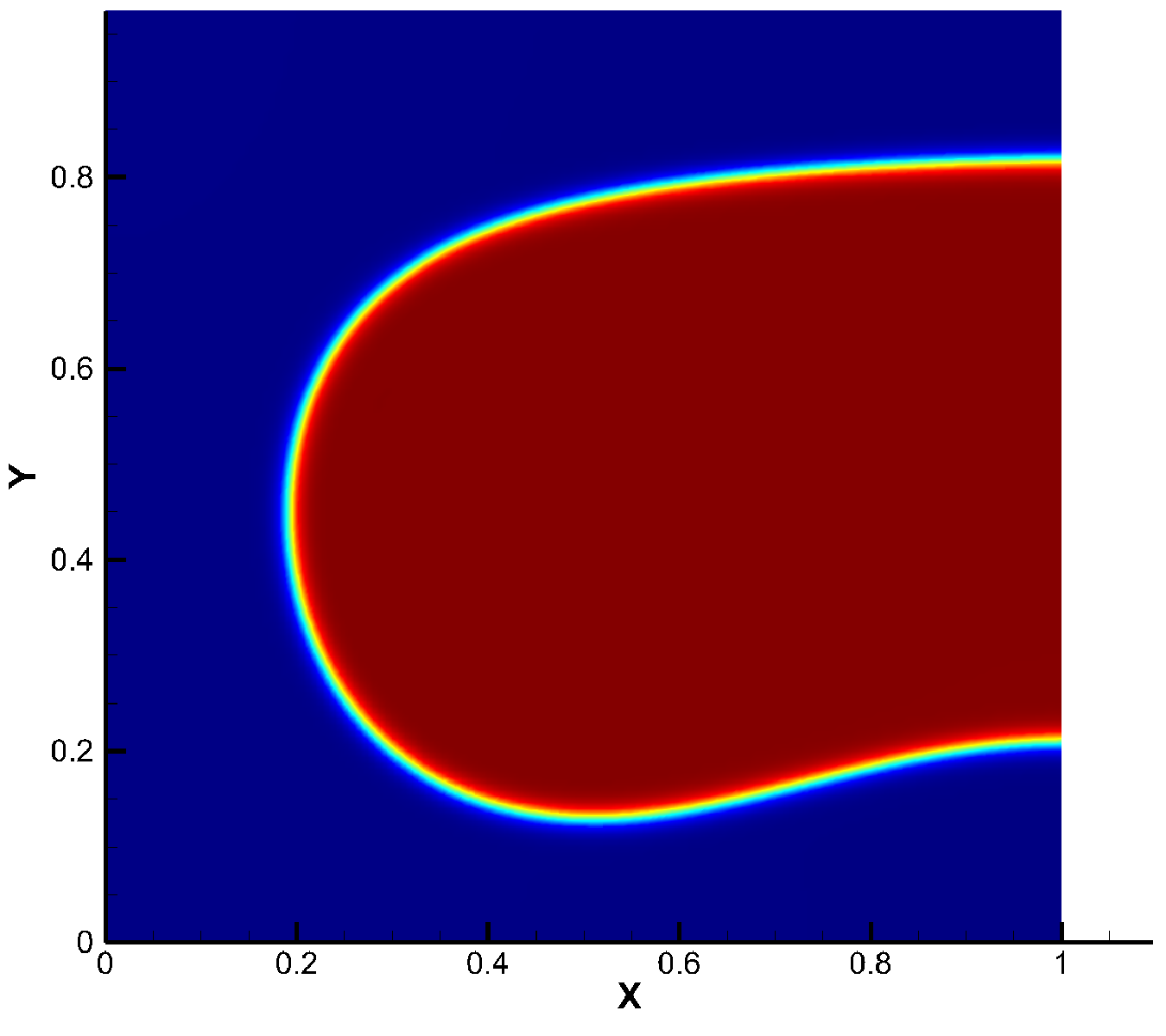}
		\end{minipage}
	}%
	\centering
	\caption{Snapshots of phase field dynamical evolution for spinodal decomposition for case I.}
\label{figure-spi}
\end{figure}

We apply the homogeneous Dirichlet boundary conditions to the velocity and magnetic fields, and enforce the homogeneous Neumann boundary conditions for the phase field and chemical potential. The time step size $\Delta t=1/1000$ and the mesh size $h=1/150$ are selected to investigate the evolution of the phase field for the case I for 2D in Figure \ref{figure-spi}. We find that over time, the phase field gradually coarsens, and the evolution is similar in both cases.

Then we conduct the system energy  (\ref{system-energy}) and the algorithm energy  (\ref{algorithm-energy}).  We fix the mesh size $h=1/64$, and set the time step size $\Delta t=1, 1/10, 1/100$, and $1/1000$ respectively.
The initial values are set according to equations (\ref{fai0}). The parameters are chosen as
\begin{equation*}
\gamma=1/100,\quad M=1,\quad \nu=1,\quad \mu=1,\quad \lambda=1/100,\quad \sigma=1.
\end{equation*}

\begin{figure}[htbp]
	\centering
\subfigure[system energy]{
		\begin{minipage}[t]{0.49\linewidth}
			\centering
			\includegraphics[width=\textwidth]{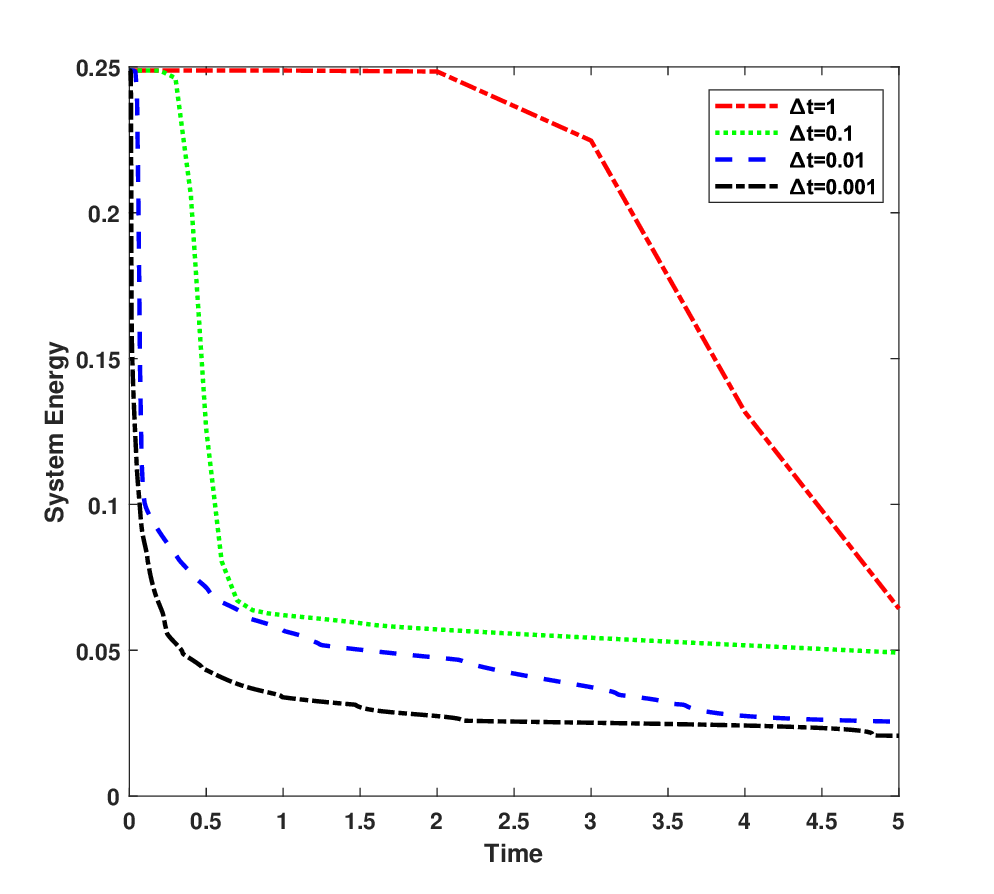}
		\end{minipage}
	}%
	\subfigure[algorithm energy]{
		\begin{minipage}[t]{0.49\linewidth}
			\centering
			\includegraphics[width=\textwidth]{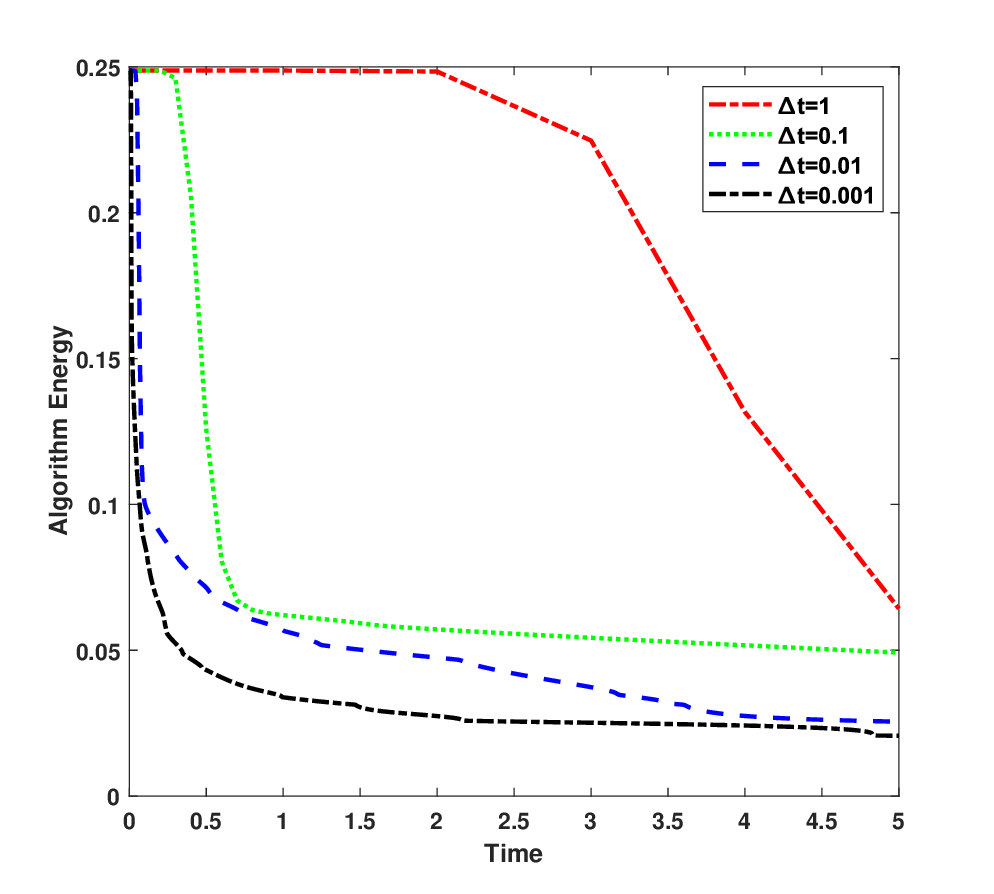}
		\end{minipage}
	}%
	\centering
	\caption{The system energy (left), algorithm energy (right)   for case I.}
\label{figure-energy-mass}
\end{figure}

In Figure \ref{figure-energy-mass} (a) and (b), the comparisons between the system energy and the numerical energy at different time steps are shown for Case I.  As the time step is refined, the energy curves gradually become flat, and the discrete masses are always conserved. This indicates good numerical consistency  in Figure \ref{figure-energy-mass} (a) and (b).  
  
\subsection{Two-phase Kelvin--Helmholtz instability problem}
 
The Kelvin--Helmholtz (K--H) instability  is a common fluid instability caused by the velocity difference at the fluid interface \cite{Lee2015Two, 2023A, Yshin2018vortex}. Because of the gravity and surface tension, the interface for which the lighter fluid is on top of the heavier fluid remains stable. Since the K--H instability has wide applications in natural and industrial fields,  we test the 2D/3D  K--H instability. The domain $\Omega=[0, 1]^{d}$, the parameter values are set to
\begin{equation}\label{K-H-parameter1}
\gamma=1/100,\quad M=1/100,\quad \nu=1/1000,\quad \mu=1,\quad \lambda=1/10000,\quad \sigma=1.
\end{equation}
The boundary conditions for $\B$ at the top ($y=1$) and bottom ($y=0$) are given by ($-1, 0$), and  the vertical component of $\u$ is $u_{2}=0$. The periodic boundary conditions for all variables are applied to the boundaries at $x=0$ and $x=1$ for both single- and double-mode sinusoidal perturbations on the K--H instability problems. 

\subsubsection{Dynamics of single mode sinusoidal perturbation}

This example illustrates the dynamics of a singe mode sinusoidal perturbation at the interface between two fluids. We consider the mesh size $h=1/150$, time step $\Delta t=1/1000$, and the following initial values:
\begin{eqnarray}\label{K-H1}
\left\{
\begin{aligned} 
\phi_{0}&= \tanh(\frac{y-0.5-0.01\sin(2\pi x)}{\sqrt{2}\gamma}),\\
\u_{0}&=\Big( \tanh(\frac{y-0.5-0.01\sin(2\pi x)}{\sqrt{2}\gamma}), 0\Big),\\
\B_{0}&=(1, 0).
\end{aligned}
\right.
\end{eqnarray}
 
\begin{figure}[htbp]
	\centering
\subfigure[$t=0.001$]{
		\begin{minipage}[t]{0.24\linewidth}
			\centering
			\includegraphics[width=\textwidth]{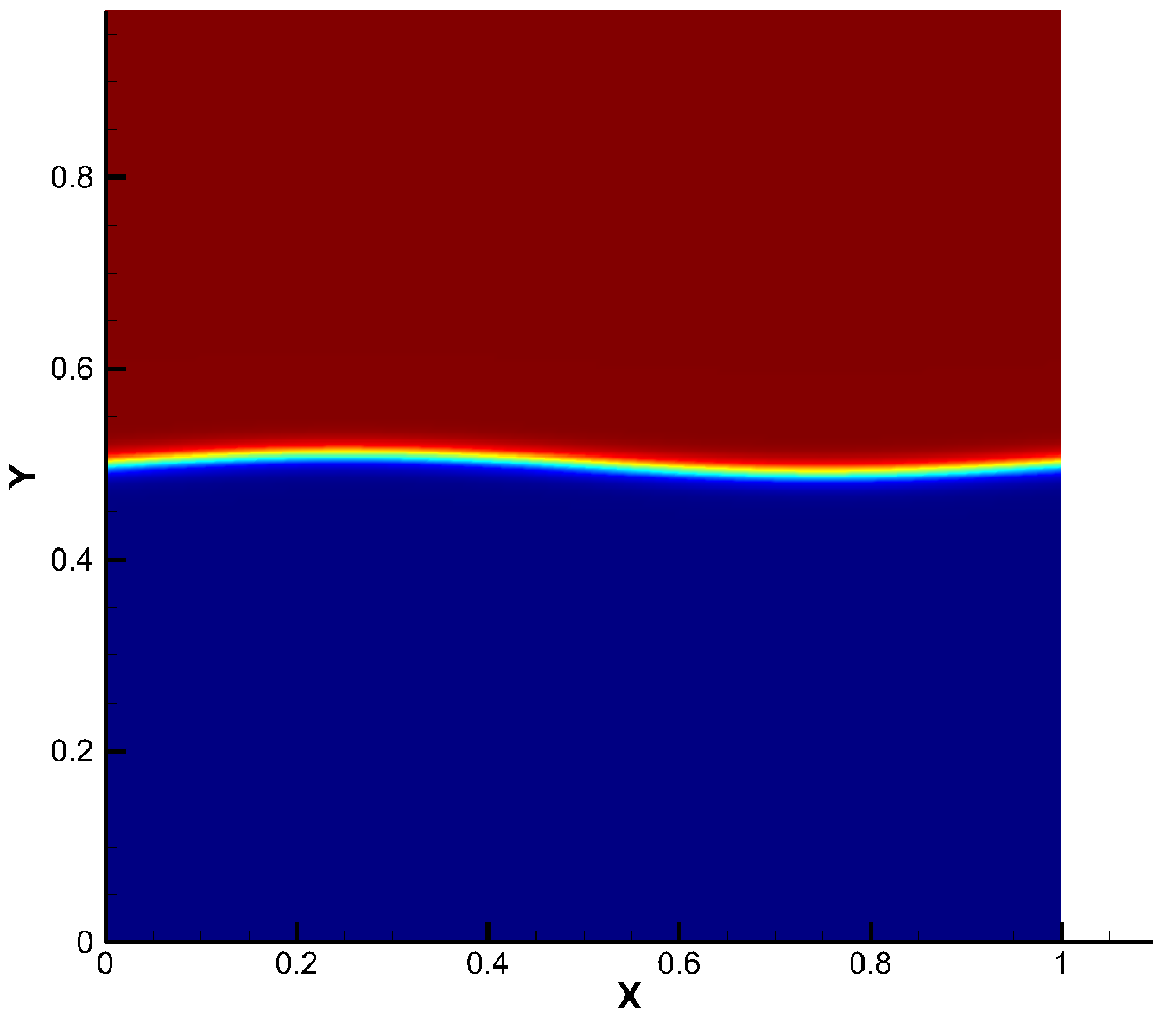}
		\end{minipage}
	}%
	\subfigure[$t=0.6$]{
		\begin{minipage}[t]{0.24\linewidth}
			\centering
			\includegraphics[width=\textwidth]{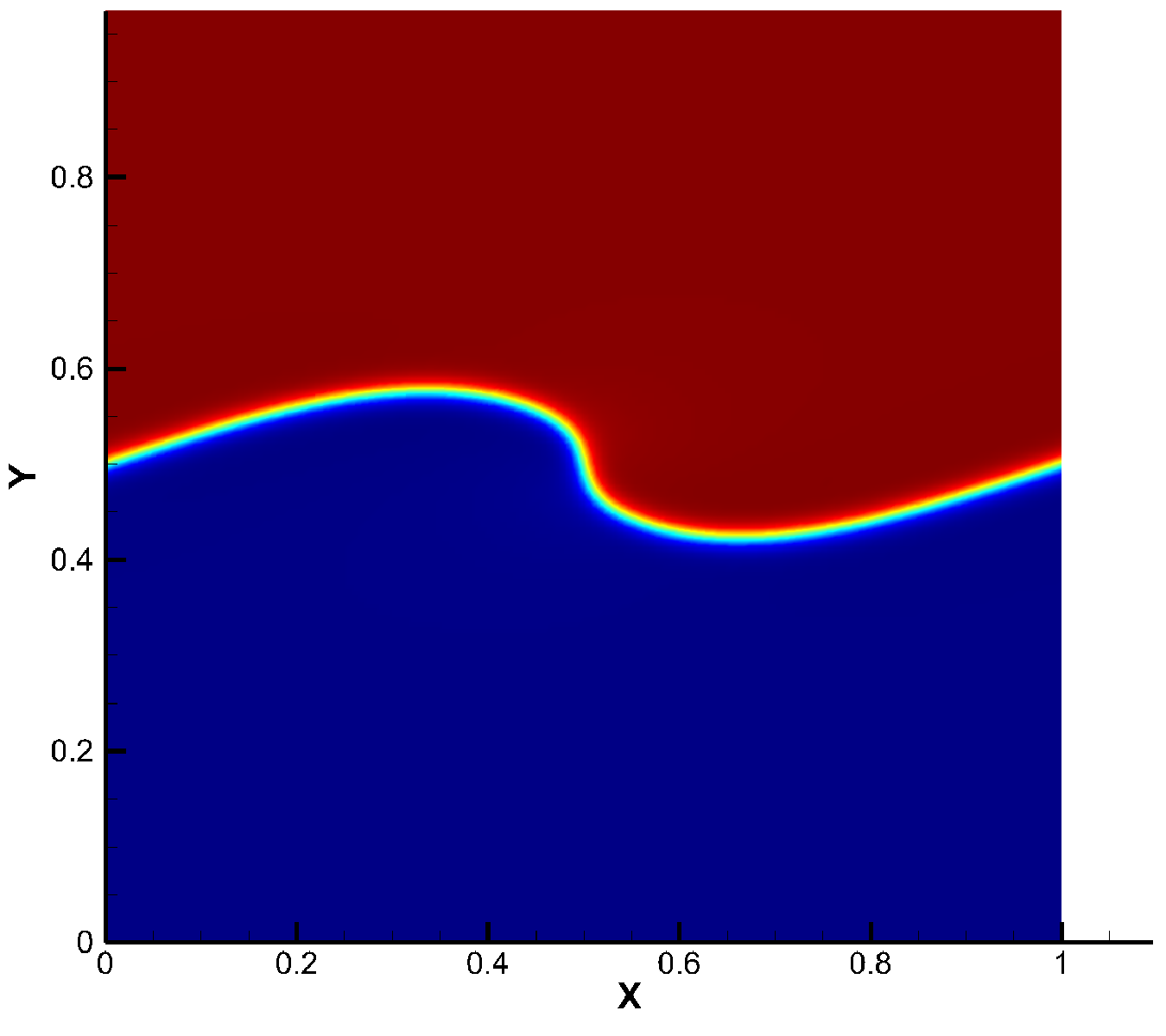}
		\end{minipage}
	}%
	\subfigure[$t=0.85$]{
		\begin{minipage}[t]{0.24\linewidth}
			\centering
			\includegraphics[width=\textwidth]{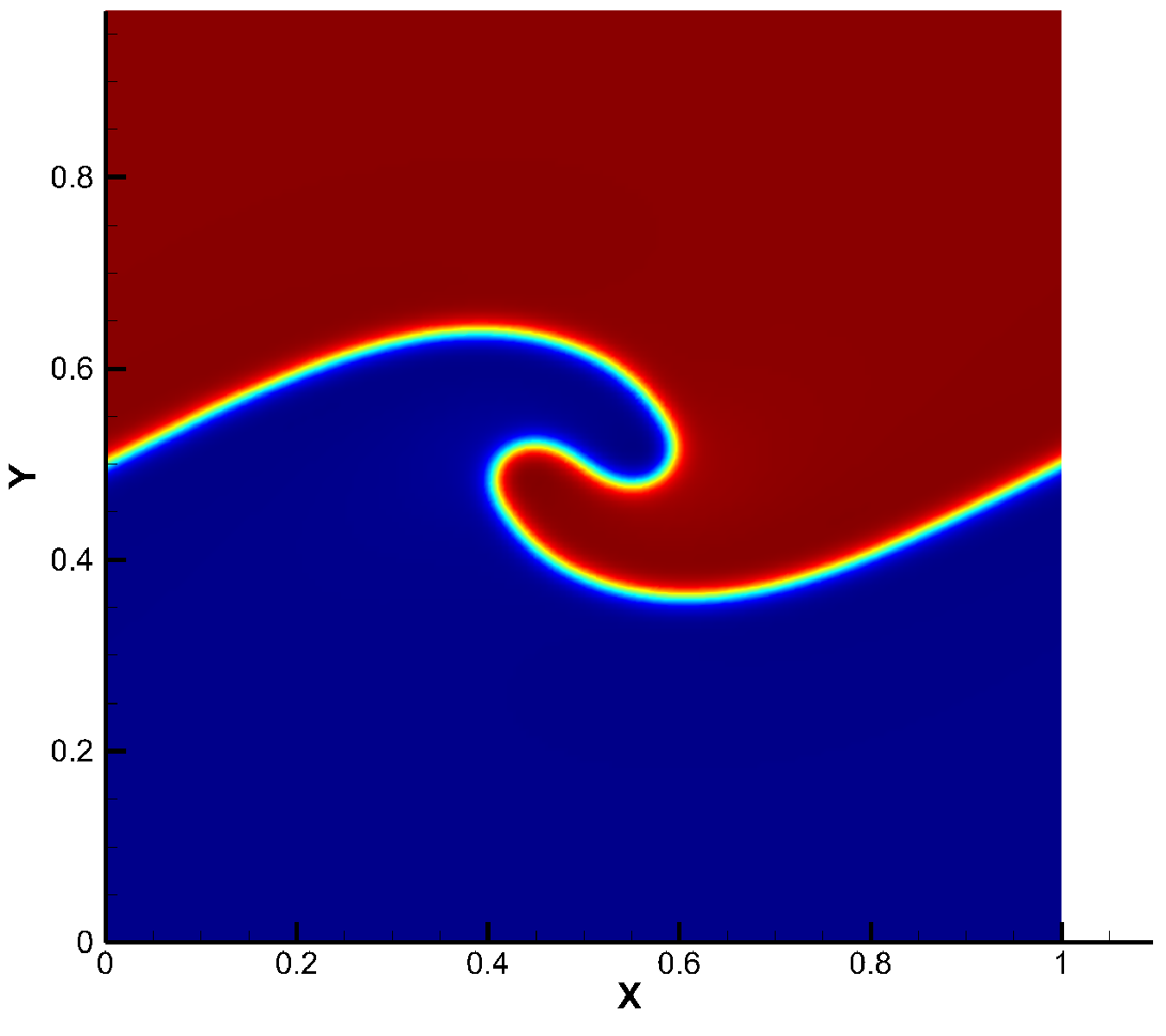}
		\end{minipage}
	}%
	\subfigure[$t=1$]{
		\begin{minipage}[t]{0.24\linewidth}
			\centering
			\includegraphics[width=\textwidth]{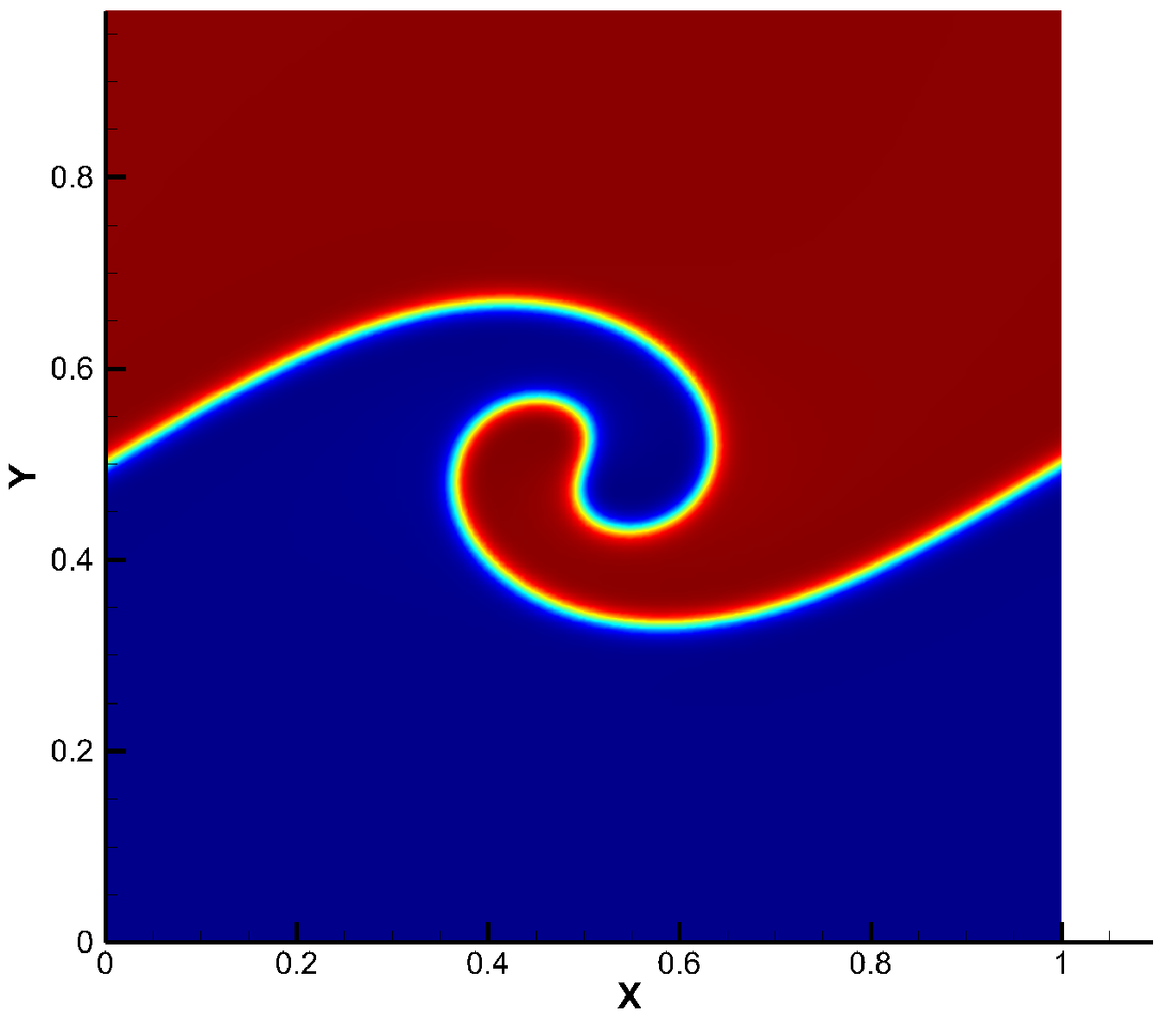}
		\end{minipage}
	}%
	\\
	\subfigure[$t=1.1$]{
		\begin{minipage}[t]{0.24\linewidth}
			\centering
			\includegraphics[width=\textwidth]{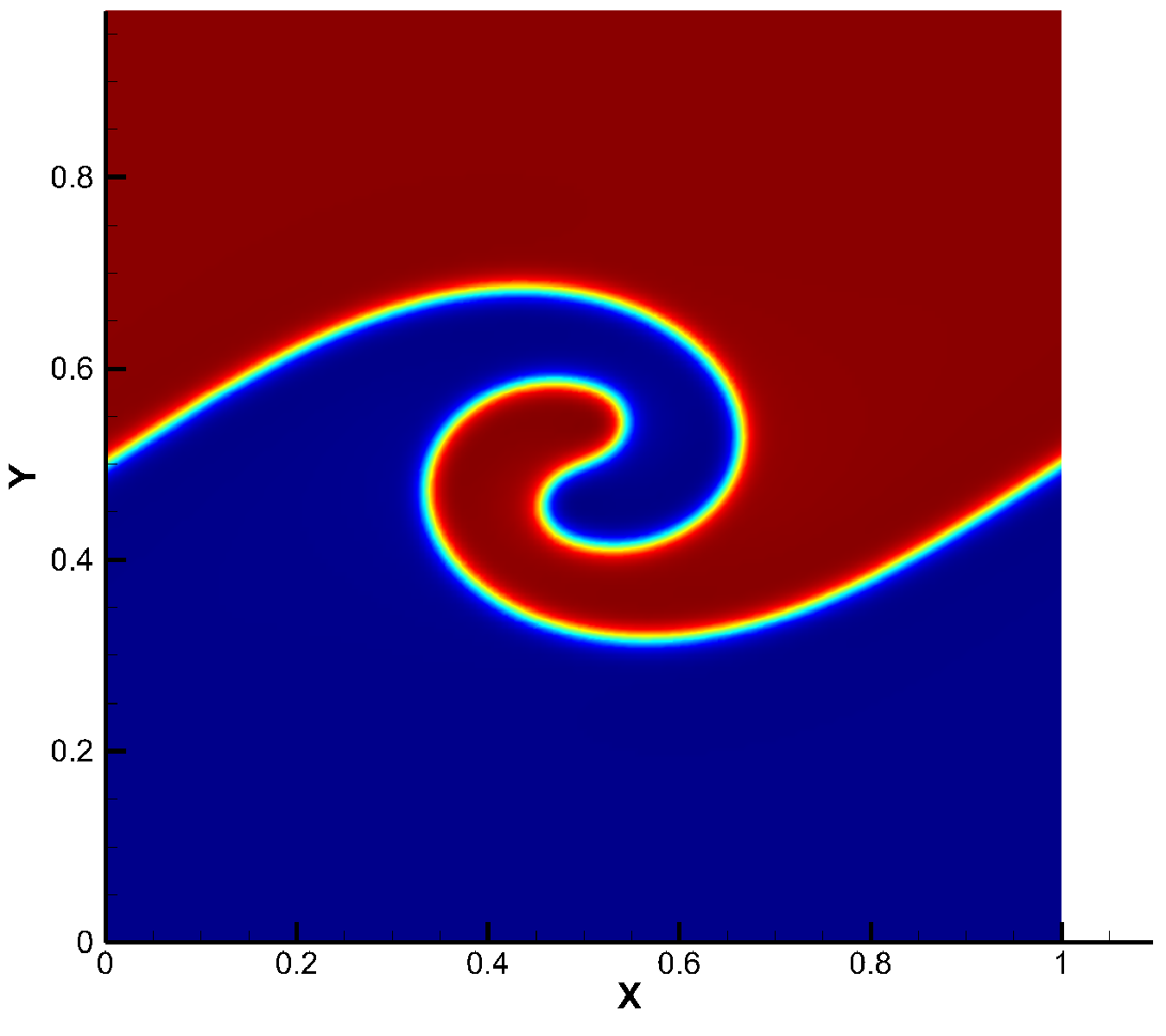}
		\end{minipage}
	}%
\subfigure[$t=1.2$]{
		\begin{minipage}[t]{0.24\linewidth}
			\centering
			\includegraphics[width=\textwidth]{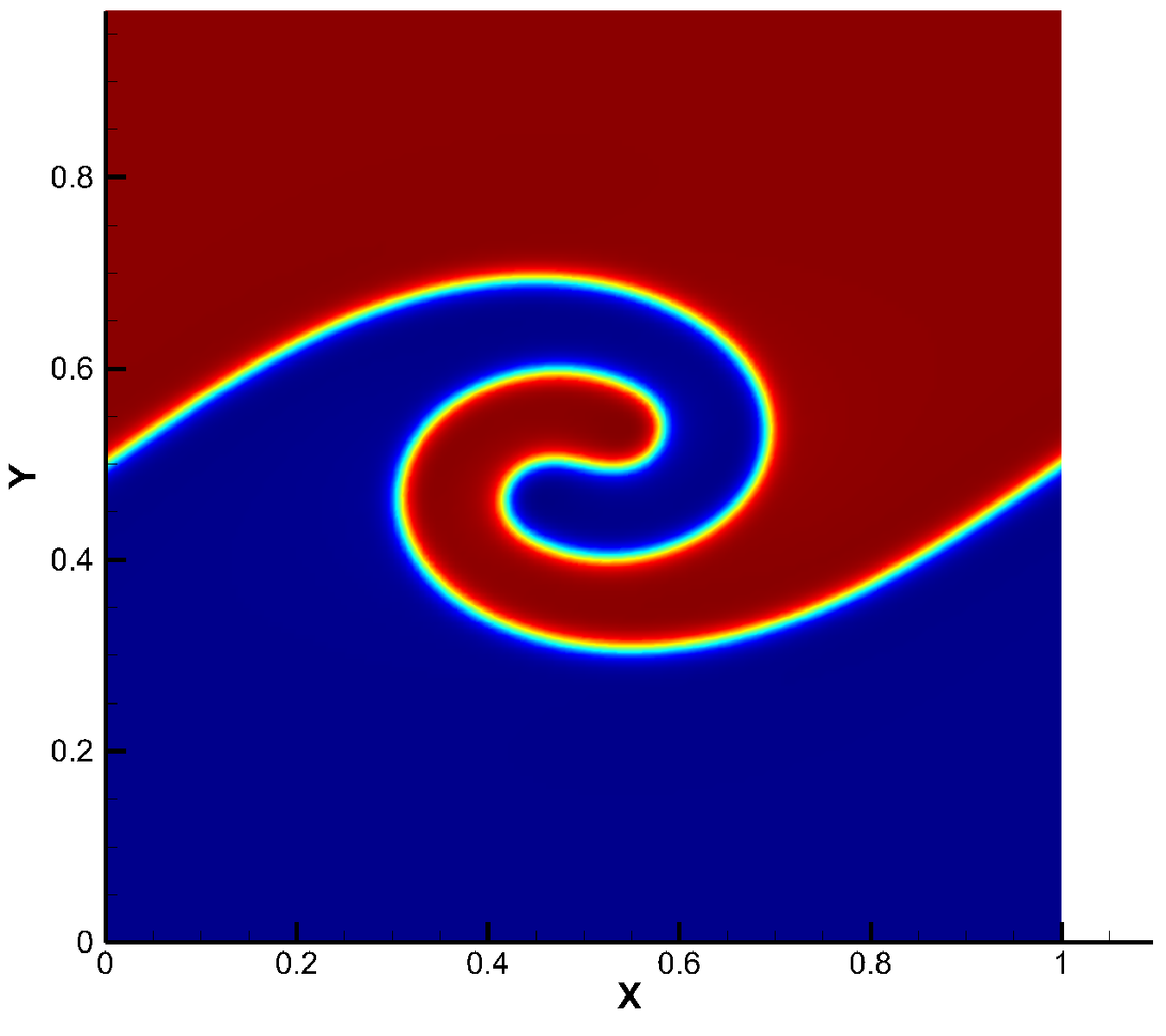}
		\end{minipage}
	}%
   \subfigure[$t=1.4$]{
		\begin{minipage}[t]{0.24\linewidth}
			\centering
			\includegraphics[width=\textwidth]{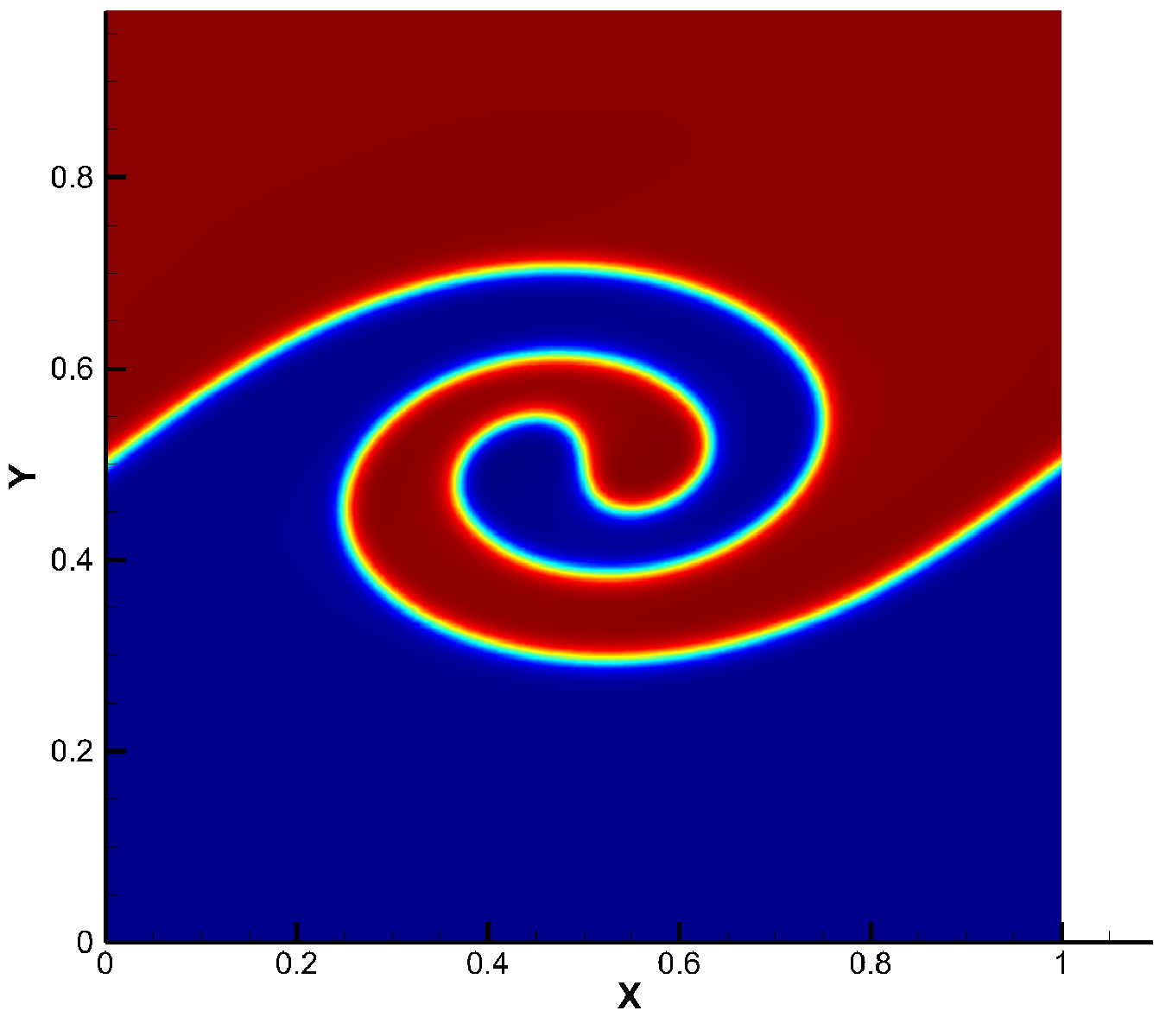}
		\end{minipage}
	}%
\subfigure[$t=1.6$]{
		\begin{minipage}[t]{0.24\linewidth}
			\centering
			\includegraphics[width=\textwidth]{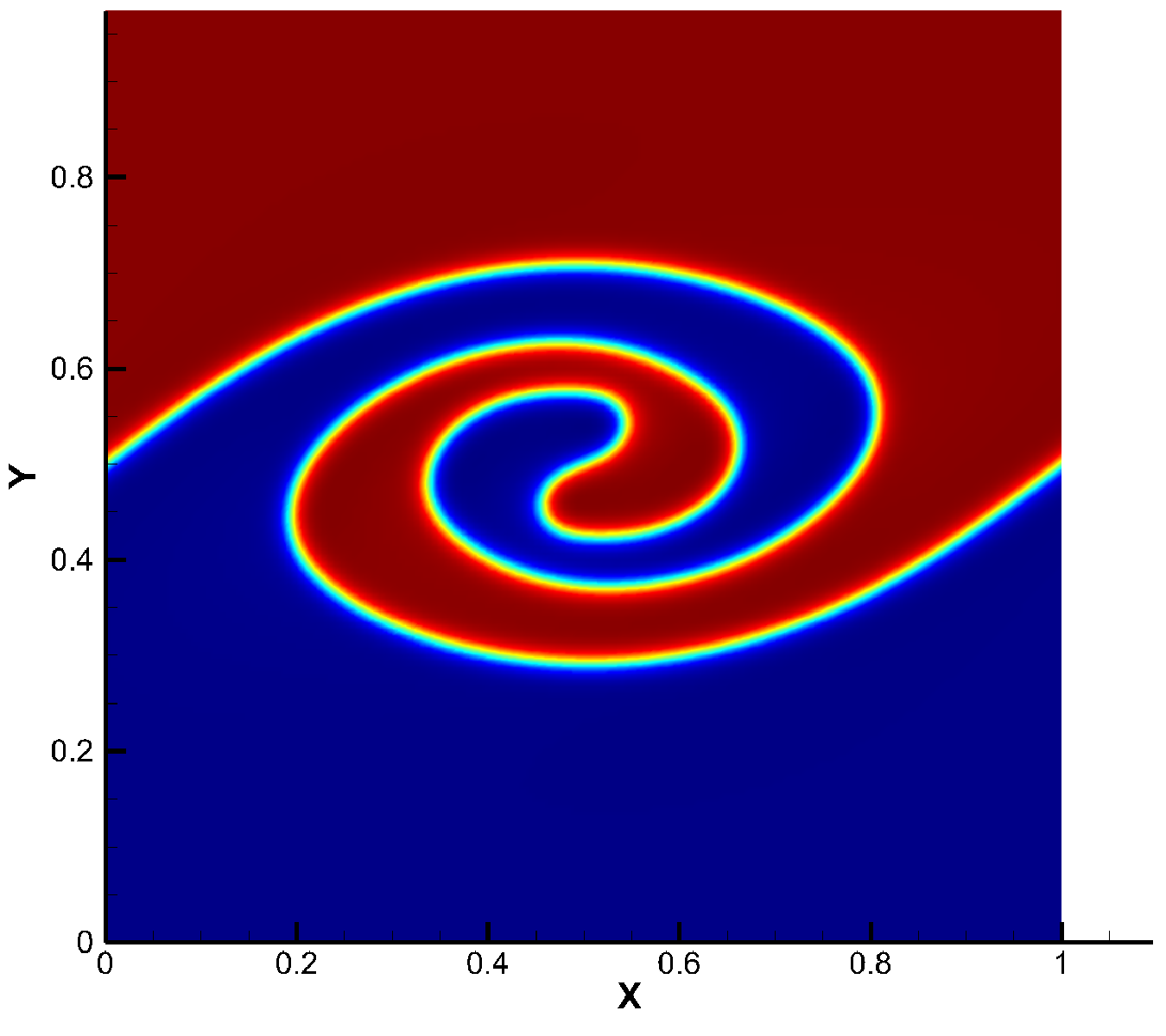}
		\end{minipage}
	}%
	\centering
	\caption{Snapshots of the phase field perturbed sinusoidal at $ t=0.001$ (a), 0.6 (b), 0.85 (c), 1 (d), 1.1 (e), 1.2 (f), 1.4 (g), 1.6 (h) for case I.}
\label{figure-single-phase}
\end{figure}

Figure \ref{figure-single-phase} shows the evolution of the phase field with a single-mode sinusoidal interface perturbation at different times. The interface undergoes a rolling up at the center of the domain at $t=0.6$. The rolling up of the interface forms a spiral shape at a later time, specifically showing  the characteristic features of K--H instability, as depicted in Figure 6.8.

\begin{figure}[htbp]
	\centering
\subfigure[$t=0.001$]{
		\begin{minipage}[t]{0.24\linewidth}
			\centering
			\includegraphics[width=\textwidth]{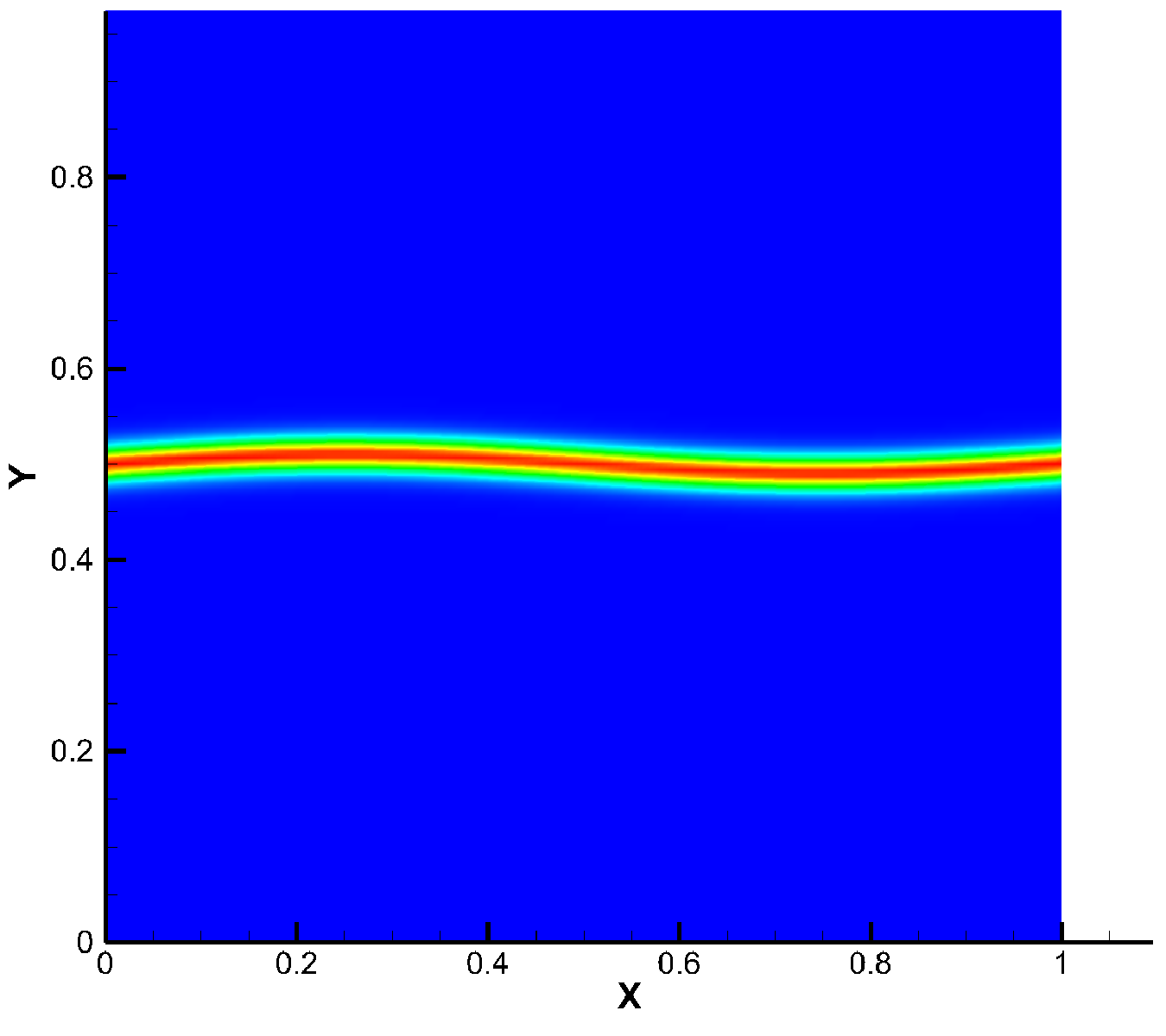}
		\end{minipage}
	}%
	\subfigure[$t=0.6$]{
		\begin{minipage}[t]{0.24\linewidth}
			\centering
			\includegraphics[width=\textwidth]{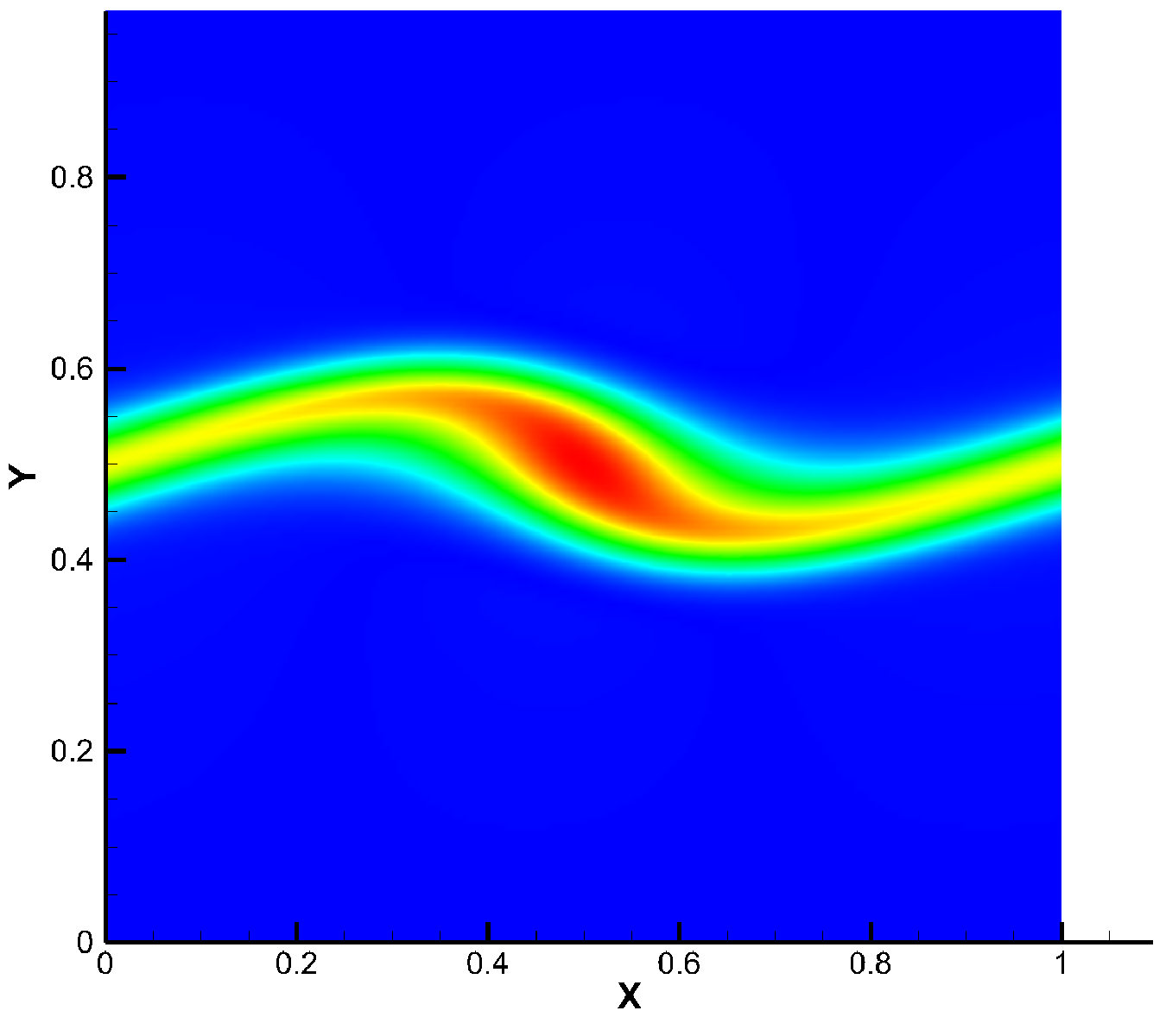}
		\end{minipage}
	}%
	\subfigure[$t=0.85$]{
		\begin{minipage}[t]{0.24\linewidth}
			\centering
			\includegraphics[width=\textwidth]{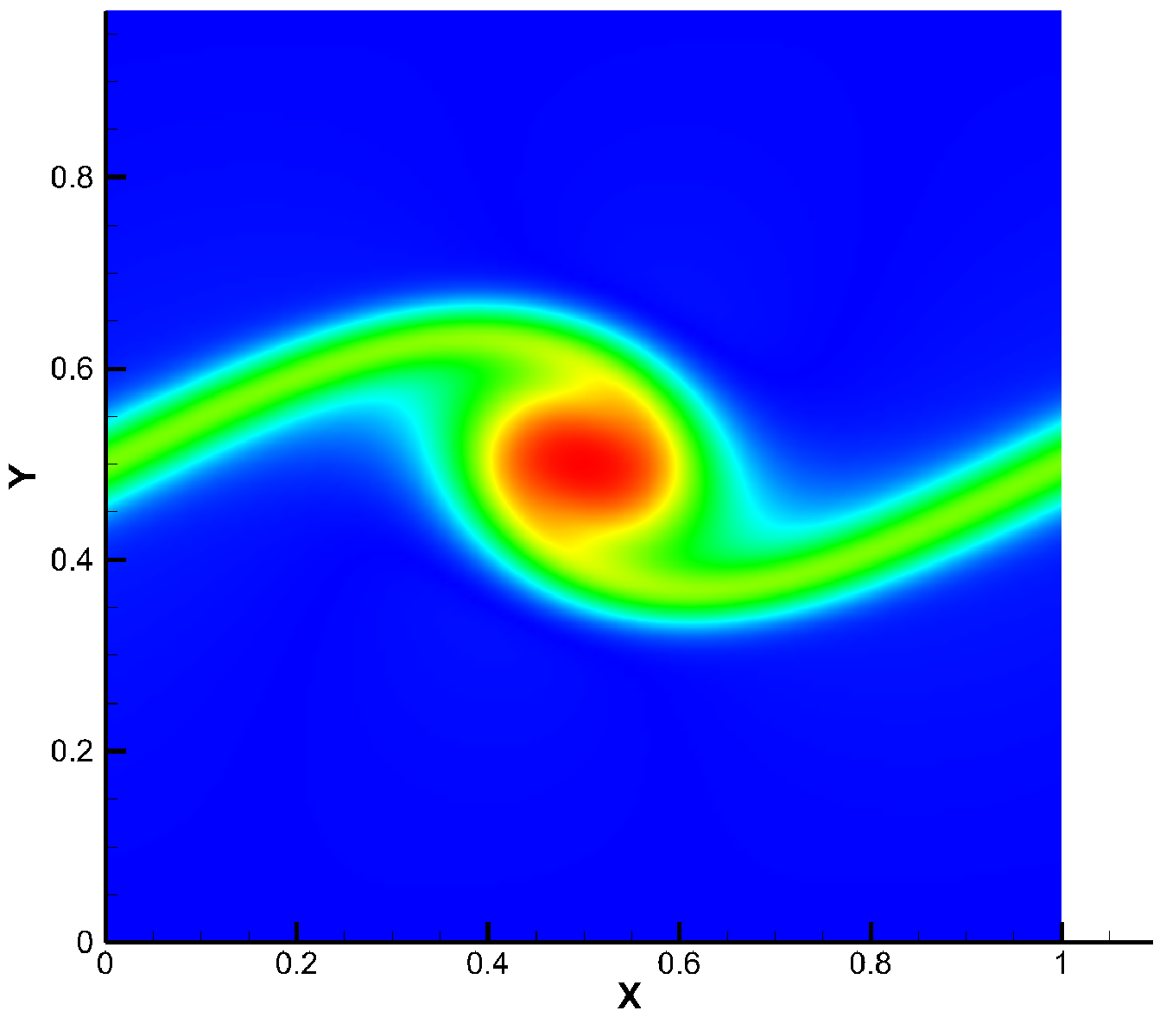}
		\end{minipage}
	}%
\subfigure[$t=1$]{
		\begin{minipage}[t]{0.24\linewidth}
			\centering
			\includegraphics[width=\textwidth]{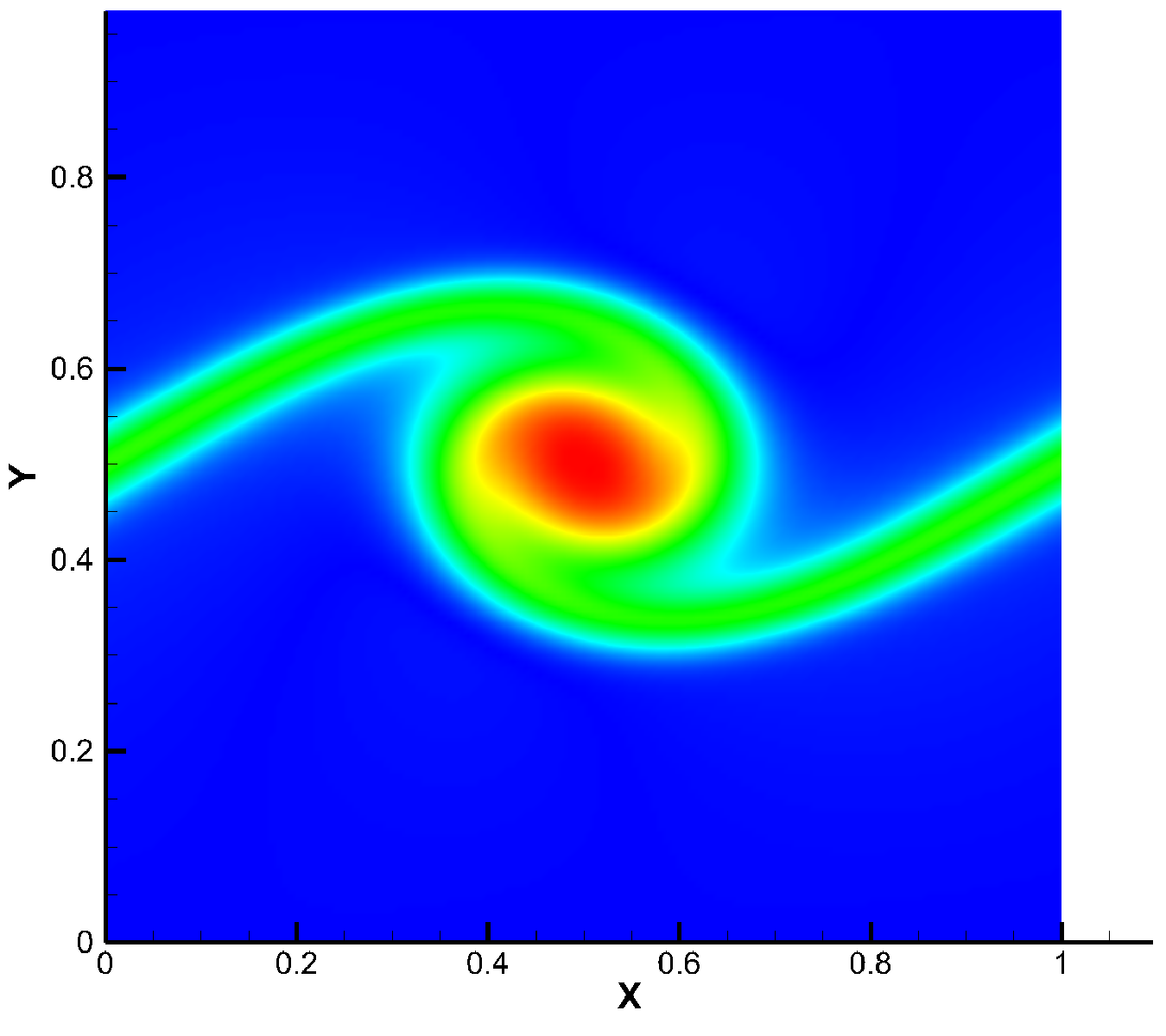}
		\end{minipage}
	}%
	\\
	\subfigure[$t=1.1$]{
		\begin{minipage}[t]{0.24\linewidth}
			\centering
			\includegraphics[width=\textwidth]{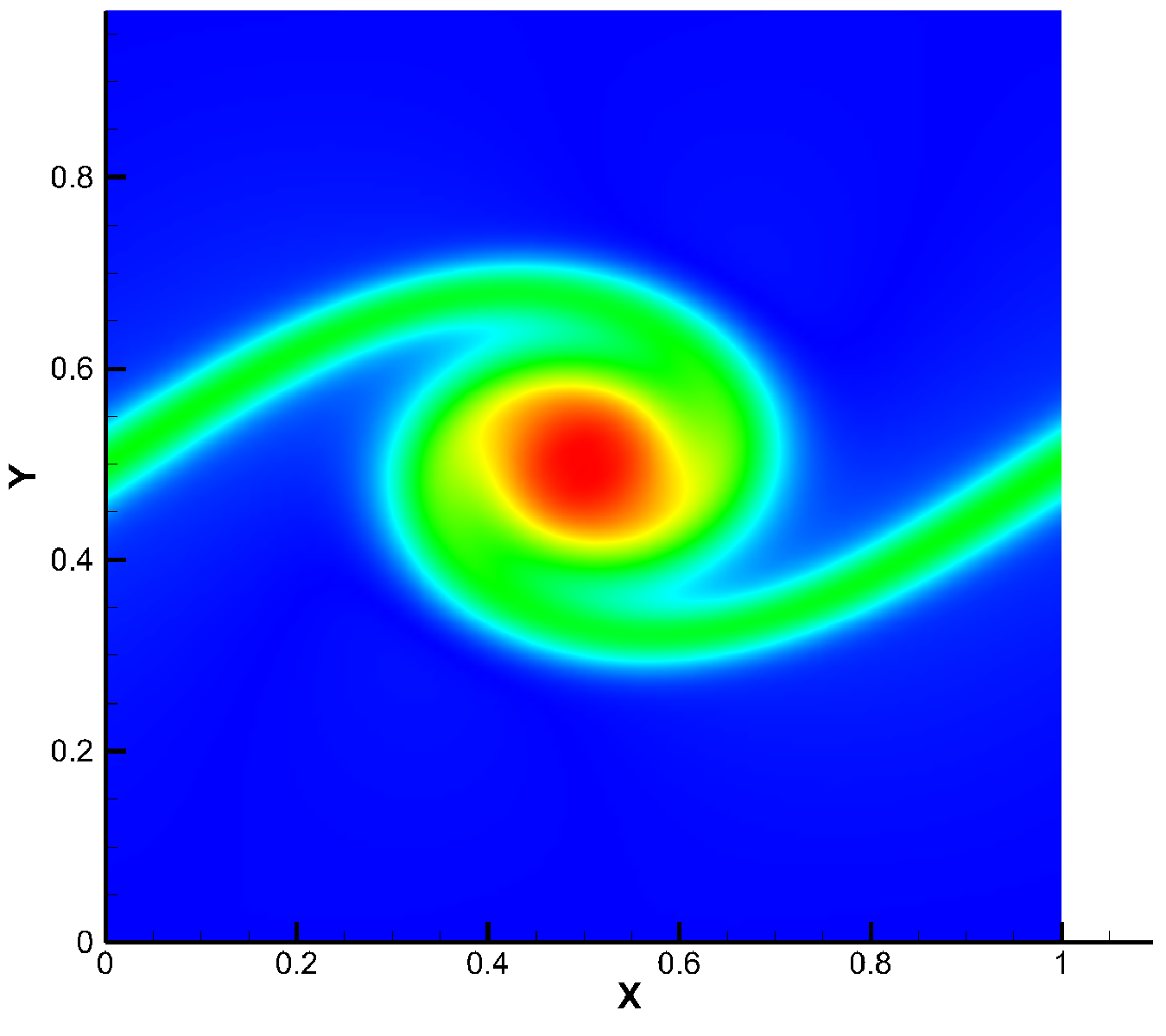}
		\end{minipage}
	}%
\subfigure[$t=1.2$]{
		\begin{minipage}[t]{0.24\linewidth}
			\centering
			\includegraphics[width=\textwidth]{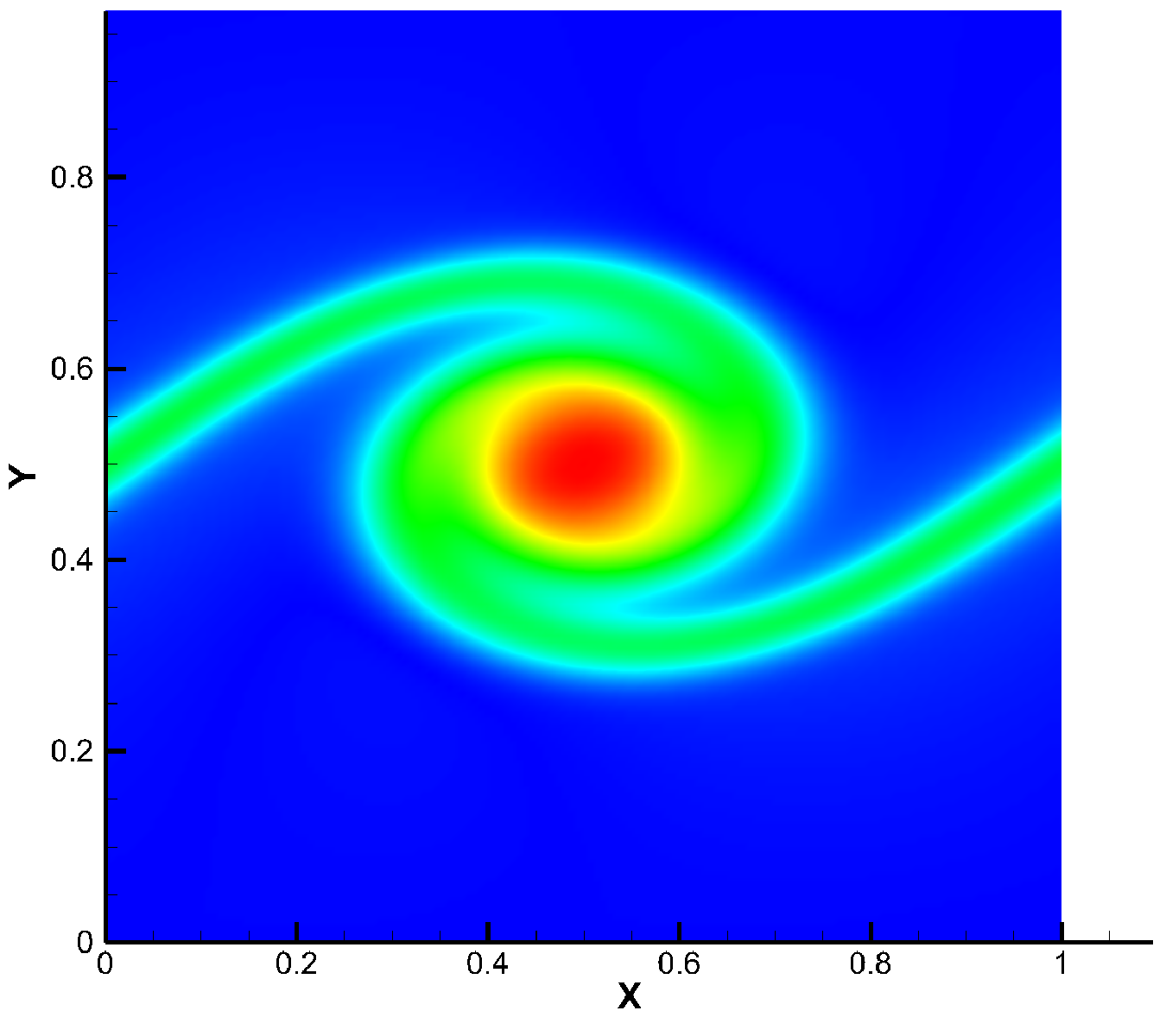}
		\end{minipage}
	}%
   \subfigure[$t=1.4$]{
		\begin{minipage}[t]{0.24\linewidth}
			\centering
			\includegraphics[width=\textwidth]{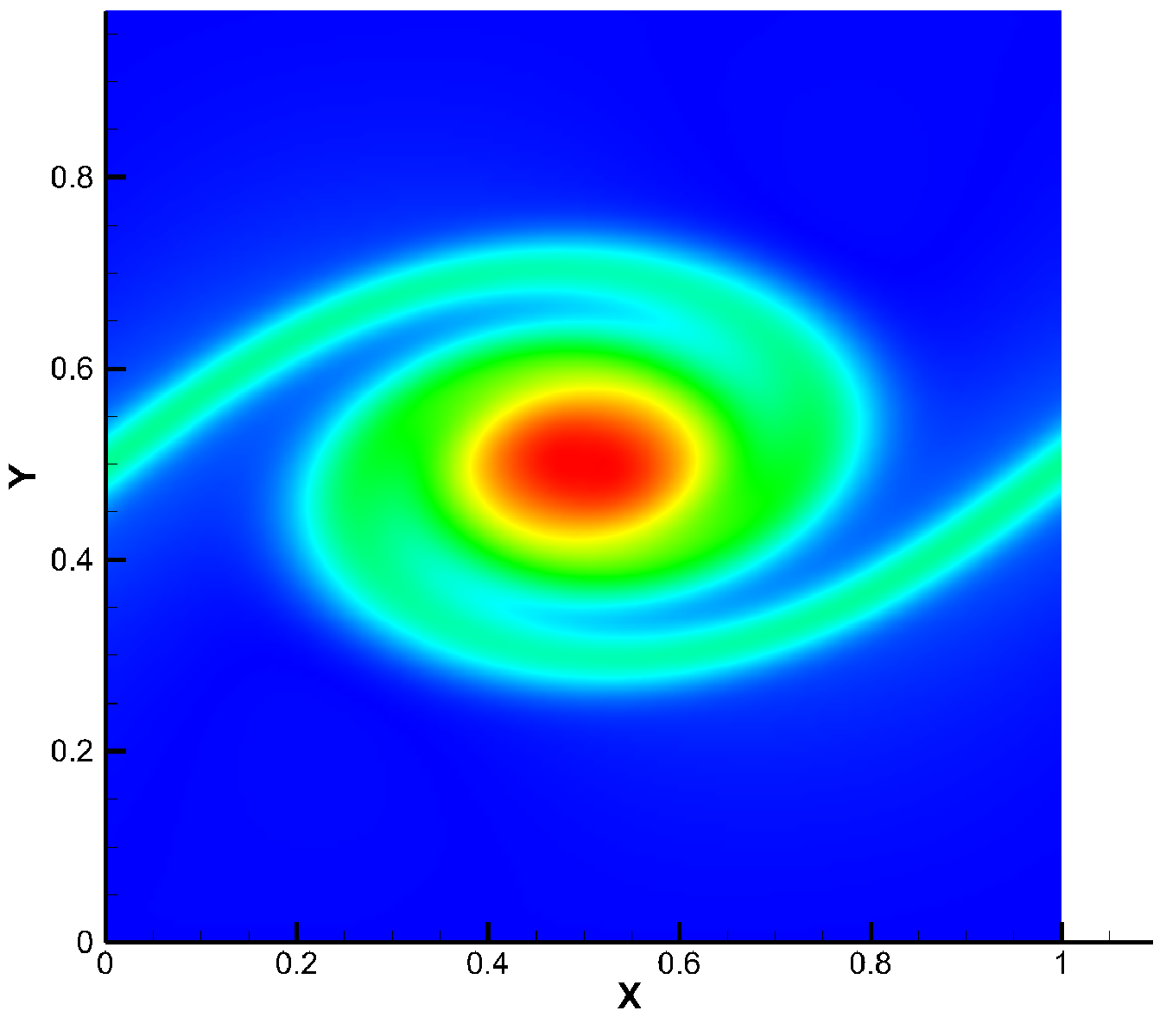}
		\end{minipage}
	}%
\subfigure[$t=1.6$]{
		\begin{minipage}[t]{0.24\linewidth}
			\centering
			\includegraphics[width=\textwidth]{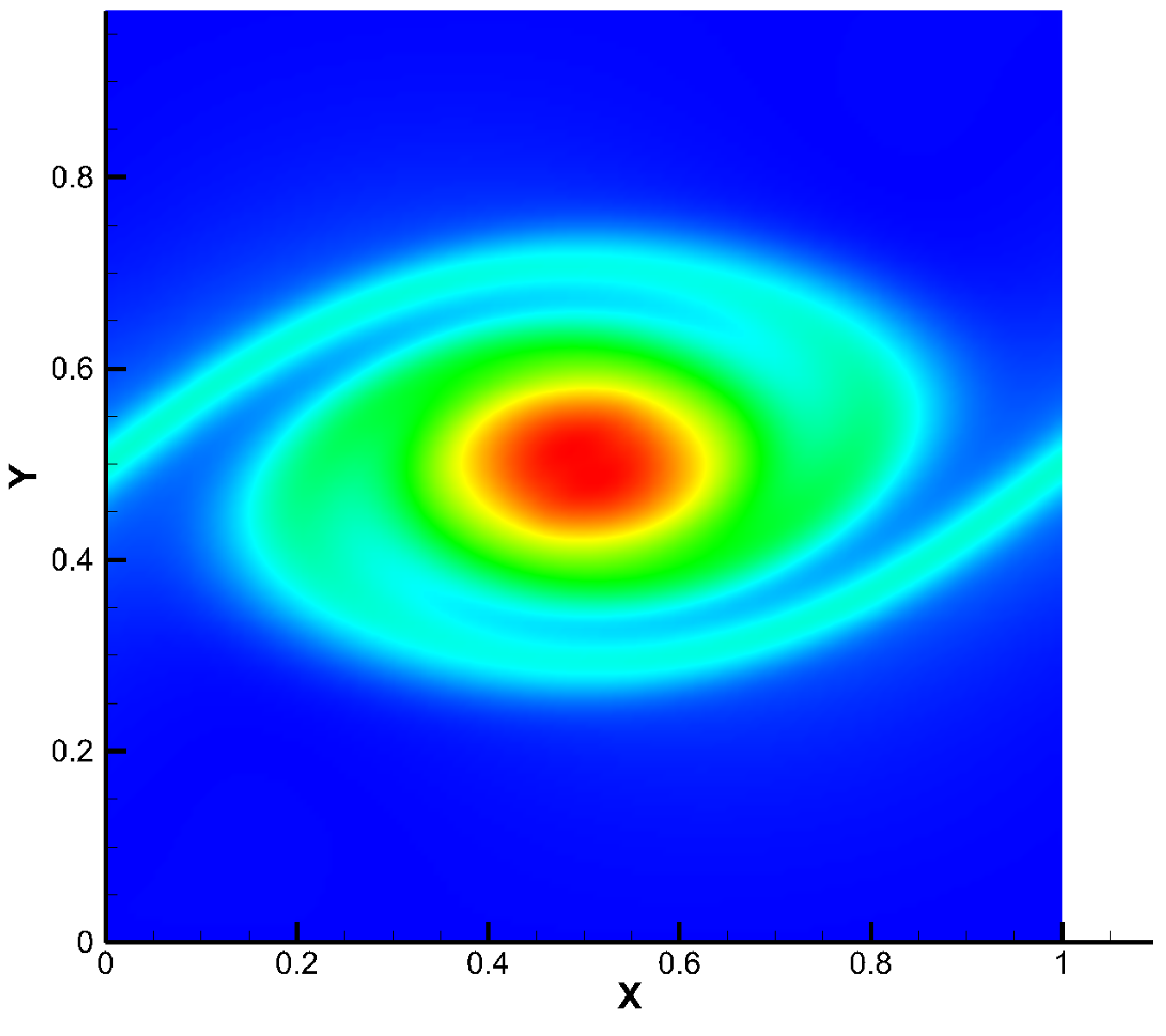}
		\end{minipage}
	}%
	\centering
	\caption{Snapshots of the  vorticity dynamics at time at $ t=0.001$ (a), 0.6 (b), 0.85 (c), 1 (d), 1.1 (e), 1.2 (f), 1.4 (g), 1.6 (h) for case I.}
\label{figure-single-vorticity}
\end{figure}

The snapshots of vorticity evolution are plotted in Figure \ref{figure-single-vorticity}. The fluids at the top and bottom flow in opposite directions, causing the vorticity to migrate towards the center of the region. As the vorticity  accumulates at the center, the interface starts to become more pronounced, and the amplitude of the instability increases. A roll-up phenomenon occurs, transforming the interface into a spiral that takes on a distinctive  ``cat's eye'' configuration.

\subsubsection{Dynamics of double mode sinusoidal perturbation}
\begin{figure}[htbp]
	\centering
\subfigure[$t=0.001$]{
		\begin{minipage}[t]{0.24\linewidth}
			\centering
			\includegraphics[width=\textwidth]{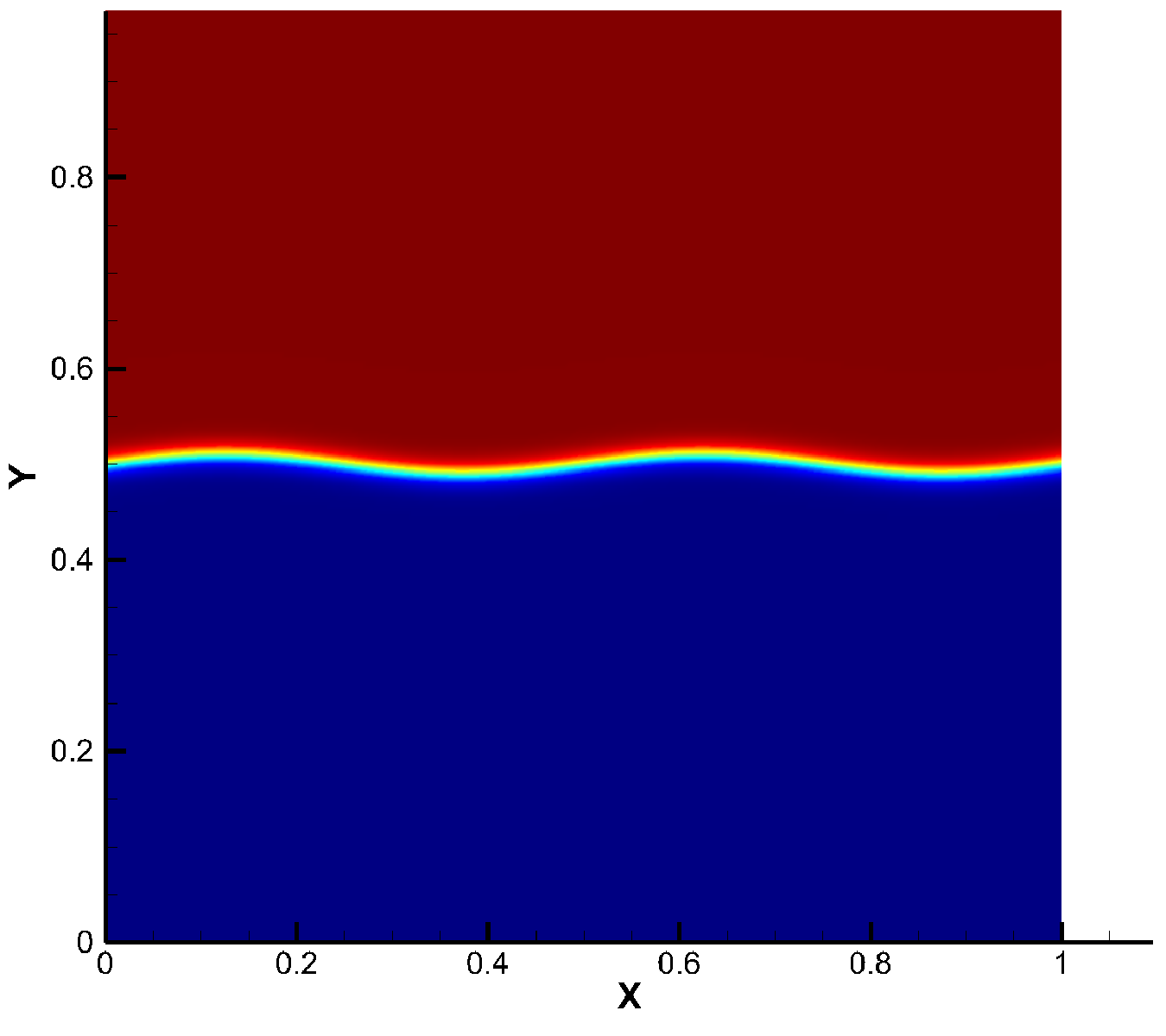}
		\end{minipage}
	}%
	\subfigure[$t=0.3$]{
		\begin{minipage}[t]{0.24\linewidth}
			\centering
			\includegraphics[width=\textwidth]{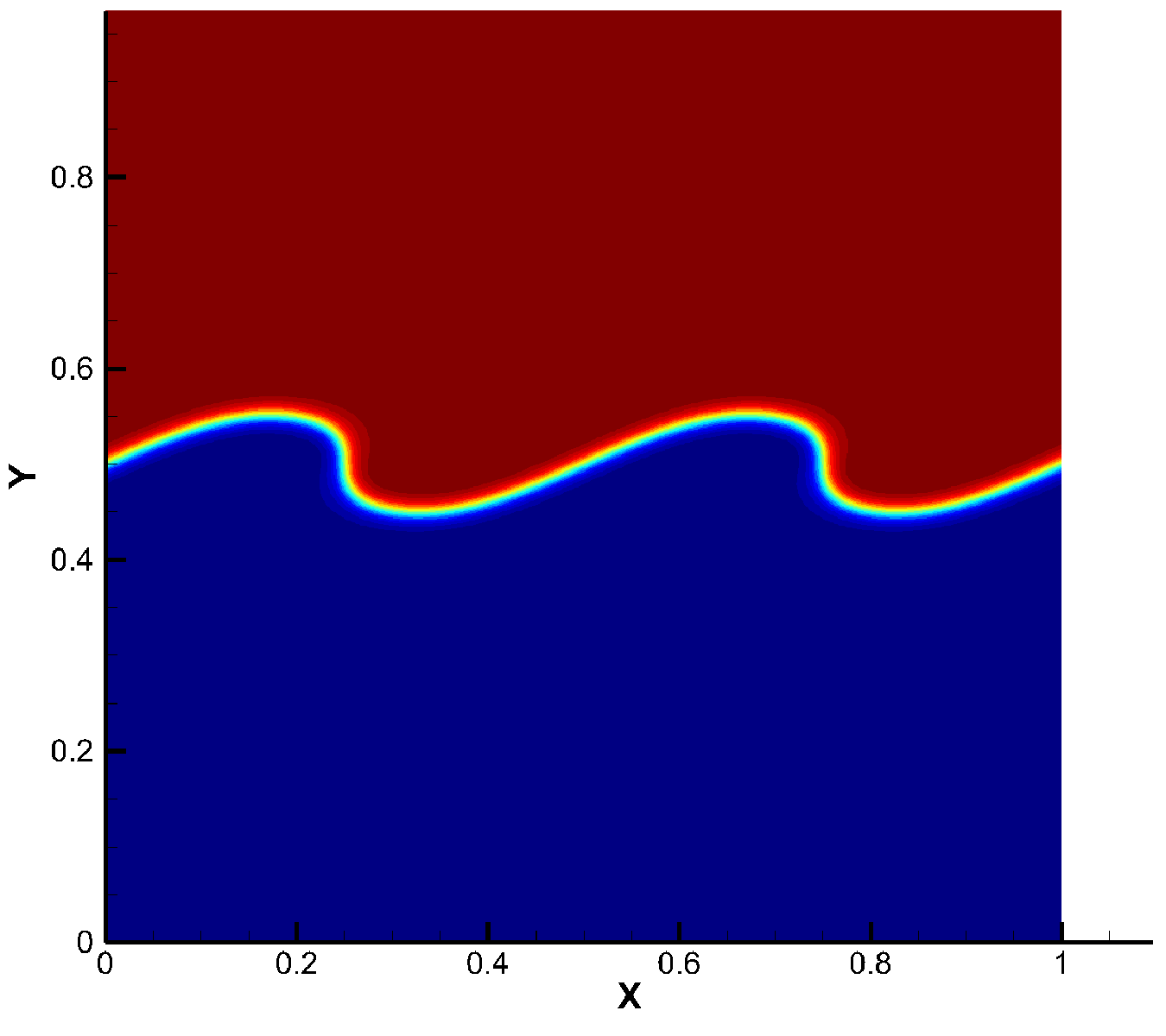}
		\end{minipage}
	}%
   \subfigure[$t=0.5$]{
		\begin{minipage}[t]{0.24\linewidth}
			\centering
			\includegraphics[width=\textwidth]{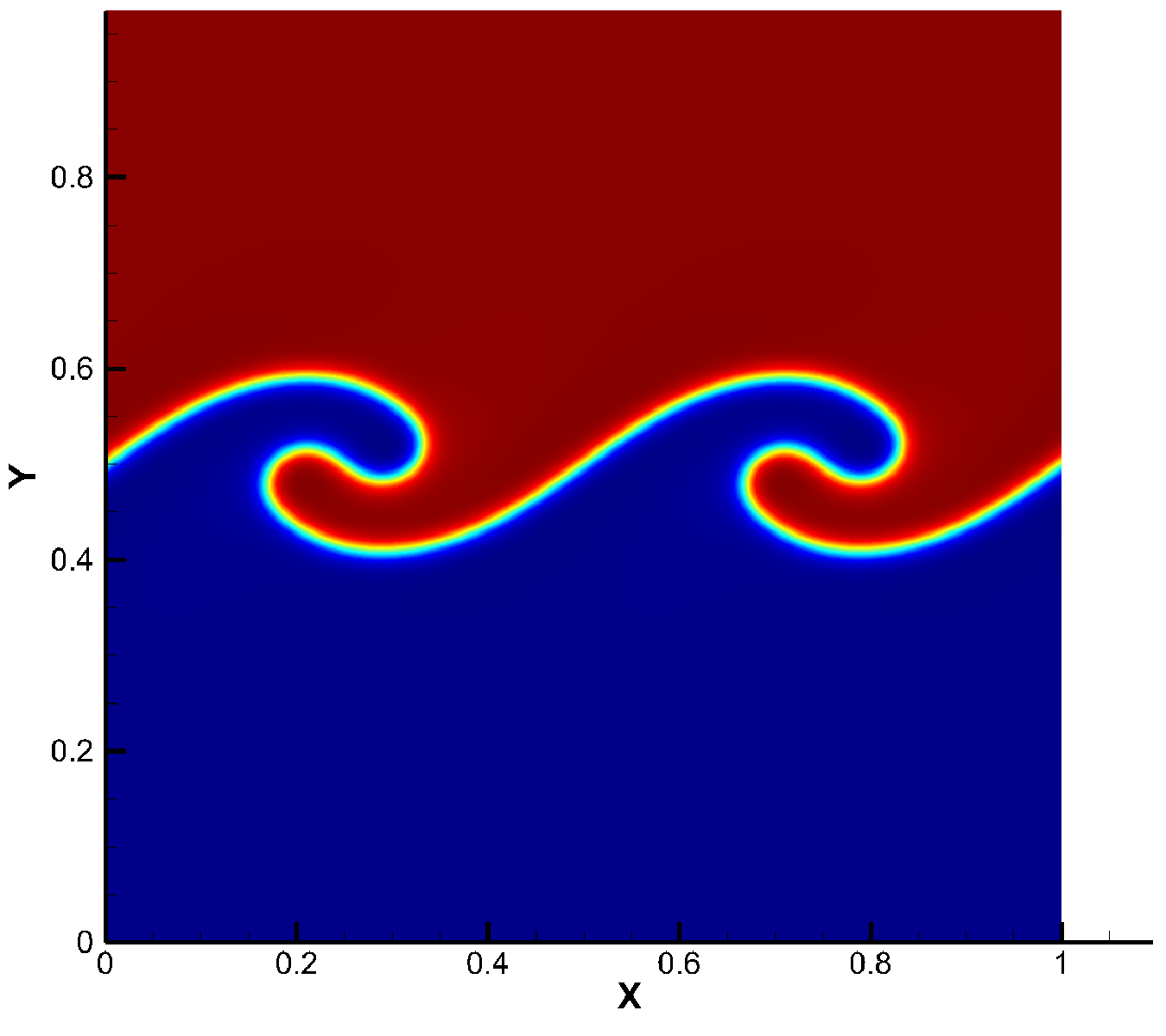}
		\end{minipage}
	}%
\subfigure[$t=0.75$]{
		\begin{minipage}[t]{0.24\linewidth}
			\centering
			\includegraphics[width=\textwidth]{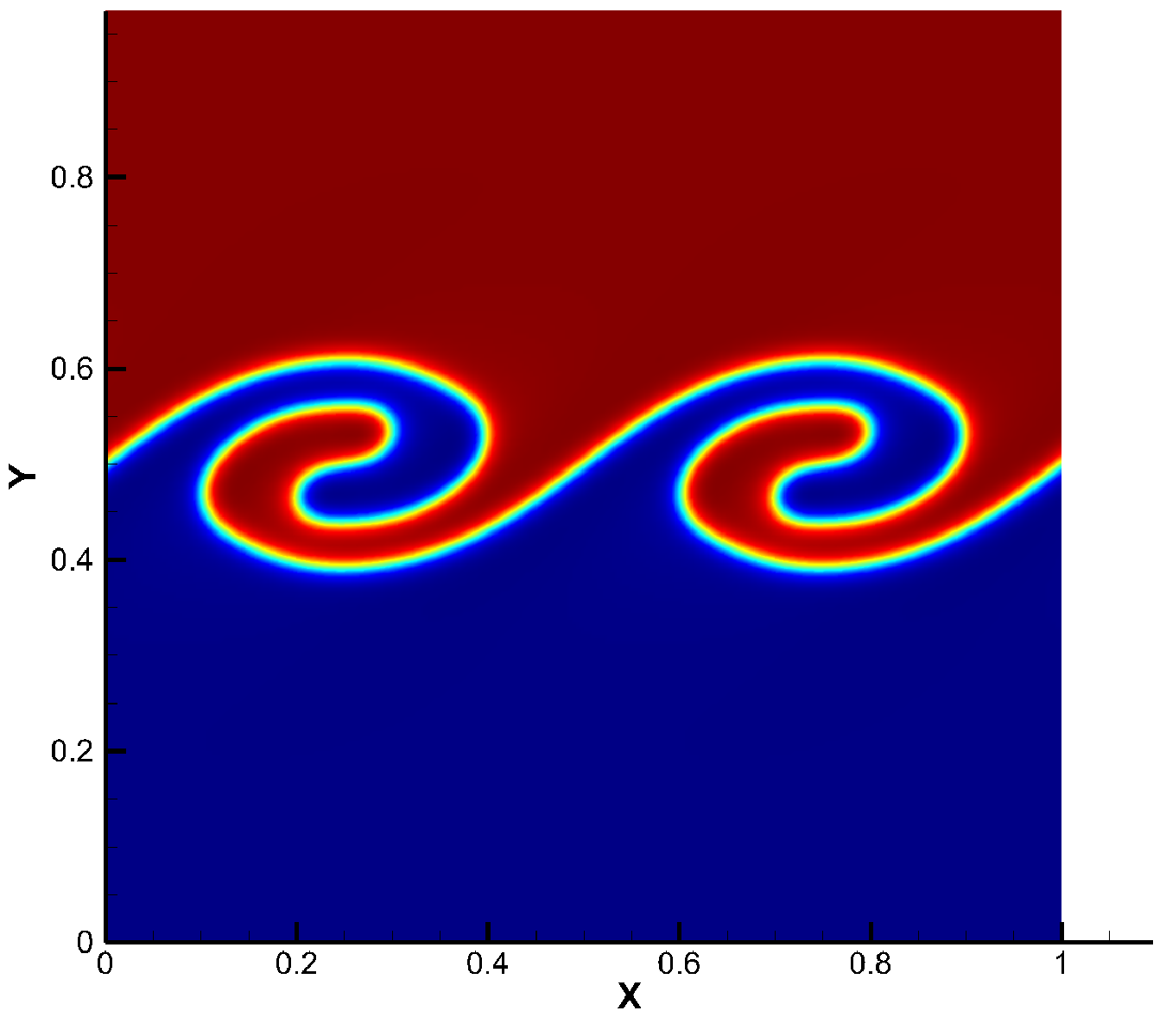}
		\end{minipage}
	}%
	\\
\subfigure[$t=0.001$]{
		\begin{minipage}[t]{0.24\linewidth}
			\centering
			\includegraphics[width=\textwidth]{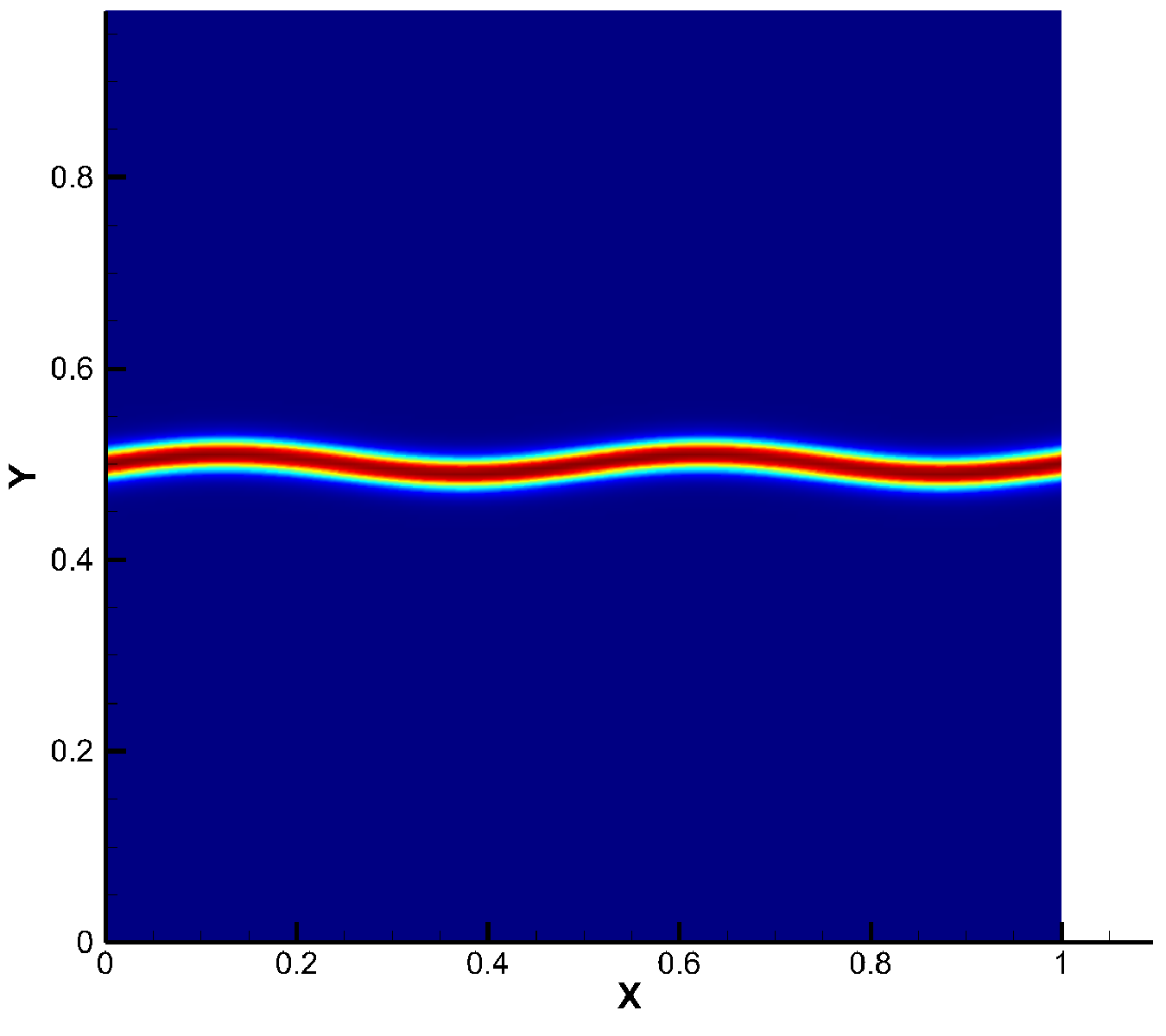}
		\end{minipage}
	}%
	\subfigure[$t=0.3$]{
		\begin{minipage}[t]{0.24\linewidth}
			\centering
			\includegraphics[width=\textwidth]{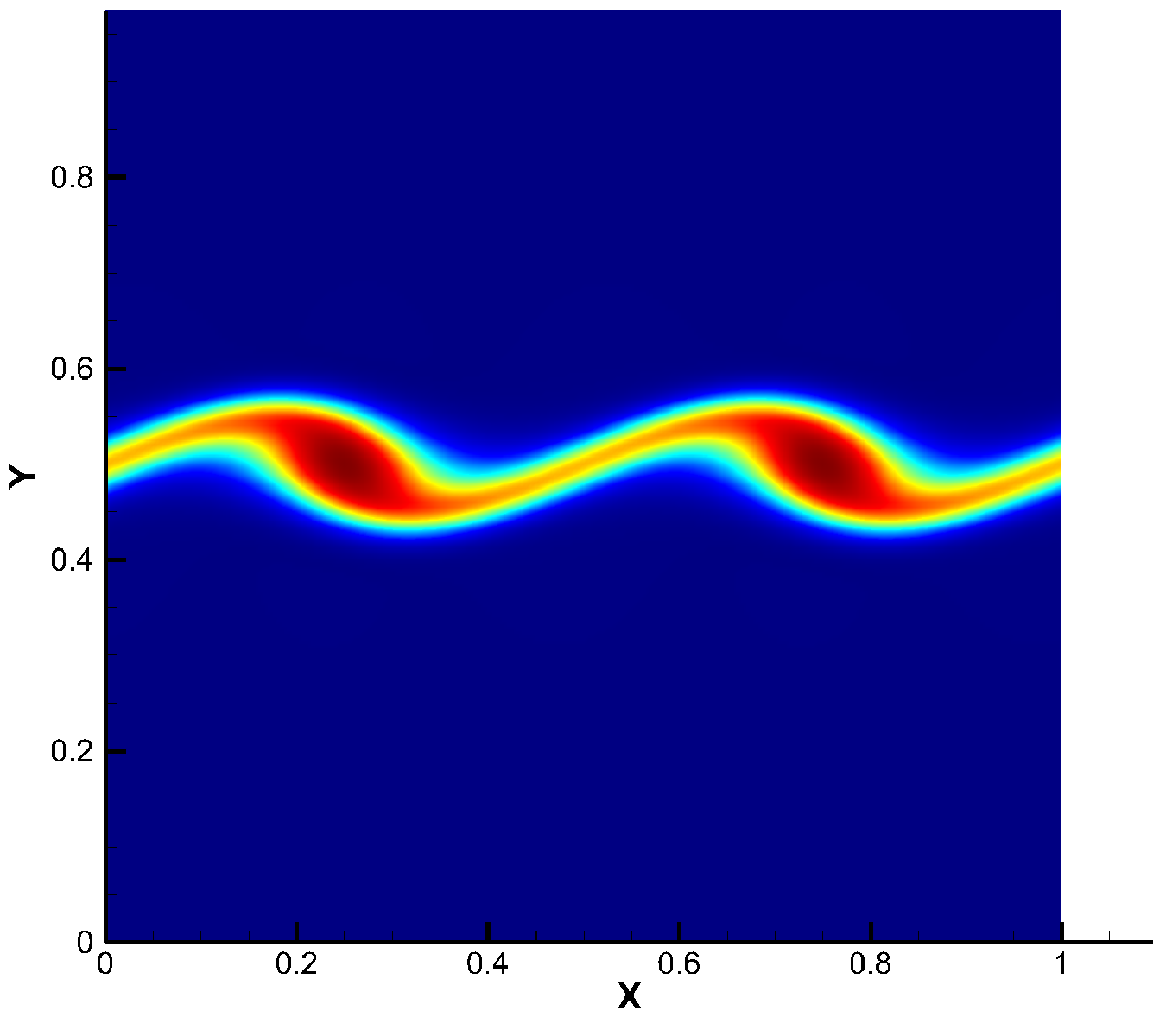}
		\end{minipage}
	}%
   \subfigure[$t=0.5$]{
		\begin{minipage}[t]{0.24\linewidth}
			\centering
			\includegraphics[width=\textwidth]{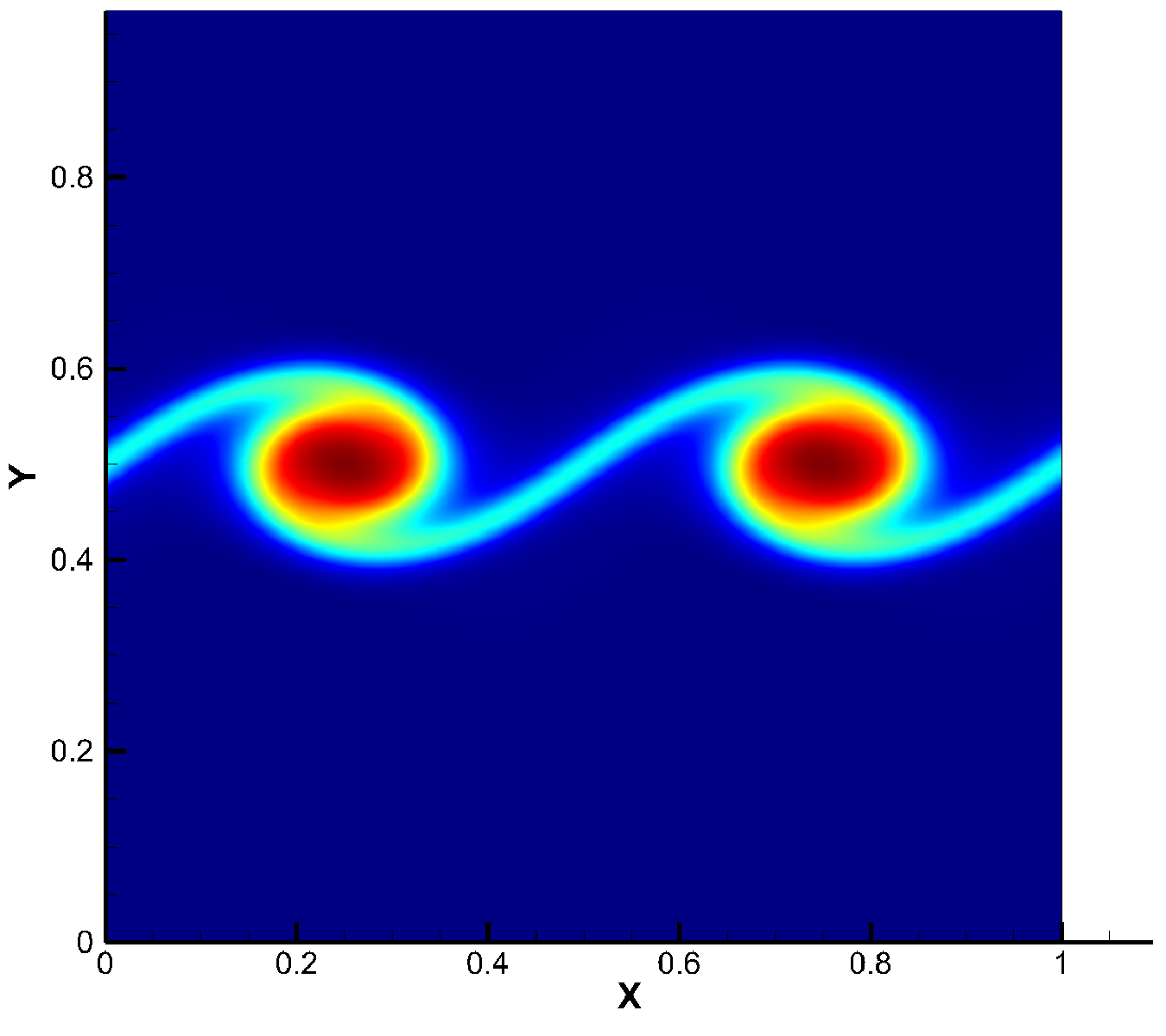}
		\end{minipage}
	}%
\subfigure[$t=0.75$]{
		\begin{minipage}[t]{0.24\linewidth}
			\centering
			\includegraphics[width=\textwidth]{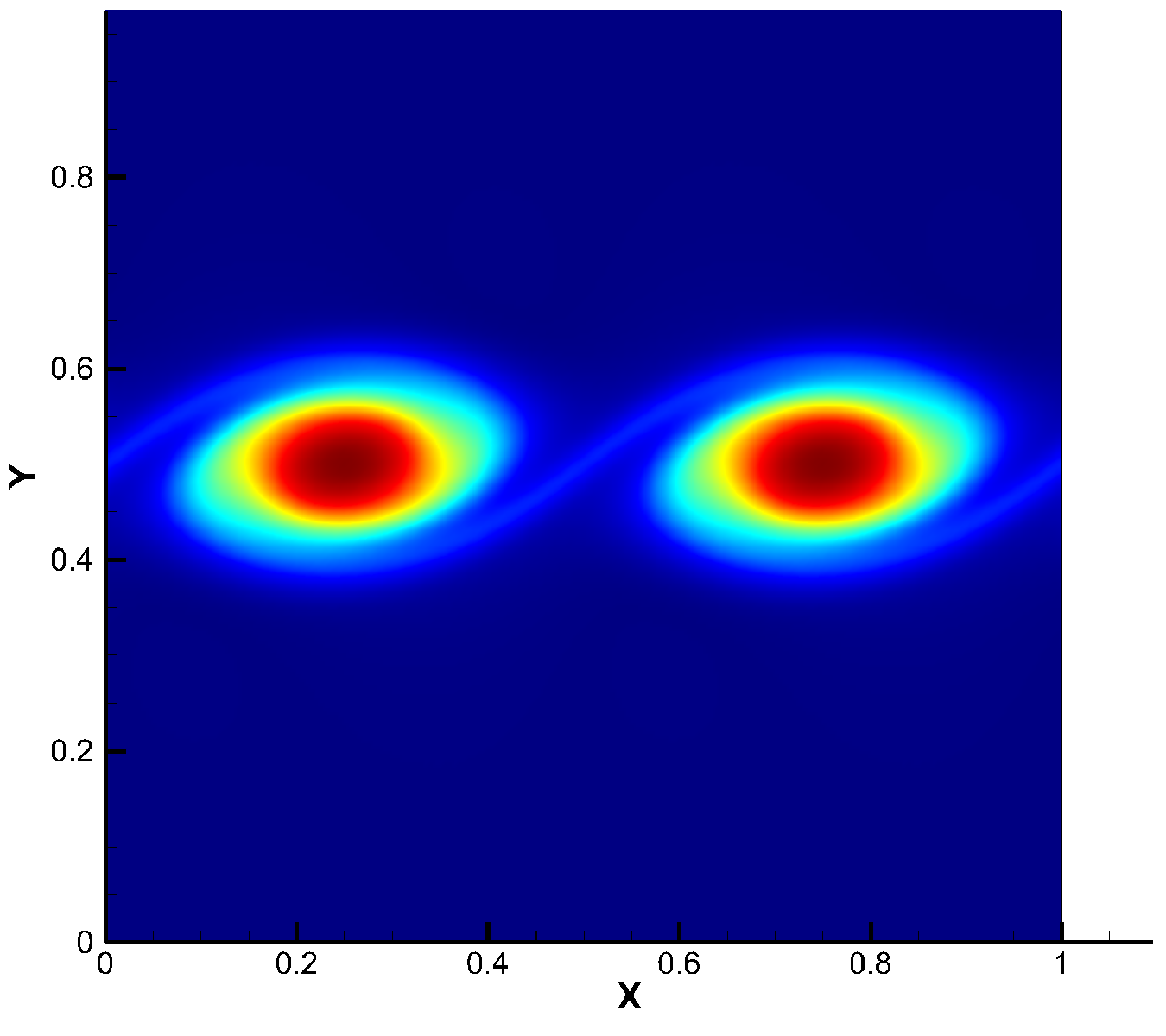}
		\end{minipage}
	}%
	\centering
	\caption{Snapshots of the phase field (upper), vorticity dynamics (lower) perturbed sinusoidal at $t=0.001$ (a), 0.3 (b),  0.5 (c), 0.75 (d) for case I. }
\label{KH-double-mu-1}
\end{figure}

\begin{figure}[htbp]
	\centering
\subfigure[$t=0.001$]{
		\begin{minipage}[t]{0.24\linewidth}
			\centering
			\includegraphics[width=\textwidth]{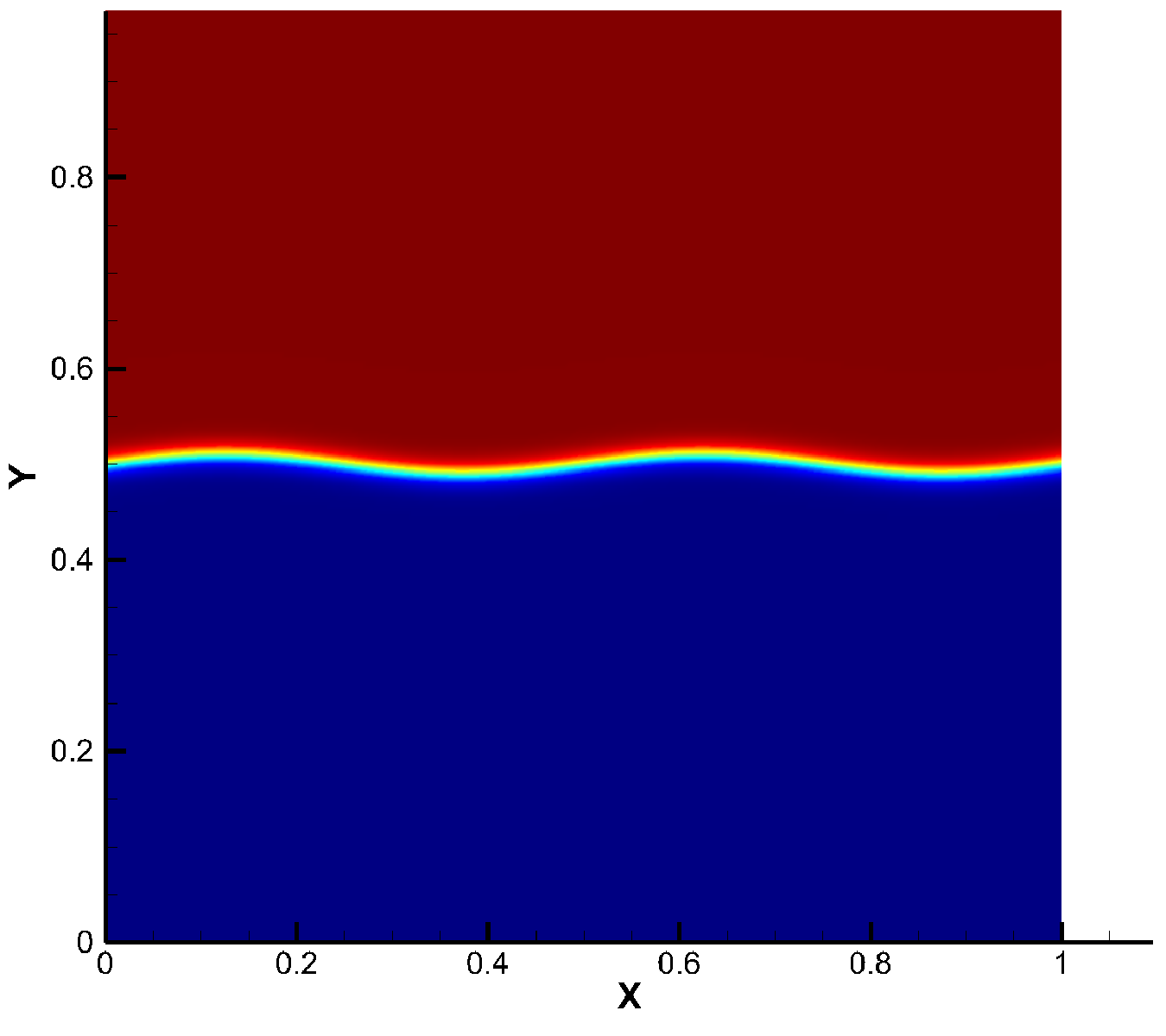}
		\end{minipage}
	}%
	\subfigure[$t=0.3$]{
		\begin{minipage}[t]{0.24\linewidth}
			\centering
			\includegraphics[width=\textwidth]{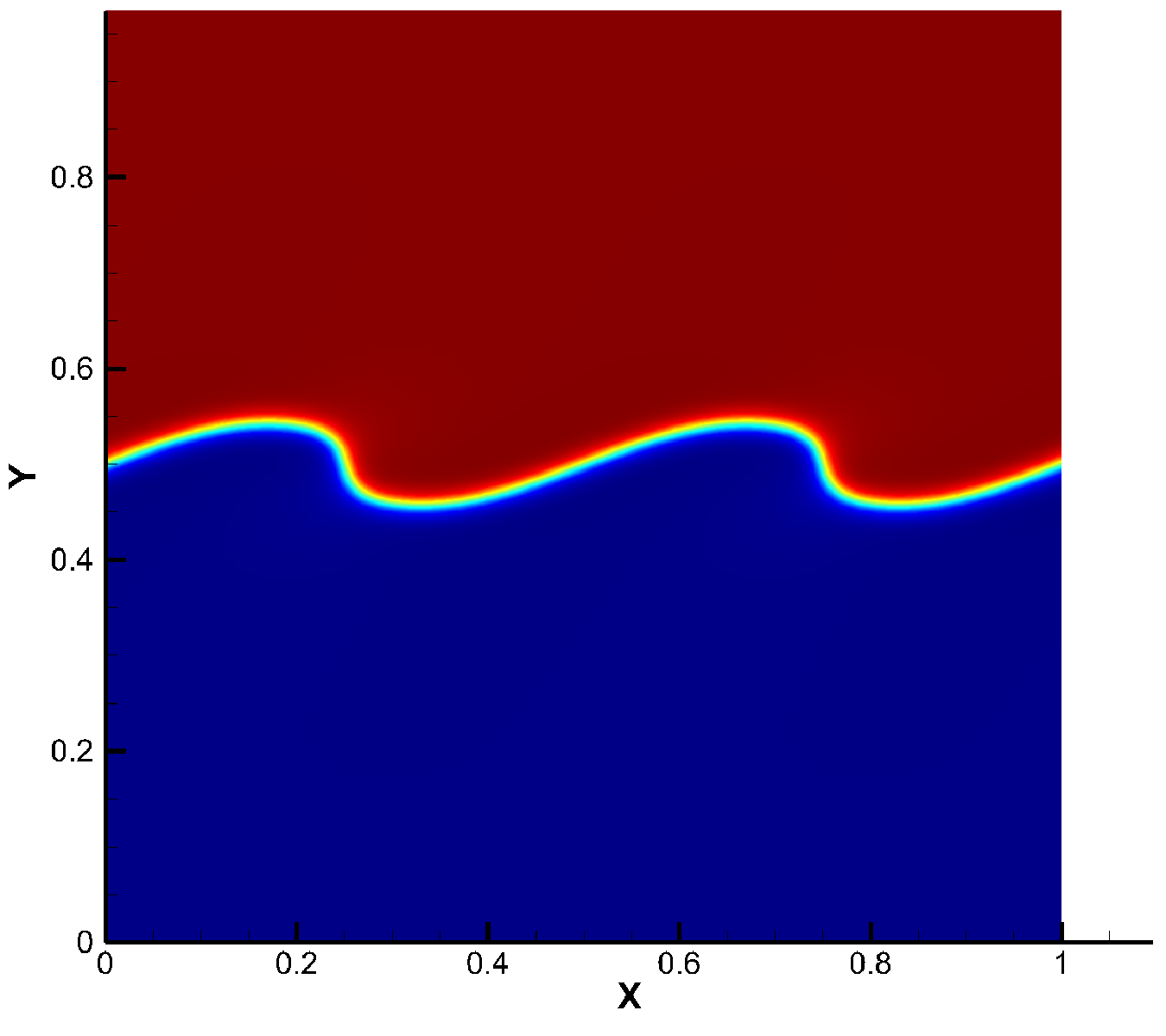}
		\end{minipage}
	}%
   \subfigure[$t=0.5$]{
		\begin{minipage}[t]{0.24\linewidth}
			\centering
			\includegraphics[width=\textwidth]{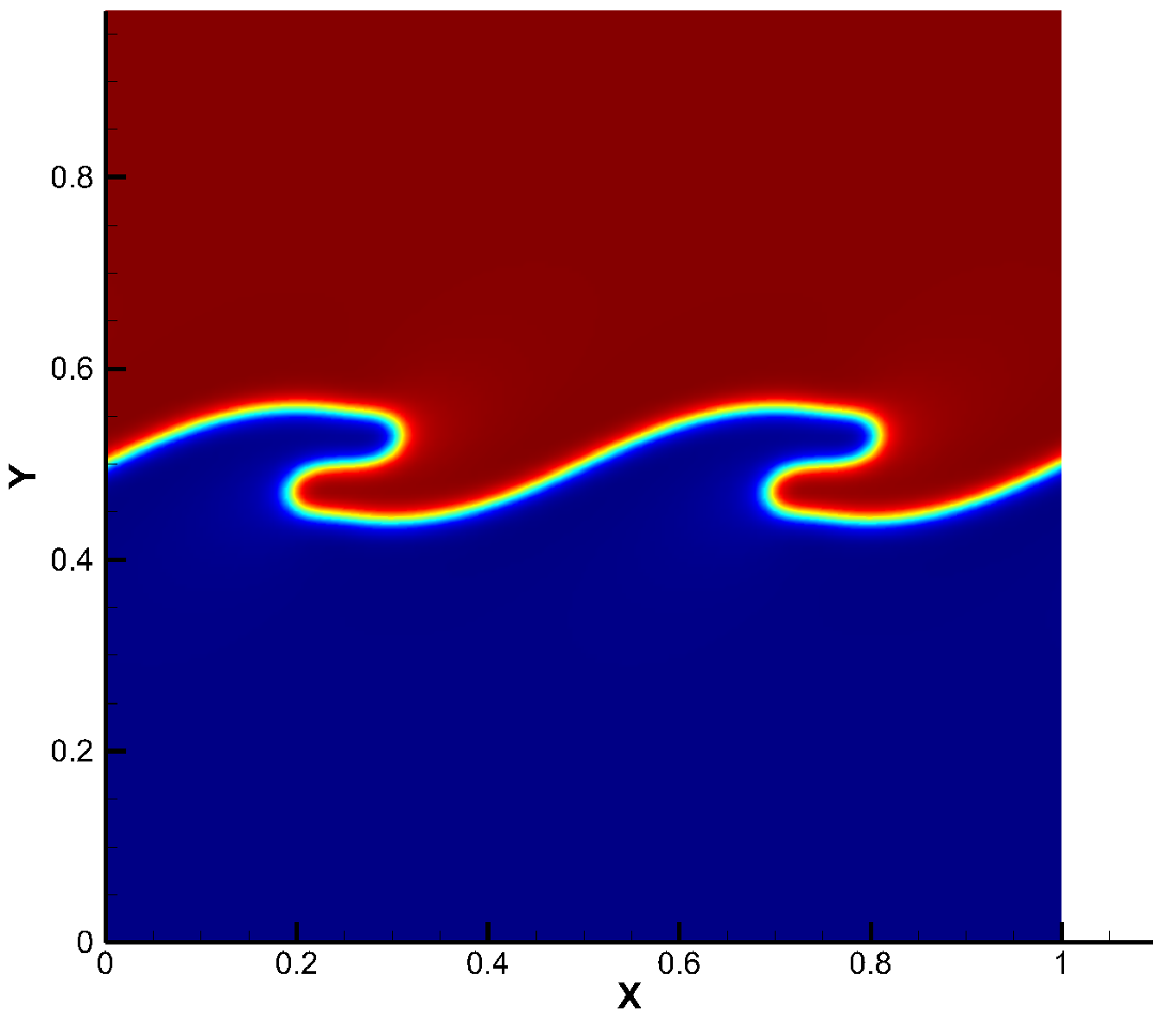}
		\end{minipage}
	}%
\subfigure[$t=0.75$]{
		\begin{minipage}[t]{0.24\linewidth}
			\centering
			\includegraphics[width=\textwidth]{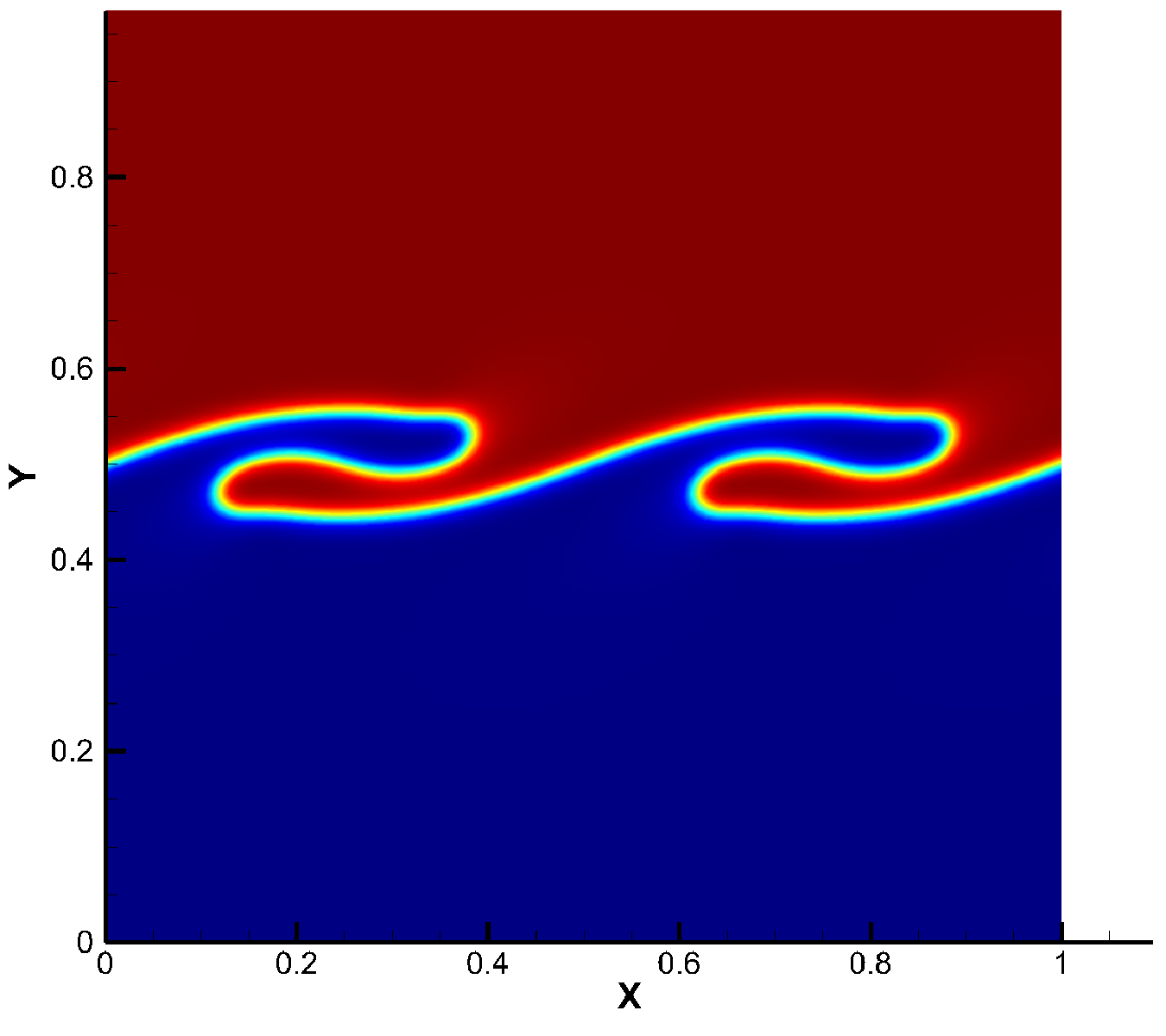}
		\end{minipage}
	}%
	\\ 
\subfigure[$t=0.001$]{
		\begin{minipage}[t]{0.24\linewidth}
			\centering
			\includegraphics[width=\textwidth]{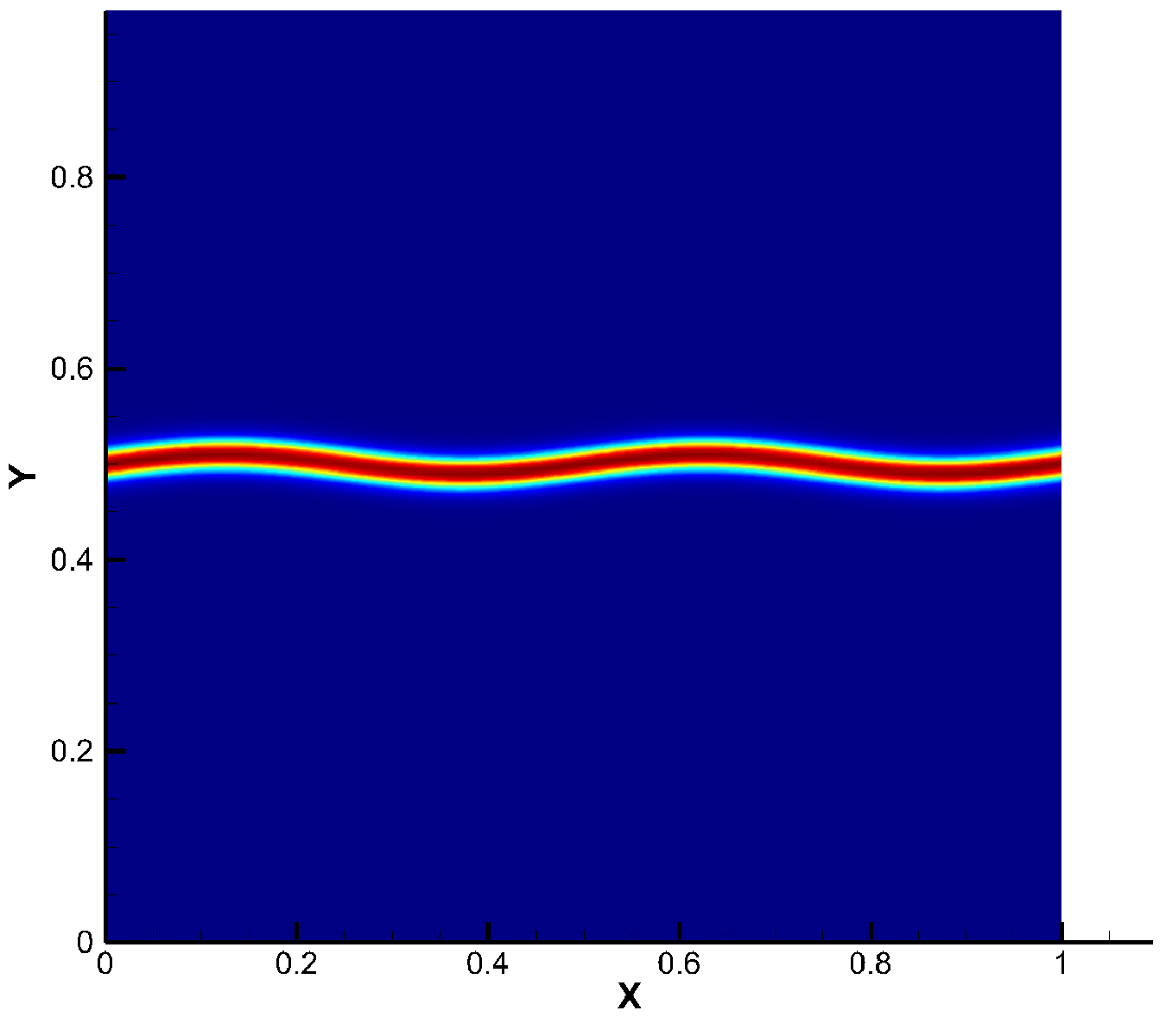}
		\end{minipage}
	}%
	\subfigure[$t=0.3$]{
		\begin{minipage}[t]{0.24\linewidth}
			\centering
			\includegraphics[width=\textwidth]{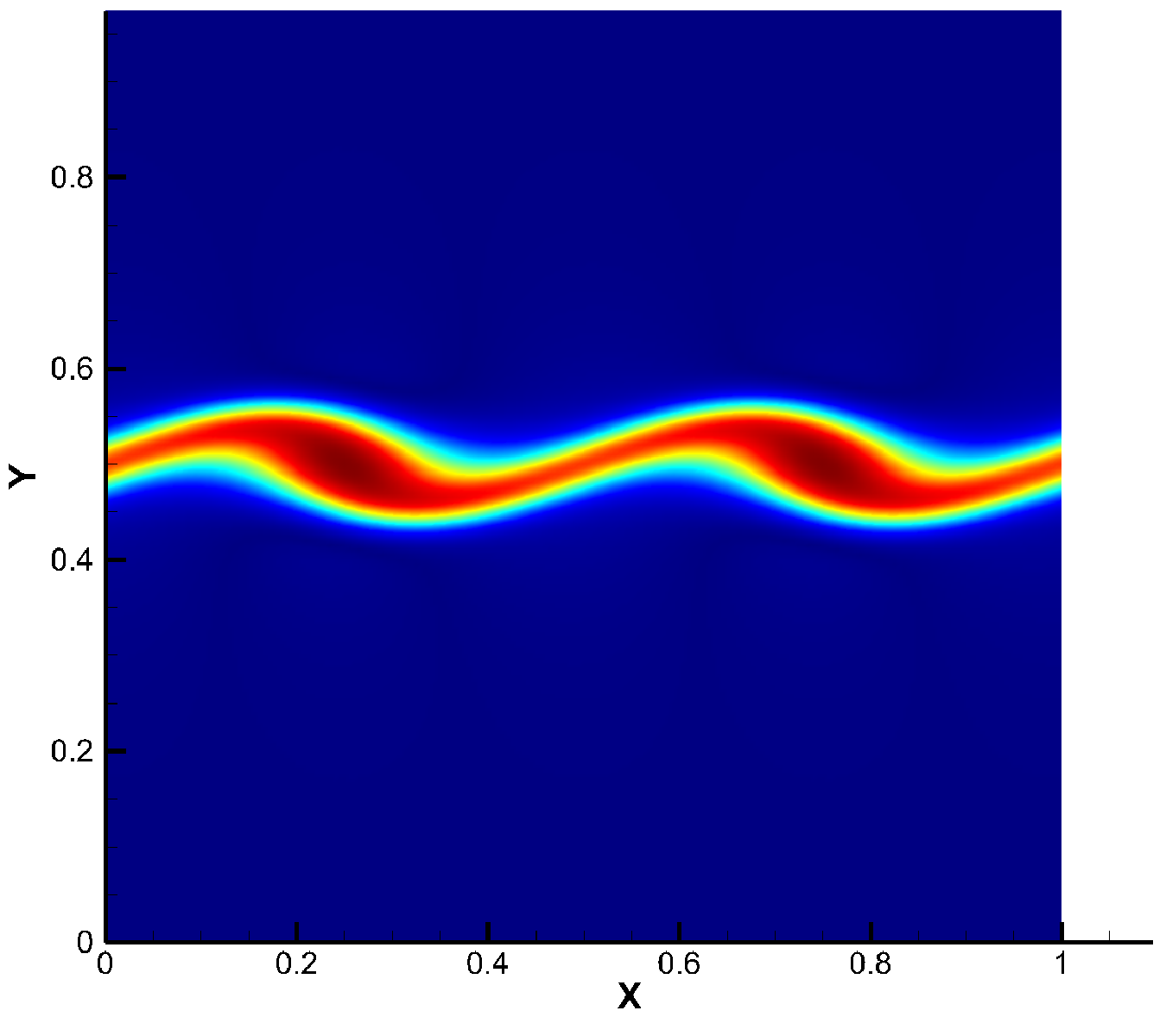}
		\end{minipage}
	}%
   \subfigure[$t=0.5$]{
		\begin{minipage}[t]{0.24\linewidth}
			\centering
			\includegraphics[width=\textwidth]{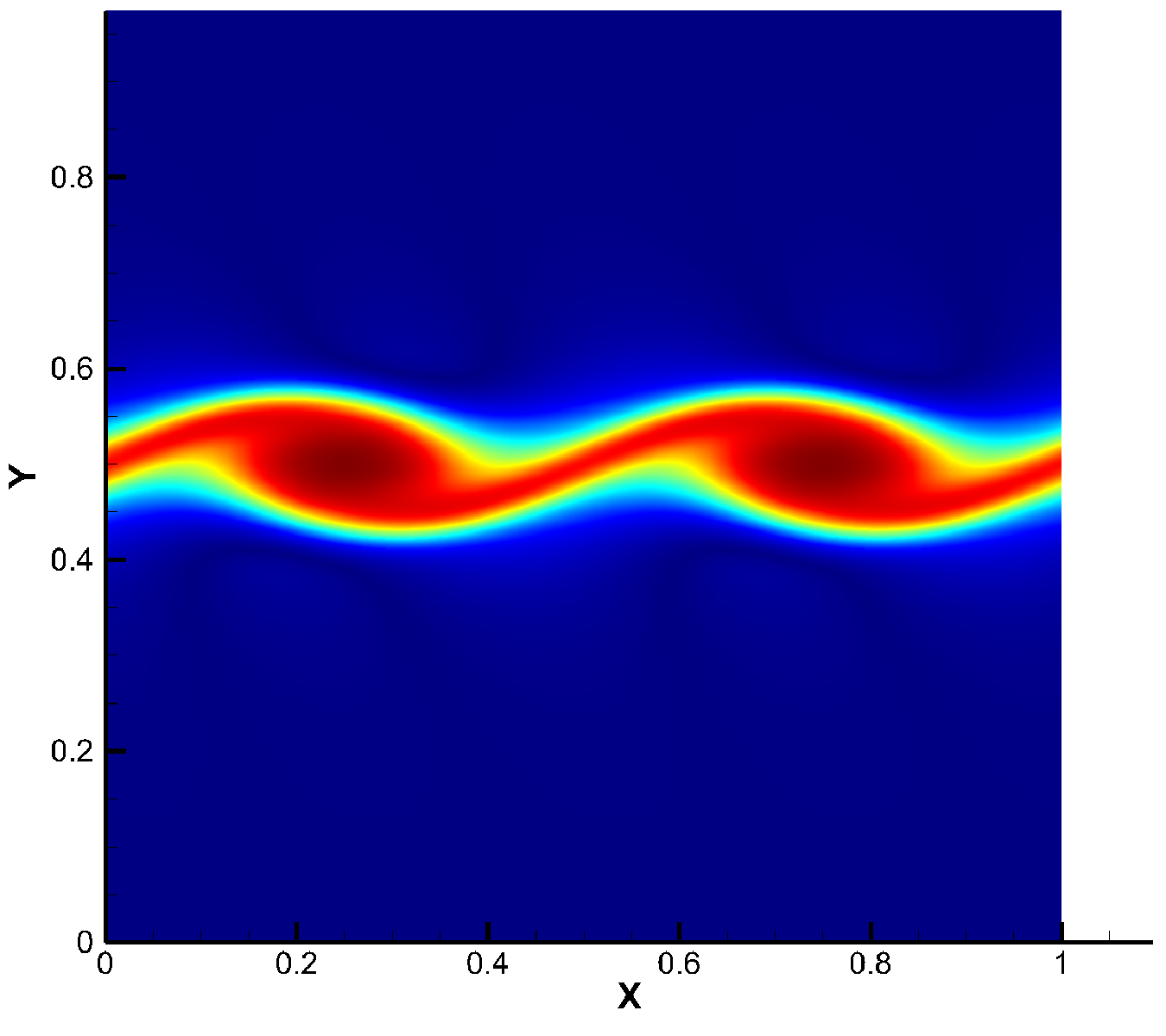}
		\end{minipage}
	}%
\subfigure[$t=0.75$]{
		\begin{minipage}[t]{0.24\linewidth}
			\centering
			\includegraphics[width=\textwidth]{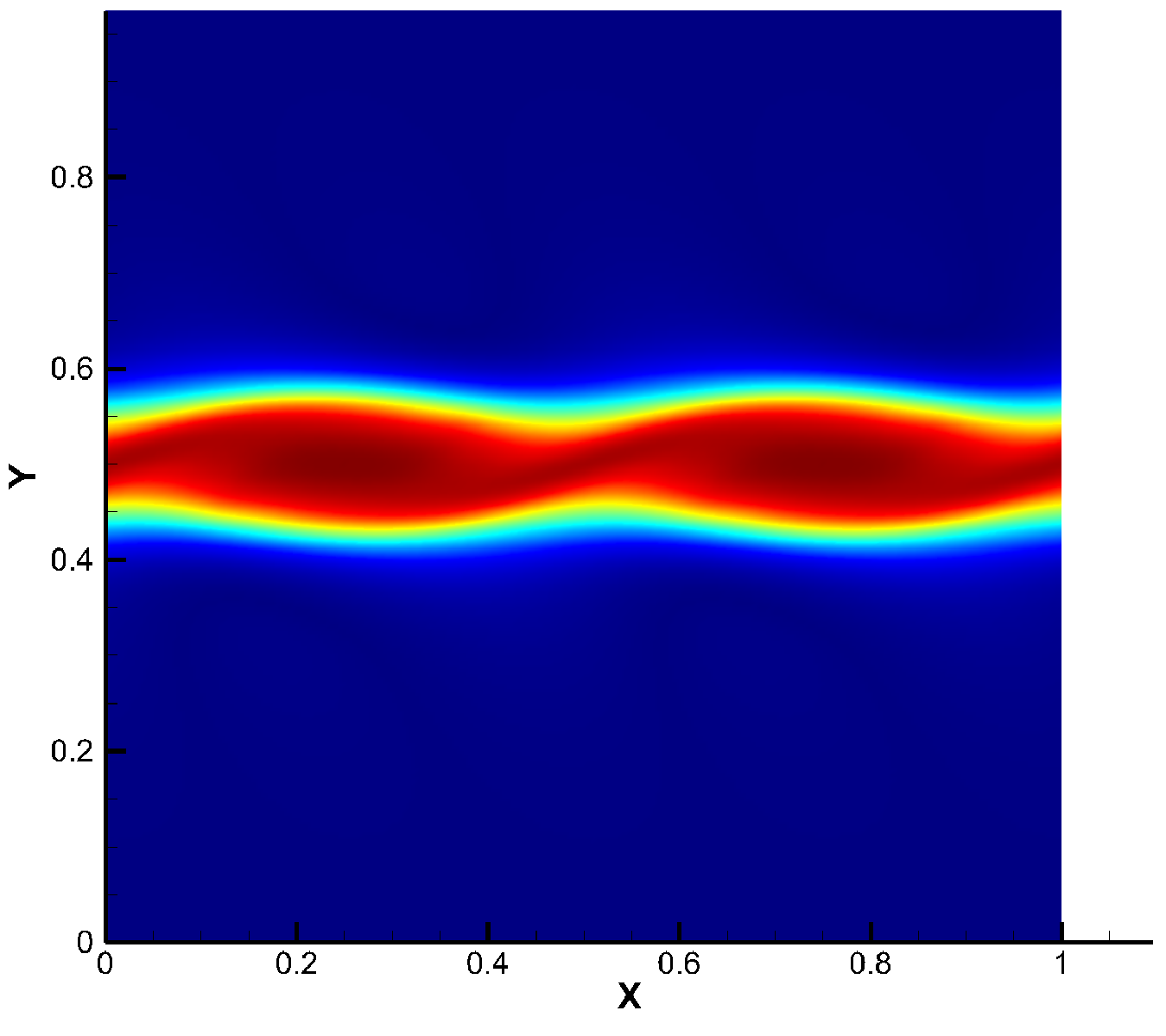}
		\end{minipage}
	}%
	\centering
	\caption{Snapshots of the phase field (upper), vorticity dynamics (lower) perturbed sinusoidal at $t=0.001$ (a), 0.3 (b),  0.5 (c), 0.75 (d) for $\mu=0.1$, $\sigma=10$. }
\label{KH-double-mu-01}
\end{figure}

\begin{figure}[htbp]
	\centering
\subfigure[$t=0.001$]{
		\begin{minipage}[t]{0.24\linewidth}
			\centering
			\includegraphics[width=\textwidth]{double-phase-0001-mu-001.eps}
		\end{minipage}
	}%
	\subfigure[$t=0.3$]{
		\begin{minipage}[t]{0.24\linewidth}
			\centering
			\includegraphics[width=\textwidth]{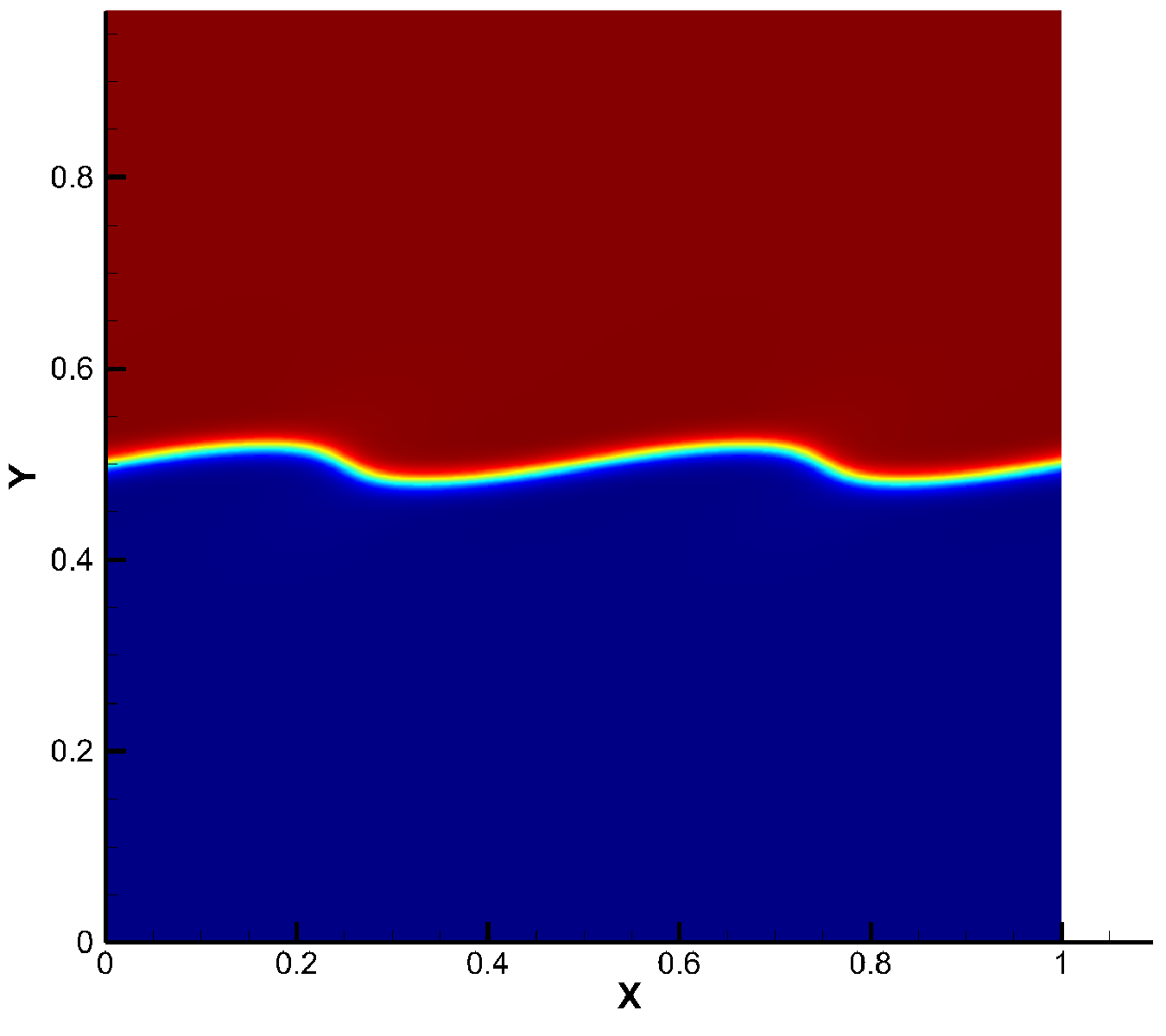}
		\end{minipage}
	}%
   \subfigure[$t=0.5$]{
		\begin{minipage}[t]{0.24\linewidth}
			\centering
			\includegraphics[width=\textwidth]{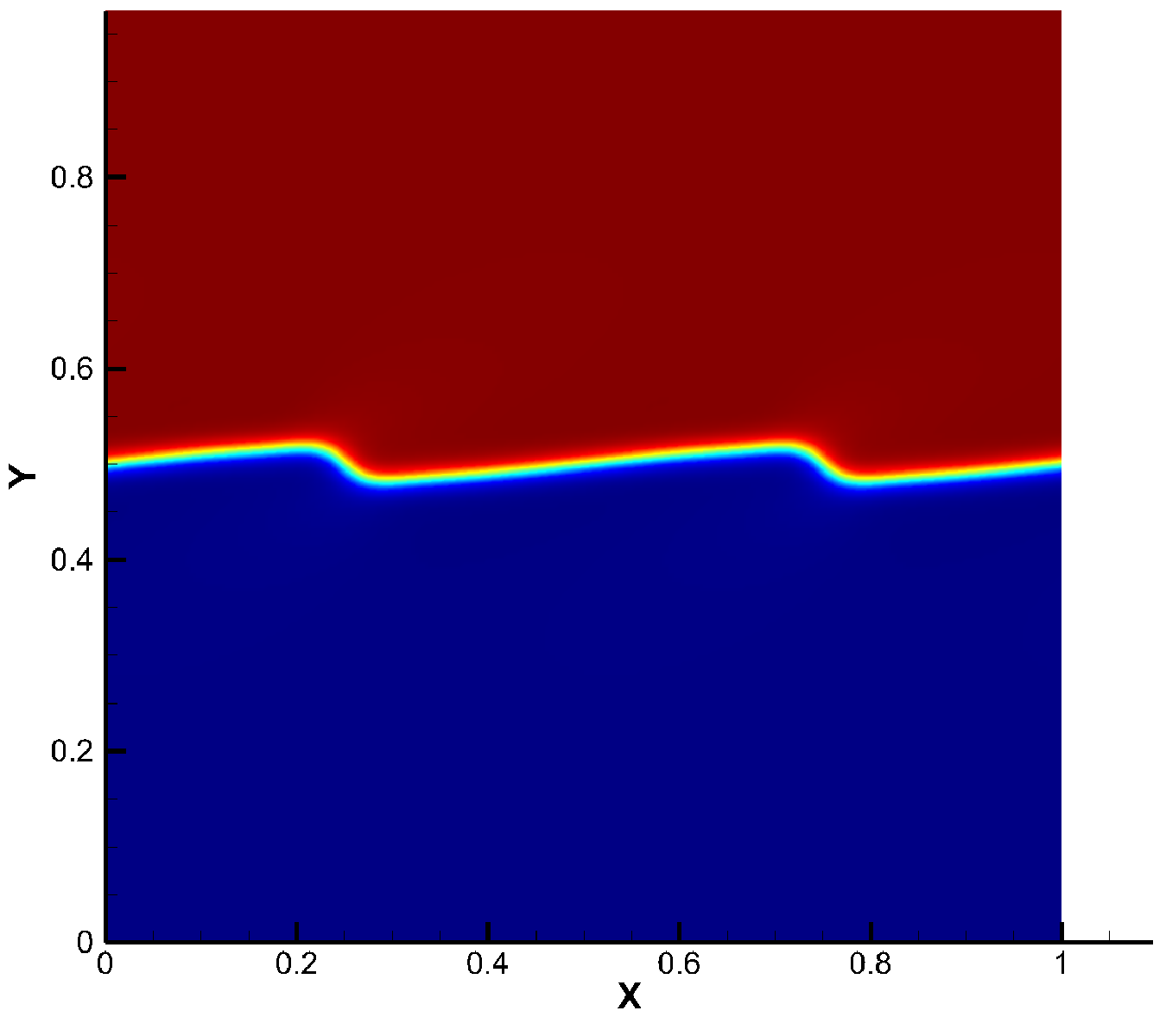}
		\end{minipage}
	}%
	\subfigure[$t=0.75$]{
		\begin{minipage}[t]{0.24\linewidth}
			\centering
			\includegraphics[width=\textwidth]{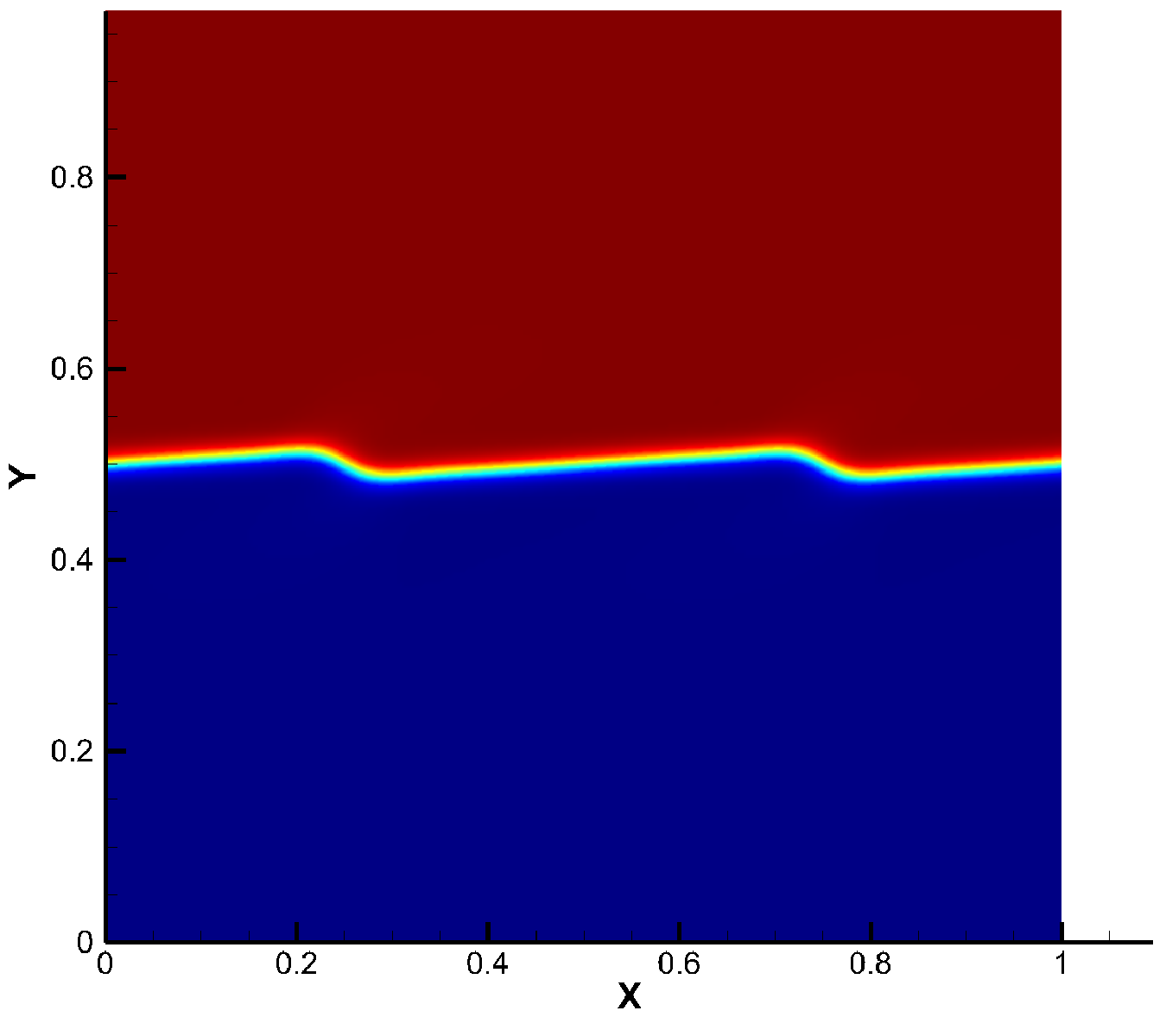}
		\end{minipage}
	}%
	\\
\subfigure[$t=0.001$]{
		\begin{minipage}[t]{0.24\linewidth}
			\centering
			\includegraphics[width=\textwidth]{double-v-0001-mu-001.eps}
		\end{minipage}
	}%
	\subfigure[$t=0.3$]{
		\begin{minipage}[t]{0.24\linewidth}
			\centering
			\includegraphics[width=\textwidth]{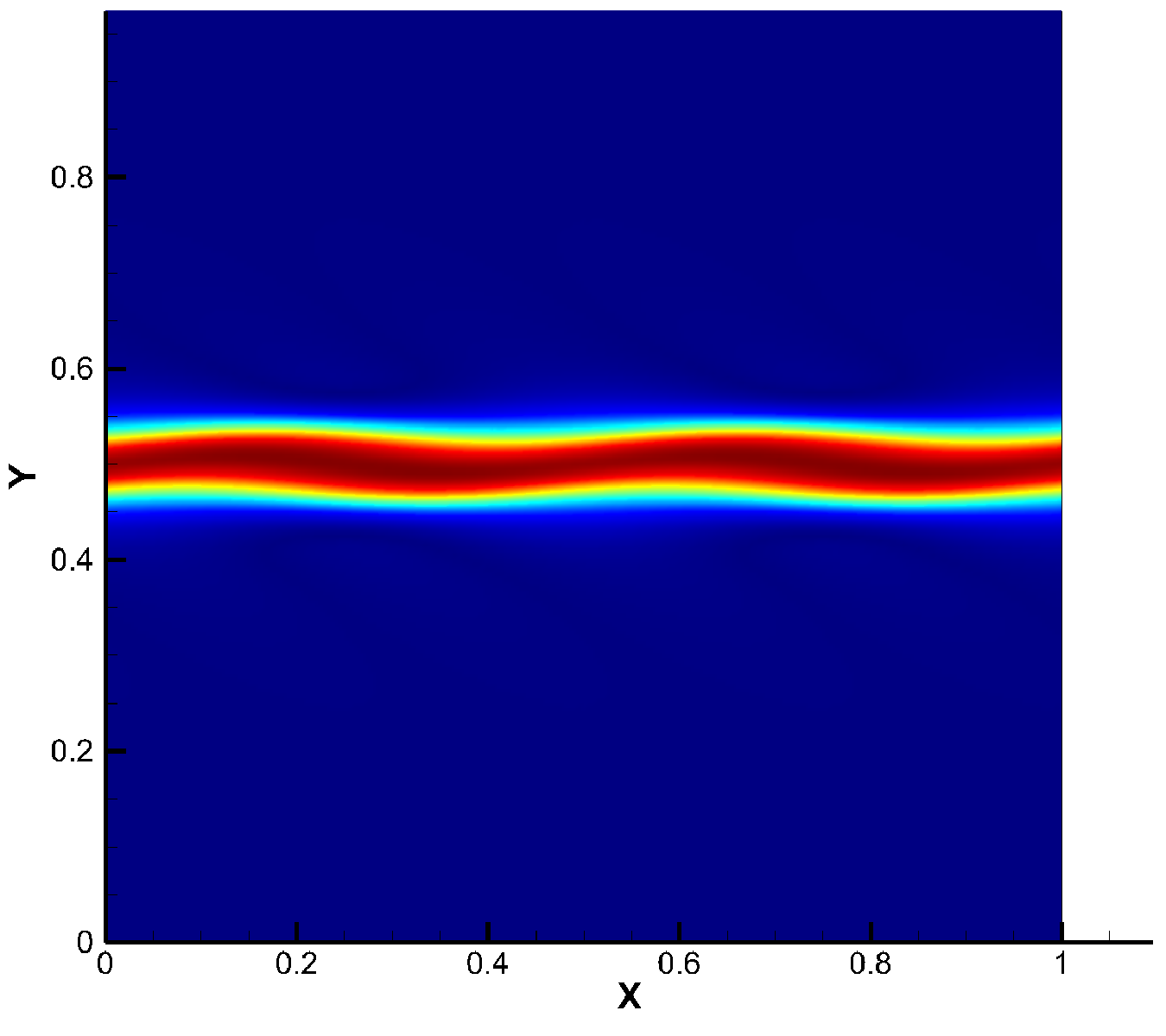}
		\end{minipage}
	}%
   \subfigure[$t=0.5$]{
		\begin{minipage}[t]{0.24\linewidth}
			\centering
			\includegraphics[width=\textwidth]{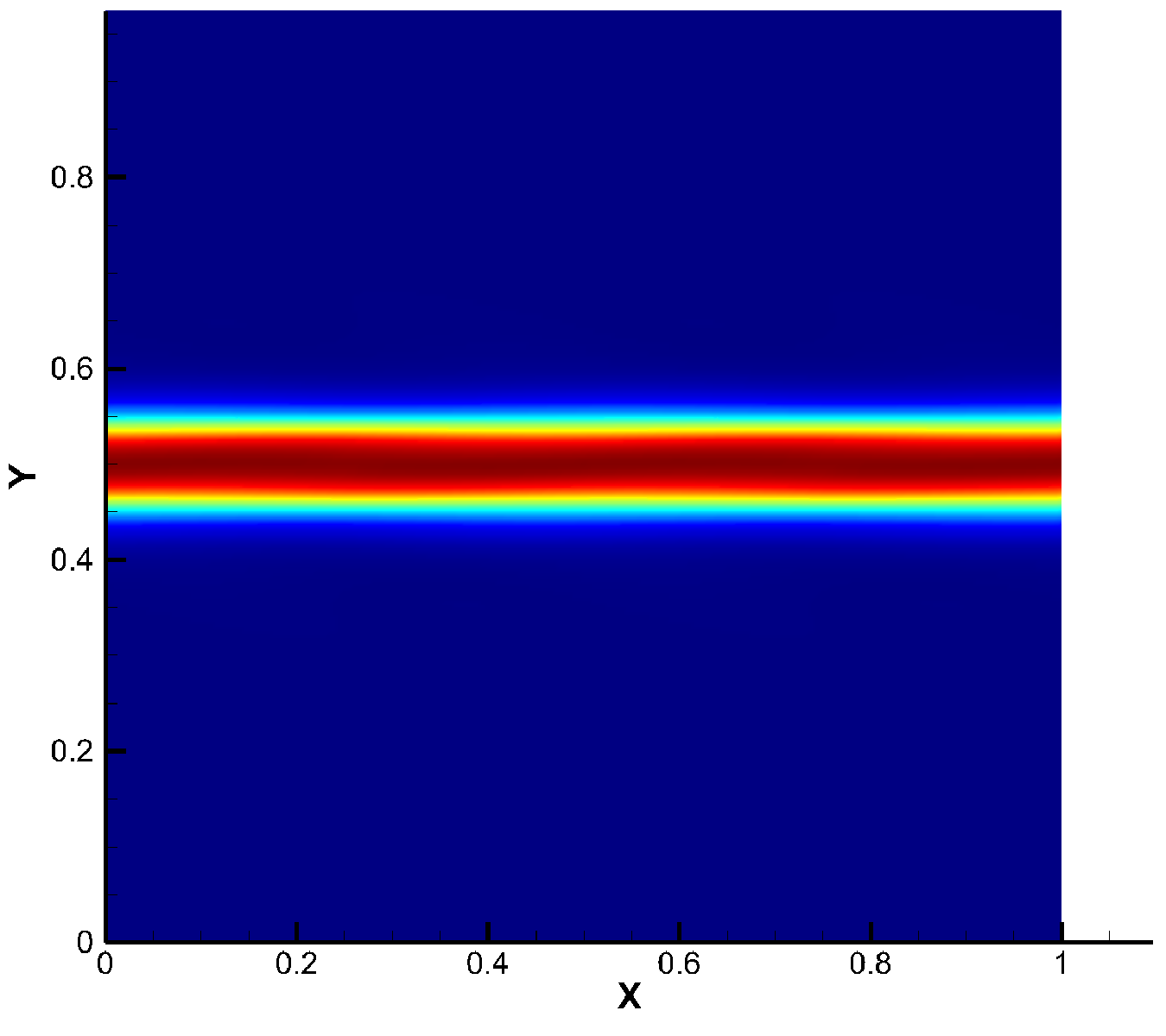}
		\end{minipage}
	}%
\subfigure[$t=0.75$]{
		\begin{minipage}[t]{0.24\linewidth}
			\centering
			\includegraphics[width=\textwidth]{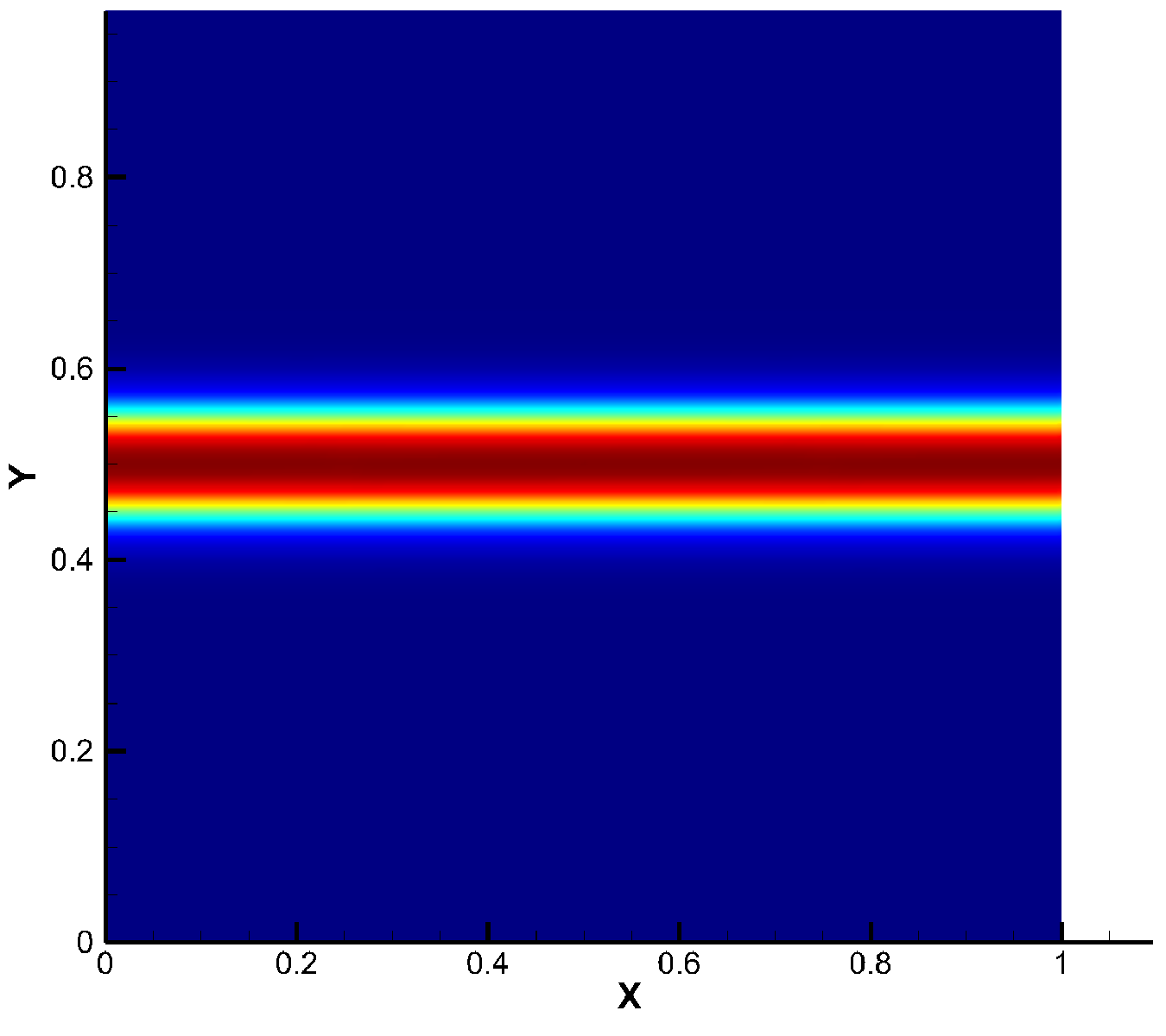}
		\end{minipage}
	}%
	\centering
	\caption{Snapshots of the phase field (upper), vorticity dynamics (lower) perturbed sinusoidal at $t=0.001$ (a), 0.3 (b),  0.5 (c), 0.75 (d) for $\mu=0.01$, $\sigma=100$. }
\label{KH-double-mu-001}
\end{figure}

In this example, we have taken a double mode sinusoidal perturbation at the interface. The domain and parameter
values are same as in the previous problem. The  initial values are given as: 
\begin{eqnarray}\label{K-H2}
\left\{
\begin{aligned} 
\phi_{0}&= \tanh(\frac{y-0.5-0.01 \sin(4\pi x)}{\sqrt{2}\gamma}),\\
\u_{0}&=\Big( \tanh(\frac{y-0.5-0.01 \sin(4\pi x)}{\sqrt{2}\gamma}), 0 \Big),\\
\B_{0}&=(1, 0).
\end{aligned}
\right.
\end{eqnarray} 
The dynamics of the interface and vorticity profiles are given in Figure \ref{KH-double-mu-1} for the parameters used in equation (\ref{K-H-parameter1}). Instead of the single mode sinusoidal perturbation, the two ``cat's eye'' patterns emerge gradually in Figure \ref{KH-double-mu-1}. To investigate the magnetic effect on the mixing fluid mixtures, we adjust the parameters related to the Lorentz force, as defined in equation   (\ref{K-H-parameter1}) to the following values: (a) $\mu=0.1$, $\sigma=10$, and (b) $\mu=0.01$, $\sigma=100$.  The results are shown in Figure \ref{KH-double-mu-01} and Figure \ref{KH-double-mu-001},  respectively.

\section{Conclusion Remarks}\label{sec-conclusion}

In this paper, we develop the optimal $L^2/\L^{2}$-norm error analysis for a convex-splitting FEM for the two-phase diffuse interface MHD model. 
We use the Ritz and Stokes quasi-projections to handle the pollution from the lower-order approximations. The following optimal error estimates 
	\begin{align*}
		&\max\limits_{0\leq k\leq K-1}\|\phi^{k+1}-\phi_{h}^{k+1} \| +\left(\Delta t \sum_{k=0}^{K-1}\|\omega^{k+1}-\omega_{h}^{k+1}\|^{2}\right)^{\frac{1}{2}}\leq C_{0}(\Delta t+h^{r+1}),\\
		&\max\limits_{0\leq k\leq K-1}\|\nabla(\phi^{k+1}-\phi_{h}^{k+1}) \|\leq C_{0}(\Delta t+h^{r}),\\
		&\max\limits_{0\leq k\leq K-1}\left(\|\u^{k+1}-\u_{h}^{k+1}\| +\|\B^{k+1}-\B_{h}^{k+1}\| \right)\leq C_{0}(\Delta t+\beta_{h}),\\
		&\left(\Delta t\sum_{k=0}^{K-1}\left( \|\nabla \cdot (\B^{k+1}-\B_{h}^{k+1})\|^{2} +\|\nabla\times(\B^{k+1}-\B_{h}^{k+1} )\|^{2} \right)\right)^{\frac{1}{2}}\leq C_{0}(\Delta t+\beta_{h}^{\star}).
	\end{align*}
However, the current work focuses only on the matched elements in velocity field and magnetic field. The development of lower-order approximations for the magnetic field will be addressed in our future research.

\bibliographystyle{elsarticle-num}
\bibliography{references}







\end{document}